\documentclass[12pt,BCOR10mm]{scrbook}

\usepackage[T1]{fontenc}
\usepackage[scaled=.90]{helvet}
\usepackage{courier}

\usepackage{amssymb,amsmath,amsthm}
\usepackage{graphicx}
\usepackage[all]{xy}
\usepackage[refpage]{nomencl}
\usepackage{multirow}

\newcommand{\tmop}[1]{\ensuremath{\operatorname{#1}}}
\newcommand{\dueto}[1]{\textup{\textbf{(#1) }}}

\theoremstyle{plain}
\newtheorem{definition}{Definition}[chapter]
\newtheorem{lemma}[definition]{Lemma}
\newtheorem{theorem}[definition]{Theorem}
\newtheorem{corollary}[definition]{Corollary}
\newtheorem{proposition}[definition]{Proposition}

\theoremstyle{remark}

\newtheorem*{note}{Note}
\newtheorem{remark}[definition]{Remark}
\newtheorem{notation}[definition]{Notation}

\begin{document}

\title{Representation Theoretical Construction of the Classical Limit and Spectral Statistics of Generic Hamiltonian Operators}
\author{von Ingolf Sch\"afer}
%\date{\today\\[2cm]\ }
\publishers{\ \\[5cm]\textbf{Dissertation}\\[1cm] 
zur Erlangung des Grades eines Doktors der Naturwissenschaften\\
-- Dr.~rer.~nat.~--}
\lowertitleback{
\begin{tabular}{lll}
\multicolumn{3}{l}{Datum des Promotionskolloquiums: 09.11.2006}\\[0.5cm] 
Gutachter: & Prof.~Dr.~E.~Oeljeklaus & (Universit\"at Bremen)\\
&  Prof.~Dr.~A.~T.~Huckleberry & (Ruhr-Universit\"at Bochum)
\end{tabular}
}
\maketitle

\tableofcontents
\listoffigures

% \chapter{Introduction}
\chapter{Introduction}

The theory of spectral statistics is concerned with the spectral
properties of ensembles of linear operators. Typically, these depend
on a parameter $N$ which is supposed to be very large or even
approaching infinity. The origin of this field is quantum physics,
where such ensembles arose as models for the energy spectra of large
atoms.

Another branch of physics, namely semiclassical physics, is also
concerned with such ensembles and their spectral statistics. In
semiclassical physics large values of $N$ should correspond to a
quantum mechanical system which approaches classical
mechanics. Details about these relations can be found in \cite{mehta}
and \cite{haake}.

Finally, spectral statistics have been studied in the context of number
theory, with the most famous example being the distribution of zeros
of the Riemann $\zeta$-function on the critical line. An
introduction to this field is given in \cite{snaith}.

Under the assumption of genericity one might hope that there exist
natural sequences of operators taken from these ensembles such that
the spectral properties of the individual operators reflect those of
the ensembles.

We are concerned here with two examples, in which spectral statistics
appear. The first being the theory of Random Matrices. In this theory
natural sequences of symmetric spaces with invariant measures on them
are given. These spaces have natural representations as matrices and
one is interested in the limit of the spectral statistics as $N\to
\infty$. An example is the sequence of unitary groups $\tmop{U}(N)$
with the Haar measure. In \cite{katzsarnak} it is proven that a limit
measure of a special kind of spectral statistics exists for this
example.

The second example, in which spectral statistics appear, is given by
the approach suggested in \cite{haakekus}. In this article the authors
consider two fixed operators in the universal enveloping algebra of
$SL(3,\mathbb{C})$ in a sequence of irreducible representations of
$SL(3,\mathbb{C})$ and study the spectral statistics by numerical
methods. The motivation from the approach stems from a previous paper
(cf.\ \cite{gnutzmannkus}) of two of the authors: Such a sequence of
irreducible representations occurs in the construction of the
classical mechanical system in the limit of a quantum mechanical
system with $SL(3,\mathbb{C})$ symmetry. We will follow this approach
in the following chapters.

Our main device in the study of spectral statistics is the nearest
neighbor statistics, i.e.\ the normalized distribution of distances of
neighboring eigenvalues (counted with multiplicity) of such linear
operators. It is frequently drawn as a histogram (see Figure
\ref{fig:histogram}). A detailed explanation of this plot can be found
in the Appendix. 
\begin{figure}
  \centering
  \includegraphics{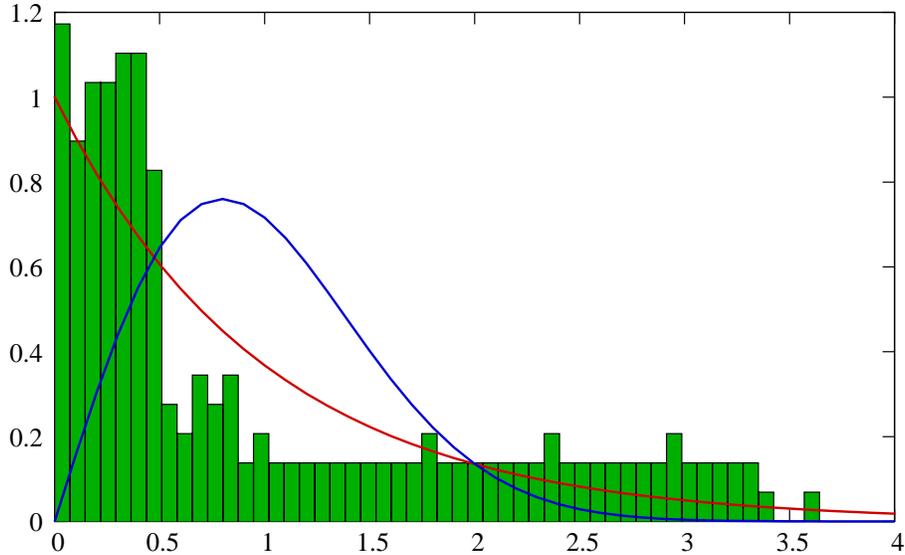}
  \caption{A sample histogram of the nearest neighbor statistics}
  \label{fig:histogram}
\end{figure}

The nearest neighbor statistics lead to Borel measures on the positive
real line by putting a Dirac measure for every occurring distance of
neighboring eigenvalues with proper normalization. Out of the
wealth of notions of convergence for such measures we choose the
weak convergence (in probability theory: convergence in distribution)
and the Kolmogorov-Smirnov convergence. The Kol\-mo\-go\-rov-Smir\-nov
distance of two measures $\mu,\nu$ is given by
\begin{equation}
  d_{KS}(\mu,\nu)=\sup_{t\in\mathbb{R}}\left|\int_{-\infty}^td\mu -
\int_{-\infty}^td\nu\right|,
\end{equation}
i.e., Kolmogorov-Smirnov convergence is uniform convergence of the
cumulative distribution functions. We will examine $d_{KS}$ for
sequences of individual operators relative to a fixed measure $\nu$,
but also average $d_{KS}$ with respect to a fixed probability measure
$\nu$ over the full ensemble. Here sequences of irreducible
representations will arise.

This text is structured into six chapters. Following the approach
in \cite{gnutzmannkus} we give a general construction of the classical
limit for semi-simple compact Lie groups in Chapter
\ref{chap:rep-theo-limit}. This can be done in a functorial way, but
the objective of Chapter \ref{chap:rep-theo-limit} is to give an
interpretation as a mathematical limit as a parameter $n$ converges to
$\infty$.

Chapter \ref{chap:spec-stat-simple} deals with the spectral statistics
of operators in the Lie algebra along sequences of irreducible
representations. It is necessary to discuss possible scalings of
these operators in this context.

The goal of Chapter \ref{chap:spec-stat-gen} is to study the spectral
statistics of exponentiated operators, which satisfy certain
conditions of genericity, in a certain completion of the universal
enveloping algebra of a semi-simple complex Lie group. The main tools
are Birkhoff's Ergodic Theorem and an estimation on $d_{KS}$ for
maximal tori of $U(N)$.

Chapter \ref{chap:poisson-spectral-stat} is devoted to the proof of 
this estimation, where we follow the structure of
\cite{katzsarnak} for the proof.

In the Appendix we collect the necessary background facts of
representation theory and symplectic geometry for the readers'
convenience. The Appendix closes with some general observations about
nearest neighbor statistics.

\section*{Acknowledgments}
I would like to mention all the people who helped me while writing
this thesis. First and foremost, I would like to thank my supervisors
Alan T.~Huckleberry and Eberhard Oeljeklaus for their support and
guidance. Furthermore, I am indebted to Marek Ku\'s, who gave me
insight on the physical motivation for the topic of this thesis and
who invited me to Warsaw for joint research.

Special thanks go to my colleagues in Bremen and Bochum, where I would
like to mention Kristina Frantzen, Daniel Greb, Christian Miebach,
Elmar Plischke, Patrick Sch\"utzdeller, Monika Winklmeier and
Christian Wyss, who always found time to discuss and helped me in
various ways.

Last but not least, I would like to thank my wife Silke Sch\"afer for
her support and patience.

Support by the Sonderforschungbereich TR 12 ``Symmetries and
Universality in Mesoscopic Systems'' is gratefully acknowledged.

%%% Local Variables: 
%%% mode: latex
%%% TeX-master: "../thesis"
%%% End: 

% moved chapter 1 to appendix

% \chapter{Construction}
\chapter[Construction of the Classical Limit]{Representation
  Theoretical Construction of the Classical Limit}
\label{chap:rep-theo-limit}
In this chapter we will give a construction of the classical limit of
Hamiltonian mechanics by a representation theoretical approach. Our method
is an abstract generalization of the method given in
\cite{gnutzmannkus} and \cite{gnutzmann} and covers systems with
compact semi-simple Lie groups as symmetry groups.

The following notation will be used without further
notice (for details cf.~Appendix): $K$ is a compact semi-simple
Lie group with complexification $G$ and the corresponding Lie
algebras are $\mathfrak{k}$ and $\mathfrak{g}$. Every representation of
$K$ will be assumed to be continuous, finite-dimensional and unitary,
where the scalar product is denoted by $\langle\cdot,\cdot\rangle$. By
convention the scalar product is complex linear in the second, and
anti-linear in the first variable.

Furthermore we assume that we have fixed a Borel subgroup $B\subset G$
and obtain a notion of positivity of roots and weights. Recall that
the choice of $B$ also determines a maximal tours $T\subset K$.

\section{The Classical Limit in the Simple Case}
\label{sec:class-limit-simple}

A guiding principle in quantum mechanics is that of correspondence.
It states that quantum mechanical systems whose size is large
compared to microscopical length scales can be described by classical
physics. The classical system attached to the quantum mechanical
system is called the classical limit (cf.~\cite{gnutzmannkus}). So
there should be some kind of functor from Hilbert spaces with
Hamiltonian operators to symplectic manifolds with Hamiltonian
functions. Actually, one might require that this functor is inverse to
so-called geometric quantization.  At least it should satisfy the
Dirac correspondence, i.e., if $\xi_{H_1}$ and $\xi_{H_2}$ are two
Hamiltonian operators with corresponding Hamiltonian functions $h_1$
and $h_2$, then the Lie bracket of $\xi_{H_1}$ and $\xi_{H_2}$ should
correspond to the Poisson bracket of $h_1$ and $h_2$:
\begin{equation}
  \label{eq:diraccondition}
  [\xi_{H_1},\xi_{H_2}] \mapsto c \{h_1,h_2\},
\end{equation}
where $c$ is a constant, usually $i\hbar$.

More often, one discusses the opposite direction, i.e., quantization
(cf.~\cite{woodhouse} Chapter~9.2).  Therefore one may call the
procedure presented here \emph{dequantization}.

Let $\rho:K\to \tmop{U}(V)$ be an irreducible representation.  Let
$\rho_*:\mathfrak{k}\to \tmop{End}(V)$ be the induced representation
of the Lie algebra.  Both $\rho$ and $\rho_*$ extend to holomorphic
resp.\ linear representations of the corresponding complexifications
$G$ and $\mathfrak{g}$. To keep notation as simple as possible we will
also denote these by $\rho$ and $\rho_*$

The map $\mu:\mathbb{P}(V)\to\mathfrak{k}^*$ given by
\begin{equation}
  \label{eq:momentummap}
  \mu^\xi([v])= -2i\frac{\langle v,\rho_*(\xi).v\rangle}{\langle v,
    v\rangle}\forall\xi\in\mathfrak{k}, v\in\mathbb{P}(V)
\end{equation}
is the momentum map with respect to the symplectic structure on
$\mathbb{P}(V)$ induced by the Fubini-Study metric (cf.\ Appendix for
details). Moreover, if $\lambda\in\mathfrak{t}^*$ is the highest
weight of $\rho$, then
\begin{equation}
  \mu([v_{max}])=\lambda
\end{equation}
for any vector $v_{max}$ of highest weight.

Since $\mu$ is an $K$-equivariant map and the stabilizers of $\lambda$ and
$v_{max}$ agree, this map is a symplectic diffeomorphism of the orbit
$K.[v_{max}]$ onto the coadjoint orbit $K.\lambda$ with the
Kostant-Kirillov form.

In the literature, this coadjoint orbit is called the \textbf{set of
  coherent states} (cf.\ \cite{perelomov}, \cite{woodhouse}). To
simplify notation we write $Z=K.\lambda$ for this set.

Equivariance implies that the map $\tilde{\mu}:\mathfrak{k}\to
C^{\infty}(Z),\ \xi \mapsto \mu^\xi (\cdot)$, satisfies
\begin{equation}
  \label{eq:equivariance_of_momentummap}
  \tilde{\mu}([\xi_1,\xi_2]) = \{\tilde{\mu}(\xi_1),\tilde{\mu}(\xi_2)\}.
\end{equation}
If we compare this equation with the Dirac condition
\eqref{eq:diraccondition}, then, up to constants, this is exactly what
we are looking for. But the Lie algebra $\mathfrak{k}$ acts by
\emph{skew} self-adjoint operators on $V$. Thus we define
$\tmop{cl}:i\mathfrak{k}\to C^\infty(Z)$ for an element $\xi_H\in
i\mathfrak{k}$ by
\begin{equation}
  \label{eq:definition_simple_classical_limit}
  \tmop{cl}(\xi_H)([x]) = \frac{1}{2}\tilde{\mu}(i\xi_H)(x) = \frac{\langle
    x,\rho_\ast(\xi_H).x
    \rangle}{\langle x,x \rangle}, 
\end{equation}
where the factor $\frac{1}{2}$ will become clear in the following. 
First note that while $i\xi_H$ is represented as a skew self-adjoint
operator, $\xi_H$ is self-adjoint. Now, we have the following version of
the Dirac correspondence for the classical limit $\tmop{cl}$:
\begin{equation}
  \label{eq:dirac_for_classical limit}
  \begin{aligned} 
    \tmop{cl}(i[\xi_{H_1},\xi_{H_2}]) & = \frac{1}{2}\tilde{\mu}(ii[\xi_{H_1},\xi_{H_2}])
    = \frac{1}{2}\tilde{\mu}([i\xi_{H_1},i\xi_{H_2}]) \\
    & =
    2\cdot\{\frac{1}{2}\tilde{\mu}(i\xi_{H_1}),\frac{1}{2}\tilde{\mu}(i\xi_{H_2})
    \} = 2\cdot\{\tmop{cl}(\xi_{H_1}),\tmop{cl}(\xi_{H_2})\}.
  \end{aligned}
\end{equation}

\section{The Classical Limit in the General Case}
\label{sec:class-limit-gener}
So far our classical limit has been defined for those self-adjoint
operators which can be expressed as the image of an element of
$i\mathfrak{k}$ under $\rho_*$. But we want to define the classical
limit for every self-adjoint linear operator on $V$. In fact, it will
be defined for all linear operators on $V$, although in general we
do not obtain real-valued functions on $Z$ if we take the classical
limit of an operator which is not self-adjoint.

Let $\mathcal{T}(\mathfrak{g})$ denote the full tensor algebra of
$\mathfrak{g}$. The Lie algebra representation $\rho_\ast$ extends
uniquely to a representation
$\rho_\ast:\mathcal{T}(\mathfrak{g})\to\tmop{End}(V)$.  This map is
surjective by the lemma of Burnside. Thus, in particular every
self-adjoint operator is contained in the image of $\rho_\ast$.  

We fix an $\mathbb{R}$-basis $\xi_1,\dots,\xi_k$ of $i\mathfrak{k}$
for the rest of this chapter. Note that this is a $\mathbb{C}$-basis
of $\mathfrak{g}$. Thus, an element $\xi_H$ of
$\mathcal{T}(\mathfrak{g})$ has a unique decomposition into
homogeneous terms consisting of sums of ``monomials'' $\xi_{\alpha_1}
\otimes \dots \otimes \xi_{\alpha_p}$ for some indices $\alpha_j \in
\{1,\dots,k\}$.  (These are not monomials in the usual sense because
of the non-commutativity.)

\begin{definition}
  The classical limit of such a ``monomial'' is
  \begin{equation}
    \label{eq:def-cl-general}
    \tmop{cl}(\xi_{\alpha_1} \otimes \dots \otimes \xi_{\alpha_p}) :=
    \tmop{cl}(\xi_{\alpha_1}) \cdot \dots \cdot \tmop{cl}(\xi_{\alpha_p}). 
  \end{equation}
  The classical limit of
  \begin{equation}
    \xi_H=\sum \alpha_I\xi_{\alpha_1} \otimes
    \dots \otimes \xi_{\alpha_{p_I}}\quad \in\mathcal{T}(\mathfrak{g})   
  \end{equation}
  is the sum of all classical limits of each ``monomial'' multiplied
  by the corresponding coefficient.
 
  We call the resulting map $\tmop{cl}:\mathcal{T}(\mathfrak{g})\to
  C^\infty(Z,\mathbb{C})$ the \textbf{classical limit map}.
\end{definition}%
\nomenclature{$\tmop{cl}$}{The classical limit map}%
Let us discuss this definition. First note that if
$\xi_H$ is abstractly self-adjoint, then $\tmop{cl}(\xi_H)$ is
real-valued. To see this, we calculate
\begin{equation}
  \tmop{cl}(\xi_H) = \tmop{cl}(\xi_H^\dagger)=\overline{\tmop{cl}(\xi_H)},
\end{equation}
where the last step is due to (\ref{eq:compatibiliyofdagger}) and
(\ref{eq:def-cl-general}). The converse is false since, in general,
$\mathcal{T}(\mathfrak{g})$ contains nilpotent elements. 

\begin{remark}
  The map $\tmop{cl}:\mathcal{T}(\mathfrak{g})\to C^\infty(Z,\mathbb{C})$ has
  a natural factorization $\tmop{cl}_S:\mathcal{S}(\mathfrak{g})\to
  C^\infty(Z,\mathbb{C})$ to the full algebra of symmetric tensors
  $\mathcal{S}(\mathfrak{g})$.
\end{remark}

In this way the classical limit map is a link between the
non-commutative algebra $\mathcal{T}(\mathfrak{g})$ and a certain
commutative subalgebra of $C^\infty(Z,\mathbb{C})$. But since
$C^\infty(Z,\mathbb{C})$ is commutative, we have to work with the tensor
algebra and cannot pass to the universal enveloping algebra
$\mathcal{U}(\mathfrak{g})$ in the definition of the classical limit,
otherwise the quotient will not be well-defined. To see this, take any
operators $\xi_a$nd $\xi_b$ such that $[\xi_a,\xi_b]\neq 0$. Then it
follows that $\tmop{cl}(\xi_a\xi_b-\xi_b\xi_a-[\xi_a,\xi_b])$ is not
equal to zero.

Let $x\in V$ be a vector of unit length. Reading 
$\tmop{cl}$ as a map to $C^\infty(V\backslash\{0\},\mathbb{C})$ we see
that
\begin{equation}
  \tmop{cl}(\rho_\ast(\xi_a\xi_b))(x) = \tmop{cl}(\xi_a)(x) \tmop{cl}(\xi_b)(x)
  = \langle x,\rho_\ast(\xi_a) x\rangle
  \cdot \langle x,\rho_\ast(\xi_b) x\rangle,
\end{equation}
which has a meaningful physical interpretation. Namely, if we think of
$\xi_a$ and $\xi_b$ as observables, then in the classical limit the
expectation value of the operator $\xi_a\xi_b$ is given by the product
of the expectation values of $\xi_a$ and $\xi_b$\footnote{This
  remark has to be taken \emph{cum grano salis}, because of the
  possible complex phases on the right-hand side. For probabilities
  one has to take the absolute value squared, which is an implicit
  convention in theoretical physics.}. But this means that the
operators $\xi_a$ and $\xi_b$ are stochastically independent in the
classical limit.

The main point of this chapter is to give an analytical realization 
of this purely algebraic construction, i.e., there will be a
parameter and we will obtain the above classical limit as an analytical
limit when this parameter goes to infinity. This will make the notion
of $\hbar\to 0$ precise in our context.  Here the theme of
non-commutativity vs. commutativity will appear again.

\section{Realizing the Classical Limit as an Analytical Limit}
The Lie algebra $\mathfrak{g}$ can be decomposed as
\begin{equation}
  \mathfrak{g} = \mathfrak{u}_- \oplus \mathfrak{t}^{\mathbb{C}}
  \oplus \mathfrak{u}_+,
\end{equation}
where $\mathfrak{t}^{\mathbb{C}}$ is the Lie algebra of the
complexified maximal torus and $\mathfrak{u}_-$ and $\mathfrak{u}_+$
are unipotent Lie subalgebras corresponding to the positive and
negative roots. We define the groups
\begin{equation}
  U_+=\tmop{exp}(\mathfrak{u}_+) ,\  U_-=\tmop{exp}(\mathfrak{u}_-)\text{, and
  }T^{\mathbb{C}}=\tmop{exp}(\mathfrak{t}^{\mathbb{C}}).
\end{equation}
Recall that the decomposition of the Lie algebra $\mathfrak{g}$ almost
yields a decomposition of $G$. ``Almost'' in this context means that it is
a decomposition of $G\backslash S$, where $S$ is a Zariski-closed set,
\begin{equation}
  \label{eq:u-tu+decomp}
  G = \text{Zarsiki closure of }U_- T^{\mathbb{C}} U_+,
\end{equation}
and even stronger
\begin{equation}
   G\backslash S \simeq U_- \times T^{\mathbb{C}} \times U_+.
\end{equation}

Let us again consider the representation
$\rho_\ast:\mathcal{U}(\mathfrak{g}) \to \tmop{End(V)}$ and choose a
vector of highest weight $v_{max}\in V$. By the definition of
$v_{max}$ we see that $U_+\subset \tmop{Stab}_G(v_{max})$ and
$\rho(T)\subset \mathbb{C}^\ast\cdot v_{max}$. Moreover, the K-orbit
through $[v_{max}]$ agrees with the $G$-orbit through this point,
i.e.\ $K.[v_{max}]=G.[v_{max}]$.

Thus, there exists a Zariski-closed set $A$ in $K.[v_{max}]$ such
that $K.[v_{max}]\backslash A$ is isomorphic to the orbit of $U_-$
through $v_{max}$ in $V$. Therefore, the $U_-$-orbit is isomorphic to
a dense, Zariski-open subset of $Z$ if we identify $Z=K.\lambda$ with
$K.[v_{max}]$ via the momentum map.

We will write $\tmop{cl}$ as composition of two maps $r$ and $s$:
\begin{equation}
  r:i\mathfrak{k}\to \tmop{Vect}(V\backslash\{0\}), \xi \mapsto -\frac{1}{2}
  X_\xi,\text{ with }
  (X_\xi f)(x) = \left.\frac{d}{dt}\right|_{t=0} f(\exp(-\xi t).x)
\end{equation}
and
\begin{equation}
  s:\tmop{Vect}(V\backslash\{0\})\to C^\infty(V\backslash\{0\},\mathbb{C}), X \mapsto
  \frac{1}{N} (XN),
\end{equation}
where $N(x)=\|x\|^2$ is the norm function squared.

Slightly changing the definition of $\tmop{cl}$  to a
map to $C^\infty(V\backslash\{0\},\mathbb{C})$ the definition of the
momentum map (\ref{eq:momentummap}) yields
the following commutative diagram:
\begin{equation}
  \xymatrix{
    i\mathfrak{k} \ar[rr]_r \ar[rrd]_{\tmop{cl}} &&
    \tmop{Vect}(V\backslash\{0\}) \ar[d]_s\\
    &&
    C^\infty(V\backslash\{0\},\mathbb{C})
  }
\end{equation}

Let us explicitly calculate the map $s$ on the $U_-$-orbit through
$v_{max}$:
\begin{equation}
  (X_\xi N)(x) = \left.\frac{d}{dt} \right|_{t=0} N(\exp(-t\xi).x).
\end{equation}
Since $x$ lies on the $U_-$-orbit, there exists a $u\in
U_-$ such that
\begin{equation}
  \label{eq:choiceofu}
  x=u.v_{max}.
\end{equation}
Now we can decompose $\exp(-t\xi)u$ uniquely as
\begin{equation}
  \exp(-\xi t)u = u_-(t)l(t)u_+(t)
\end{equation}
for $t$ in a neighborhood of $0$, where $u_-(t)\in U_-$, $l(t)\in
T^{\mathbb{C}}$ and $u_+(t)\in U_+$. To see this note that we can
decompose the identity and the set of decomposable elements is a
Zariski open set by (\ref{eq:u-tu+decomp}). Using the chain rule and
self-adjointness of $\xi$, we obtain
\begin{equation}
  (X_\xi N)(x) = 2\left\langle x,\left.\frac{d}{dt}\right|_{t=0}\exp(-\xi t).x
  \right\rangle
  = 2\left\langle x,\left.\frac{d}{dt}\right|_{t=0}u_-(t)l(t)u_+(t).v_{max}
  \right\rangle.
\end{equation}
But since $u_+(t)\in U_+\subset \tmop{Stab}_G(v_{max})$ for all $t$ we
have
\begin{equation}
  (X_\xi N)(x) = 2\left\langle x,\left.\frac{d}{dt}\right|_{t=0}u_-(t)l(t).v_{max}
  \right\rangle.
\end{equation}
According to the product rule and using $l(0)=Id$, $u_-(0)=u$ we find
\begin{equation}
  (X_\xi N)(x) = 2\langle x,u\left.\frac{d}{dt}\right|_{t=0}l(t).v_{max} 
  \rangle + 2\left\langle x,\left.\frac{d}{dt}\right|_{t=0} 
  u_-(t).v_{max} \right\rangle.
\end{equation}
Due to the fact that $l(t)\in T$ acts as scalar on $v_{max}$ this can be
simplified as follows
\begin{equation}
  (X_\xi N)(x) = 2\dot{l}(0)\langle x,x\rangle + 2\left\langle
    x,\left.\frac{d}{dt}\right|_{t=0}u_-(t).v_{max} \right\rangle.
\end{equation}

Thus, we can read the right hand side as a differential operator
applied to the norm function. This operator consists of a
multiplication part with $2\dot{l}(0)$ and a vector field part which
is tangential to the $U_-$-orbit. Let $\mathcal{D}(U_-.v_{max})$
denote the algebra of linear differential operators on
$U_-.v_{max}$. We claim that the above procedure affords a map
\begin{equation}
  \label{eq:def_rtilde1}
  \tilde{r} : i\mathfrak{k} \to \mathcal{D}(U_-.v_{max}), 
  \xi \mapsto m_\xi + \xi_{tan}, 
\end{equation}
where $\xi_{tan}$ is the vector field tangent to the $U_-$ orbit
whose one parameter group at $x$ is given by $2\left.\frac{d}{dt}
\right|_{t=0} u_-(t)$ with respect to the above decomposition, and
$m_\xi$ is a smooth function on the $U_-$-orbit with
$m_\xi(x)=2\dot{l}(0)$. The only thing we have to show is that the
construction is independent of the choice of $u$ in
(\ref{eq:choiceofu}). But if we choose $u'$ with
\begin{equation}
  x = u.v_{max} = u'.v_{max}
\end{equation}
then $u'u^{-1}\in\tmop{Stab}_G(v_{max})$. So, $u'=ug$, where
$g\in\tmop{Stab}_G(v_{max})$. But as $g$ acts trivially on $v_{max}$ the
calculation does not change.

The map $\tilde{r}$ will be the crucial point in the following. We
will discuss it from an abstract point of view later on, but first we
extend $\tilde{r}$ to $\mathcal{T}(\mathfrak{g})$ in the following
manner
\begin{equation}
  \label{eq:def_rtilde2}
  \tilde{r}(\xi_{\alpha_1}\otimes\dots\otimes\xi_{\alpha_p}) =
  \tilde{r}(\xi_{\alpha_1}) \circ 
  \dots \circ \tilde{r}(\xi_{\alpha_p}).
\end{equation}
This is well-defined because the $\tilde{\xi_j}$ are linear
differential operators, so they respect scalar multiplication and
addition.

Before we go into the details of the convergence, we need a fact about
the norm.
\begin{theorem}
  \label{thm:norm-multiplicative}
  Let $\lambda$ be the highest weight of the representation $\rho$
  with decomposition into fundamental weights $f_j$ as follows
  \begin{equation}
    \lambda = \sum_{j=1}^r \lambda_j f_j.
  \end{equation}
  Then the squared norm function $N$ on the $U_-$-orbit decomposes
  as
  \begin{equation}
    N(u.v_{max}) = c\cdot N_1(u.v_{max})^{\lambda_1} \cdot \dots \cdot
    N_r(u.v_{max})^{\lambda_r},
  \end{equation}
  where $r$ is the rank of $\mathfrak{g}$ and $N_1,\dots,N_r$ are the
  squared norms of the fundamental unitary representations
  corresponding to the fundamental weights $f_1,\dots,f_r$.
\end{theorem}

\begin{proof}
  For every fundamental representation $\rho_{(j)}$ we have a
  holomorphic line bundle $L_j\to G/B_-$ such that the representation
  of $G$ on $\Gamma_{\tmop{hol}}(G/B_-,L)$ is equivalent to
  $\rho_{(j)}$ (cf.~Appendix Theorem \ref{thm:borel-weil}).

  By induction and Lemma~\ref{tensorbuendellemma},
  we find that the representation with highest weight $\lambda = \sum
  \lambda_j f_j $ is given by the action on the sections of
  \begin{equation}
    L = L_{(1)}^{\lambda_1} \otimes \ldots \otimes 
    L_{(r)}^{\lambda_r}.
  \end{equation}

  Let $h_j$ denote the induced $K$-invariant, hermitian bundle metric
  on $L_j$, which is given in Lemma \ref{lem:norm-lemma--for-bundles},
  and $h$ the induced metric for $L$. 

  Choose a common open covering $\{W_k\}$ of $G/B_-$, such that $L$
  and all $L_j$ are trivializable over each $W_k$. Without loss of
  generality we may assume that $W_1=U_-\cdot[v_{\tmop{max}}]$. Each
  hermitian bundle metric $h_j$ is given by a family $\{m_{k,j}:W_k\to
  \mathbb{R}_+\}$, $h$ by the family $\{m_k:W_k\to \mathbb{R}_+\}$. 
  
  A direct calculation shows that the family $\{m_k':W_k\to
  \mathbb{R}_+\}$ given by
  \begin{equation}
    m_k':= m_{k,1}^{\lambda_1}\cdot\ldots\cdot m_{k,r}^{\lambda_r} 
  \end{equation}
  represents a hermitian, $K$-invariant bundle metric $h'$ on $L$.
  Thus, $h'=ch$ for some positive constant $c$.  Using
  (\ref{buendelmetrik}) we see that the norm on $W_1$ is defined by
  the bundle metric up to this scalar.

  This completes the proof of Theorem \ref{thm:norm-multiplicative}.
\end{proof}

In the following we will consider a highest weight $\lambda=\sum_j
\lambda_jf_j$. If we are given a function like
\begin{equation}
  \frac{u+\bar{u}}{1+\|u\|^2}\lambda_1 + 17\lambda_2
\end{equation}
then we can think of the function as a polynomial in
$\lambda_1,\lambda_2$ where the coefficients are smooth functions. It
is even a homogeneous polynomial of degree $1$.

\begin{notation}
  The ring of smooth functions on $U_-.v_{max}$ is denoted by the
  symbol $R$, i.e.\ $R:=C^\infty(U_-.v_{max},\mathbb{C})$, and the
  ring of polynomials in the $\lambda_j$ with coefficients in $R$ by
  $R[\lambda]$.
\end{notation}
 
The key result of this chapter is the following:
\begin{theorem}
  \label{thm:cl_main_thm}
  Let $\lambda=\sum_j \lambda_j f_j$ be the highest weight of $\rho$
  and assume that at least one $\lambda_j > p$ for a fixed natural
  number $p$.  Furthermore, let $\alpha=\xi_{\alpha_1}\otimes \dots
  \otimes \xi_{\alpha_p}$ be a ``monomial'' element of degree $p$ in the
  generators $\xi_j$ of $\mathfrak{g}$ as chosen above.

  Then $f(\lambda):=\frac{1}{N}\tilde{r}(\alpha)(N)\in R[\lambda]$ and
  $\tmop{deg} f = p$. The homogeneous part of degree $p$ of $f$ is, up
  to a real, multiplicative constant, given by
  $\tmop{cl}(\xi_{\alpha_1})\cdot\dots\cdot\tmop{cl}(\xi_{\alpha_p})$,
  where we view the $\tmop{cl}(\xi_{\alpha_j})$ as elements of
  $R[\lambda]$. Moreover, the constant does not depend on $\alpha$.
\end{theorem}

\begin{proof}
  By definition, every $\tilde{r}(\xi_{\alpha_j})$ is a first order
  partial differential operator. Hence the summands in the
  derivative of $N=N_1^{\lambda_1} \cdot \dots \cdot N_r^{\lambda_r}$,
  after dividing by $N$, are polynomials in $\lambda$ of degree at
  most $p$.  On the other hand, at least one such summand must be a
  polynomial of degree at least $p$. If all were of lesser degree, one of the
  $\xi_{\alpha_j}$ would be multiplication by a constant, which is not
  the case, or the partial derivatives would lower every exponent
  $\lambda_j$ to 0, which yields a contradiction because at least one
  $\lambda_j$ is larger than $p$. This proves the first part of the
  theorem.

  For the second part, we consider the case $p=1$ first.Then there
  is no degree zero term in the polynomial
  $\frac{1}{N}\tilde{r}(\alpha)(N)$ since
  \begin{equation}
    \frac{1}{N}\tilde{r}(\alpha)(N) = \frac{1}{N} r (\alpha)(N)
  \end{equation}
  in the above construction. But $r(\alpha)$ is a vector field
  and contains no multiplicative part, so we have only partial
  derivatives turning $N$ into a homogeneous polynomial of degree 1
  after dividing by $N$. This proves the second statement for $p=1$.

  Let $p \geq 2$ and $\xi_{\alpha_{1}} \otimes \ldots \otimes
  \xi_{\alpha_{p}}$ be given.  We have $\tilde{r}(\xi_a)=c+\sum
  a_j\frac{\partial}{\partial z_j}$ for some $a_j$ and $c$ in some
  coordinate system $\{z_j\}$ on $U_-.v_{max}$. By the induction
  hypothesis
  \begin{equation}
    \tmop{cl}(\xi_{\alpha_2}\otimes \dots \otimes \xi_{\alpha_p}) =
    \tmop{cl}(\xi_{\alpha_2}) \cdot \ldots \cdot
    \tmop{cl}(\xi_{\alpha_p}) + q,  
  \end{equation}
  where $q$ is a polynomial of degree less than $p-1$. Using the
  product rule of differentiation we calculate explicitly
  \begin{multline}
     \tmop{cl}(\xi_{\alpha_{1}} \otimes \ldots \otimes
     \xi_{\alpha_{p}})  = \frac{1}{N} \tilde{r}(\xi_{\alpha_{1}}
     \otimes \ldots \otimes
     \xi_{\alpha_{p}})(N)  
    = 
     \frac{1}{N} \tilde{r}(\xi_{\alpha_{1}})(N
     \tmop{cl}(\xi_{\alpha_2}) \ldots 
     \tmop{cl}(\xi_{\alpha_{p}}) + Nq)  \\ \displaybreak[1] 
    =
     \frac{1}{N} \left( c+\sum a_j\frac{\partial}{\partial z_j}
     \right)(N \tmop{cl}(\xi_{\alpha_2}) \ldots 
     \tmop{cl}(\xi_{\alpha_{p}}) + Nq)  \\ \displaybreak[1] 
    =
     c\cdot\tmop{cl}(\xi_{\alpha_2}) \ldots 
     \tmop{cl}(\xi_{\alpha_{p}}) + q + \frac{1}{N}
     \tmop{cl}(\xi_{\alpha_2})  \ldots
     \tmop{cl}(\xi_{\alpha_{p}}) \left( \sum
       a_j\frac{\partial}{\partial z_j}
     \right)(N) \\
    + 
     \frac{1}{N} \left( \sum a_j\frac{\partial}{\partial z_j}
     \right) (\tmop{cl}(\xi_{\alpha_2}) \ldots
     \tmop{cl}(\xi_{\alpha_p})) + \frac{1}{N} \left( \sum
       a_j\frac{\partial}{\partial z_j} \right) (Nq)
      \\
    = 
     \frac{1}{N} (\tmop{cl}(\xi_{\alpha_2}) \cdot \ldots \cdot
     \tmop{cl}(\xi_{\alpha_p}))\left( \sum a_j\frac{\partial}{\partial
        z_j} +c \right)(N) + \text{ terms of degree less than $p$}.
      \\
    = 
     \tmop{cl}(\xi_{\alpha_1})+
    (\tmop{cl}(\xi_{\alpha_2}) \cdot \ldots \cdot
    \tmop{cl}(\xi_{\alpha_p})) + \text{ terms of degree less than
      $p$}. 
  \end{multline}
   
  Here the first summand is a homogeneous polynomial of degree $p$, as
  claimed.  The remaining summands are certainly of lower degree,
  because each $\tmop{cl}(\xi_b)$ is of degree one and taking the partial
  derivatives can only lower the degree.
\end{proof}

After these preparations we define the classical limit along a ray in
the following way.
\begin{definition}
  Let $\rho:K\to \tmop{U}(V)$ be a non-trivial, irreducible, unitary
  representation of a semisimple, compact Lie group $K$ on a
  finite-dimensional vector space $V$ corresponding to the highest
  weight $\lambda$.

  We call a sequence $\left(\rho_n:K\to \tmop{U}(V_n)\right)_{n\in
    N^*}$ of irreducible, unitary representations, each $\rho_n$ corresponding
  to the highest weight $n\cdot \lambda$, the \textbf{ray through}
  \boldmath $\rho$\unboldmath. For simplicity, we shall always assume
  that $\rho_1=\rho$.
\end{definition}

Let $\xi_H \in \mathcal{T}(\mathfrak{g})$ be an abstractly hermitian
operator and $\xi_1,\dots,\xi_k$ be a basis of $i\mathfrak{k}$.  We
have a unique decomposition into ``monomials'' of $\xi_H = \sum_j a_j
\xi_{j_1}\otimes \dots \otimes \xi_{j_{d(j)}}$, where each $a_j$ is a
complex number. (Keep in mind that these are not monomials in the
usual sense because of the non-commutativity.)

\begin{definition}
  The \textbf{n-th approximation of the classical limit} %
  \nomenclature{$\tmop{cl}_n$}{$n$-approximation of the classical
    limit map}%
  is
  \begin{equation}
    \tmop{cl}_n(\xi_H) = \sum_j a_j \frac{1}{n^{d(j)}} \frac{1}{N}
    \tilde{r}_n(\xi_{j_1}\otimes \dots \otimes 
    \xi_{j_{d(j)}})(N).
  \end{equation}
  Here $\tilde{r}_n$ is defined as in (\ref{eq:def_rtilde1}) and
  (\ref{eq:def_rtilde2}) with respect to the representation $\rho_n$, 
  i.e.\ we substitute every $\lambda_j$ in the resulting polynomials by
  $n\cdot\lambda_j$.
\end{definition}

\begin{theorem}
  Along a ray through the non-trivial, irreducible representation
  $\rho$ the $n$-th approximations of the classical limit converge to
  the classical limit uniformly on compact subsets of
  $U_-.v_{max}$ for every fixed $\xi_H \in
  \mathcal{T}(\mathfrak{g})$, i.e.\
  \begin{equation}
    \tmop{cl}_n(\xi_H) \to \tmop{cl}(\xi_H) \text{ uniformly on compact
      subsets, as }n\to \infty. 
  \end{equation}
\end{theorem}

\begin{proof}
  Decompose $\xi_H$ into its homogeneous parts:
  \begin{equation}
    \xi_H = \sum_j \xi_j
  \end{equation}
  where each $x_j$ is homogeneous of degree $j$. Since $\rho$ is a
  non-trivial representation, at least one $\lambda_j$ in the
  decomposition $\lambda=\sum_j \lambda_j f_j$ is not zero. Because
  $\xi_H$ has only a finite degree, the conditions of Theorem
  \ref{thm:cl_main_thm} are satisfied for all $n$ sufficiently big.
  Applying this theorem to each ``monomial'' in every $\xi_j$ implies
  \begin{equation}
    \tmop{cl}_n(\xi_j) = \tmop{cl}(\xi_j) + \frac{1}{n}
    (\text{terms of lower degree}).
  \end{equation}
  It follows that for any compact set $M$
  \begin{equation}
    \tmop{cl}_n(\xi_H)(x) \to \tmop{cl}(\xi_H)(x) \text{ as }n\to\infty 
  \end{equation}
  for all $x\in M$ uniformly.
\end{proof}

This completes the construction of the classical limit as a
mathematical limit.  The reader might wonder whether the convergence
on a dense, open subset of $Z$ suffices.  Note that $\tmop{cl}$ is
defined on the whole of $Z$, but our $U_-$ chart is
not. Unfortunately, it is not clear that every approximation can be
extended to $Z$, but nevertheless the limit does extend continuously.

Let us now discuss the procedure a more abstractly. The main
step is the substitution of $\tilde{r}$ for $r$ in the definition of the
classical limit.  After this, the other theorems follow from Theorem 
\ref{thm:norm-multiplicative}. But what are these
deformed vector fields $\tilde{r}(\xi)$? In a way this is at least 
in a formal sense similar to a connection in a line bundle plus multiplicative 
function, like in geometric quantization.
Indeed, we have a line bundle here. It is the tautological
bundle $V\backslash\{0\} \to \mathbb{P}(V)$ restricted to
$K.[v_{max}]$.  Furthermore, the $U_-$-orbit can be thought of as a
section of this bundle over the dense open set $U_-.[v_{max}]$.
Since $U_-$ is biholomorphic to some $\mathbb{C}^p$, we get a
chart for the bundle here. In this chart $\tilde{r}$ is in fact just
a connection plus a multiplicative part.

\begin{figure}[ht]
  \centering
  \includegraphics{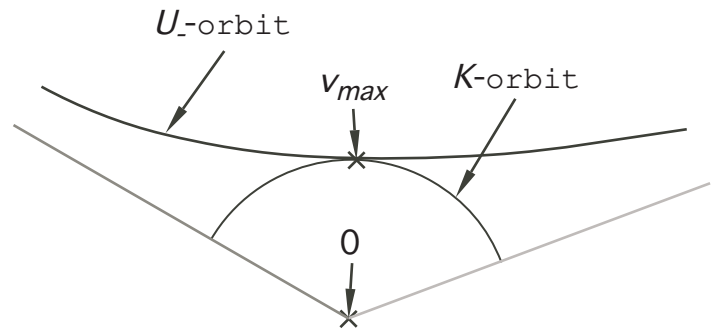}
  \caption{A picture of the $U_-$-section.}\label{fig:bundle}
\end{figure}

A visualization of the situation is provided by Figure
\ref{fig:bundle}. Here we see the origin in $V$ and $v_{max}$. Since
$K$ acts unitarily, the $K$-orbit preserves the metric and is drawn as
a circular arc. The $U_-$-orbit is non-compact and drawn as a very
flat parabola. If we look at this in $\mathbb{P}(V)$, we see that the
$U_-$-orbit is not a global section of the tautological bundle
because the horizontal axis has no intersection with the
$U_-$-orbit.

% LocalWords: complexification adjointness unipotent Zarisiki semisimple Borel
%%% Local Variables: 
%%% mode: latex 
%%% TeX-master: "../thesis" 
%%% End:
% LocalWords:  biholomorphic unitarily Weil equivariant

% \chapter{Spectral statistics of simple Hamotilonian operators}
\chapter[Spectral Statistics of Simple Operators]{Spectral Statistics of Simple Hamiltonian Operators}
\label{chap:spec-stat-simple}
The spectral statistics of simple Hamiltonian operators, i.e., the
nearest neighbor statistics for elements of some semi-simple Lie
algebra, are discussed in this chapter. The main interest is in the
behavior of the spectral statistics in irreducible representations as
the dimension goes to infinity. Thereafter, the notion of rescaling is
introduced and some consequences of the choice of rescaling are given.

\section{A Convergence Theorem for Simple Operators}
In this section we give an estimation on the number of weights of
irreducible representations and in certain cases deduce from it the
convergence of the spectral statistics for simple operators.

Here $K$ always denotes a semi-simple, compact Lie group with a fixed
maximal torus $T$ and a fixed notion of positivity of roots. We write
$W$ for the Weyl group of $K$ with respect to $T$. Further, let $G$ be
the complexification of $K$ and denote the corresponding Lie
algebras by $\mathfrak{g}$ and $\mathfrak{k}$. For any hermitian
matrix $A$ we write $\mu_A$ for the nearest neighbor statistics of
$A$, i.e.,
\begin{equation}
  \mu_A := \mu(X(A))
\end{equation}
as defined in (\ref{eq:def_nearestneighbordistribution}). If $U$ is a
unitary matrix we will write $\mu_U$ for the nearest neighbor
statistics of unitary matrices
(\ref{eq:def_nearestneighbordistribution_circle}), i.e.
\begin{equation}
  \mu_U := \mu_c(X(A)).
\end{equation}
It is clear by the subscript which kind of statistics is meant, so we
use the same abbreviation.

We start with a lemma.
\begin{lemma}
\label{lem:weight-estimate}
Let $\rho_\lambda:K\to U(V_\lambda)$ be an irreducible, unitary
representation with highest weight $\lambda$. Let $\lambda=\sum
\lambda_j f_j$ be the decomposition of $\lambda$ into the basis of
fundamental weights $f_j$.
Then the number $n_\lambda$ of possible weights of $\rho_\lambda$ is
bounded as follows
\begin{equation}
  n_\lambda\ \leq\ \tmop{ord}(W)\cdot\prod_{j}(\lambda_j+1).
\end{equation}
\end{lemma}

\begin{proof}
  Starting from $\lambda$ we get all other weights by subtracting
  multiples of the roots. The lattice of roots is a sublattice of the
  lattice of weights, so we can reach every weight by subtracting
  multiples of the fundamental weights $f_j$.

  There are at most $\prod_j (\lambda_j+1)$ of the such possible
  substractions that give positive weights and every weight is in the
  $W$-orbit of a positive weight, which has at most $|W|$ elements.
\end{proof}

Now we give a rough estimate for the dimension of an irreducible
representation.
\begin{lemma}
\label{lem:dim-estimate}
Under the assumptions of Lemma \ref{lem:weight-estimate} we have the
following inequality for the dimension of $\rho_\lambda$:
\begin{equation}
  \tmop{dim} \rho_\lambda \geq \prod_{\alpha\in\Pi^+,\langle\lambda,\alpha\rangle>0} \frac{\langle \lambda, 
    \alpha \rangle}{\langle\delta,\alpha\rangle}, 
\end{equation}
where $\Pi^+$ denotes the set of positive roots and
$\delta=\frac{1}{2}\sum_{\alpha\in\Pi^+}\alpha$.
\end{lemma}

\begin{proof}
Weyl's dimension formula reads
\begin{equation}
  \tmop{dim} \rho_\lambda = \prod_{\alpha\in\Pi^+} \frac{\langle \delta+\lambda, 
  \alpha \rangle}{\langle\delta,\alpha\rangle} = \prod_{\alpha\in\Pi^+}\left( 1 +\frac{\langle \lambda, 
  \alpha \rangle}{\langle\delta,\alpha\rangle} \right).
\end{equation}
Now, $\langle \lambda,\alpha \rangle \geq 0$ and $\langle \delta,
\alpha \rangle>0$ for all positive roots $\alpha$.  Thus, the
inequality is clear.
\end{proof}

We write $\delta_{\tmop{Dirac}}$ for the Dirac measure with mass $1$
at $0$ and apply these lemmas to the situation of Chapter
\ref{chap:rep-theo-limit} where we looked at rays to infinity. 

\begin{theorem}
  \label{thm:spectr-stat-simpl-ops}
  Let $\rho:K\to \tmop{U}(V)$ be an irreducible representation with highest
  weight $\lambda=\sum \lambda_j f_j$ and the sequence
  $(\rho_n:K\to\tmop{U}(V_n))_{n\in\mathbb{N}^\ast}$ be a ray through $\rho$.

  If $r:=\tmop{rank}(K)\geq 2$ and 
  \begin{equation}
    r< \ ^\#\{\alpha\in\Pi^+ : \langle \alpha,\lambda \rangle>0\} 
  \end{equation}
  then for every 
  $\xi\in i\mathfrak{k}\backslash\{0\} $
  \begin{equation}
    \mu_{\rho_{\ast,m\lambda}(\xi)} \to \delta_{\tmop{Dirac}}\text{
      in $d_{KS}$, as }m\to\infty.
  \end{equation}
\end{theorem}

\begin{proof}
  Let $\xi\in i\mathfrak{k}$ be given. The element
  $i\xi\in\mathfrak{k}$ is conjugated to an element
  $\eta\in\mathfrak{t}=\tmop{Lie}(T)$. We will show that
  \begin{equation}
    p_m:=\frac{\text{number of (different) eigenvalues of }
      \rho_{\ast,m\lambda} (\eta)}{\tmop{dim} \rho_{m\lambda}} \to 0 
  \end{equation}
  as $m\to \infty$. This implies the convergence to $\delta_{Dirac}$
  since the value of $\lim_{s\to 0} \int_0^s d\mu_A$ is
  \begin{equation}
    1-\frac{\text{number of (different) eigenvalues}}{
      \text{number of rows of $A$}}
  \end{equation}
  for any hermitian matrix $A$ by the definition of the nearest
  neighbor statistics. Thus, $\mu_{\rho_{\ast,m\lambda} (\xi)}$ has mass
  $1-p_m$ at zero, which proves the convergence.

  It remains to show the claim about $p_m$. To do so, note that the
  eigenvalues of $\rho_{m\lambda}(\xi)$ are just the values of the
  weights of the representation evaluated at $\xi$. So, it is
  sufficient to prove that the ratio of the different weights and
  the dimension of $\rho_{m\lambda}$ converges to zero.
   
  To show this we combine the inequalities of Lemma
  \ref{lem:weight-estimate} and \ref{lem:dim-estimate}, but first we
  simplify the notation a bit. We denote by $Q$ the set of $\alpha\in
  \Pi^+$, such that $\langle \alpha,\lambda \rangle >0$ and by $q$ the
  cardinality of $Q$. Finally, the number of different weights in
  $\rho_{m\lambda}$ is $n_{m\lambda}$.

  We obtain
\begin{equation}
  \label{eq:convergencerate_for_torus}
  \frac{n_{m\lambda}}{dim{\rho_{m\lambda}}} \leq \frac{\left(\tmop{ord}(W)\cdot\prod_{j=1}^r(\lambda_j+1)\right) m^r}
  {\left(\prod_{\alpha\in Q} \frac{ \langle \lambda,\alpha \rangle }{ \langle\delta,\alpha\rangle }\right) m^q } 
  = c(\lambda)\, m^{r-q}.
\end{equation}
  Here $c(\lambda)$ is a constant, depending only on $\lambda$, and, since
  $r<q$ by the hypothesis, the ratio converges to zero as
  promised. This proves the theorem.
\end{proof}

\begin{remark}
  The number $q$ in the above proof is the complex dimension of the
  coadjoint orbit through $\lambda$, i.e., the complex dimension of
  the classical phase space in the classical limit of Chapter
  \ref{chap:rep-theo-limit}.   
\end{remark}

\begin{corollary}
  The conditions of the above theorem will be automatically satisfied
  if $K$ is simple, $\tmop{rank}{K}\geq 2$, and $\lambda$ lies in the
  interior of the Weyl chamber.
\end{corollary}
\begin{proof}
  First, we remark that $r$ equals the number of positive
  roots for any representation whose highest weight is in the
  interior of the Weyl chamber, since the interior is defined by the
  condition $\langle \lambda,\alpha\rangle >0$ for every simple root
  $\alpha$. But positive roots are positive integer combinations
  of simple roots $\langle \lambda,\alpha\rangle >0$ for all positive
  roots $\alpha$. This completes the proof.
\end{proof}

We now give another corollary.

\begin{corollary}
  Under the assumptions of the theorem let $t_1,\dots,t_p \in
  \mathfrak{g}$ be given such that $\xi= t_1 \otimes \dots
  \otimes t_p\in \mathcal{T}(\mathfrak{g})$ is abstractly
  hermitian in the sense of definition \ref{def:formaladjoint_on_tg}.
  Furthermore, let $\rho_{\ast,m\lambda}$ be the induced Lie algebra
  representation with highest weight $m\lambda$ extended to the
  full tensor algebra. Then
\begin{equation}
    \mu_{\rho_{\ast,m\lambda}(\xi)} \to \delta_{\tmop{Dirac}}\text{  weakly as } m\to\infty,
  \end{equation}
  if $p \cdot r <\ ^\#\{\alpha\in\Pi : \langle\alpha,\lambda\rangle>0\} $.
\end{corollary}

\begin{proof}
  We can assume without loss of generality that all $t_j$ are always
  represented as diagonal matrices and we proceed by induction. From
  (\ref{eq:convergencerate_for_torus}) it follows that for each $t_j$
  the number of its eigenvalues $n_{j,m\lambda}$ divided by the
  dimension is smaller than $c(\lambda)m^{r-q}$.  But the maximal
  number of eigenvalues in a product of diagonal matrices is just the
  product of the number of eigenvalues of each matrix.  Thus, we have
  a numerator $m^{rp}$ here instead of $m^r$ in
  (\ref{eq:convergencerate_for_torus}). But by assumption $rp<q$,
  i.e.\ the number of eigenvalues of the product divided by the
  dimension is decreasing faster than $1/m$.

  This proves the corollary.
\end{proof}

\section{Rescaling}
In this section we discuss the notion of rescaling. This concept
appeared already in Chapter \ref{chap:rep-theo-limit}.  There the
classical limit along rays
$(\rho_m:K\to\tmop{U}(V_m))_{m\in\mathbb{N}^\ast}$ through a given
representation $\rho$ was considered and the scaling was given by
substituting $\frac{1}{m}\xi_j$ for $\xi_j$. Since we are interested
in the problem of scaling in general, we define the notion of a
rescaling map abstractly.

Let $\mathcal{U}(\mathfrak{g})$ denote the universal enveloping
algebra of $\mathfrak{g}$ and $\dagger$ the formal adjoint
(cf. Appendix). We choose a fixed basis $\xi_1,\dots,\xi_n$ of
$\mathfrak{g}$ and write the elements of $\mathcal{U}(\mathfrak{g})$
as ordered polynomials in the $\xi_j$. Furthermore the multiindex notation
$\Xi^I$ will be used for $\xi_1^{i_1}\dots \xi_n^{i_n}$. 

The basic problem can be seen if one considers the hermitian
operators $\xi$ and $\xi\eta$ in a sequence of irreducible
representations. As the dimensions of the representations increase the
maximal eigenvalues of $\xi\eta$ will in general grow faster than
those of $\xi$. In principle, we would like the rate of growth to be
the same, including the option of no growth at all. This motivates the
following definition.
\begin{definition}
  A \textbf{rescaling map} $r_\rho$ for the irreducible representation
  $\rho:K \to \tmop{U}(V)$ is given by a map
  \begin{equation}
    r_\rho:\mathcal{U}(\mathfrak{g}) \to \mathcal{U}(\mathfrak{g}),
      \sum a_I \Xi^I \mapsto \sum \frac{1}{s^{|I|}} a_I \Xi^I,
  \end{equation}
  where $s$ a positive integer number.
\end{definition}

\begin{lemma}
  Every rescaling map $r_{\rho}$ is linear, injective and
  compatible with $\dagger$.
\end{lemma}

\begin{proof}
  This follows directly from the definition of $r_{\rho}$.
\end{proof}

Of all possible scalings the most natural one is the scaling by
inverse dimension since we have no other natural quantity associated
to arbitrary sequences of irreducible representations.
\begin{definition}
  Let $\tmop{Irr}(K)$ denote the set of equivalence classes of
  irreducible, unitary representations of $K$.  The \textbf{rescaling by
  inverse dimension} is the family of rescaling maps $(i_\rho)_{\rho\in
    \tmop{Irr}(K)}$ given by
  \begin{equation}
    i_\rho:\mathcal{U}(\mathfrak{g}) \to \mathcal{U}(\mathfrak{g}),
     \ \sum a_I X^I \mapsto \sum \frac{1}{(\tmop{dim} \rho)^{|I|}} a_I X^I
  \end{equation}
  for each $\rho\in\tmop{Irr}(K)$.
\end{definition}

If we are considering rays through a fixed irreducible representation
with highest weight $\lambda$, then we have another natural
quantity: the parameter $m$ for each $\rho_{m\lambda}$.

\begin{definition}
  Let $\rho:K\to \tmop{U}(V)$ be an irreducible representation with highest
  weight $\lambda=\sum \lambda_j f_j$ and the sequence
  $(\rho_m:K\to\tmop{U}(V_m))_{m\in\mathbb{N}^\ast}$ be a ray through $\rho$.

  The \textbf{rescaling by inverse parameter} is the family
  of rescaling maps $(p_{\rho_{m\lambda}})$ given by
  \begin{equation}
    p_{\rho_{m}}:\mathcal{U}(\mathfrak{g}) \to \mathcal{U}(\mathfrak{g}),
    \ \sum a_I X^I \mapsto \sum \frac{1}{m^{|I|}} a_I X^I.
  \end{equation}
\end{definition}

\subsection{Rescaling and Spectral Statistics}
In the first section we considered simple operators only, i.e.
Lie algebra elements. Rescaling has no effect in this case since for
any self-adjoint matrix $A$ and any $c>0$
\begin{equation}
  \mu_A = \mu_{c\cdot A}.
\end{equation}
But rescaling has an effect if we consider operators whose monomial
parts have different degrees, e.g.\
\begin{equation}
  \xi + \eta^2 \in \mathcal{U}(\mathfrak{g}).
\end{equation}

Recall that for a highest weight $\lambda$ the set $Q$ is defined as
$Q=\{\alpha\in\Pi_+\,:\,\langle \alpha,\lambda\rangle\}$ and $q=\,^\#
Q$. We state the following lemma:
\begin{lemma}
  \label{lem:resc-spectr-stat}
  Let $\xi_H=\sum_{I} a_{I}\Xi^I \in \mathcal{U}(\mathfrak{g})$ be
  given with $\xi_H^\dagger=\xi_H$ and consider the ray 
  $(\rho_m:K\to\tmop{U}(V_m))_{m\in\mathbb{N}^\ast}$ through an
  irreducible representation $\rho:K\to \tmop{U}(V)$ of highest weight
  $\lambda$.

  Then
  \begin{equation}
    \label{eq:resc-spectr-stat}
    \| \rho_{\ast,m} (i_{\rho_m}(\xi_H)) \|_{\tmop{End}(V_m)} \leq
    c_1(\lambda) \sum_{I} |a_{I}| c_2(\lambda)^{|I|} \cdot m^{|I|-q|I|}
  \end{equation}
  where the $c_j(\lambda)$ are constants depending only on $\lambda$ and
  $\|\cdot \|_{\tmop{End}(V_m)}$ denotes the operator norm on
  $\tmop{End}(V_m)$.
\end{lemma}

\begin{proof}
  We use the explicit construction of irreducible representations by
  Borel-Weil. For this let
  \begin{equation} 
       S_j = (s_1^{(j)}, \ldots, s_{d (j)}^{(j)}), j = 1, \ldots, r 
  \end{equation}
  denote a basis of the $j$-th fundamental representation. These are
  holomorphic sections in a holomorphic line bundle
  \begin{equation} 
     L_j \rightarrow G / B 
  \end{equation}
  where $B$ is a Borel subgroup of $G$ and $L = G \times_{\chi_{_j}}
  \mathbb{C}$, such that $\chi_j : B \rightarrow \mathbb{C}$ is the
  exponentiated character of the fundamental weight $\lambda_j$. The
  irreducible representation with highest weight $\lambda$ is then given by
  the action on sections of the line bundle
  \begin{equation} 
     L = L_1^{\otimes \lambda_1} \otimes \ldots \otimes L_j^{\otimes \lambda_j}
     \rightarrow G / B . 
  \end{equation}
  By the theorem of Borel-Weil the tensors of the form
  \begin{equation}
    \label{eq:tensors_in_representation}
    S_1^{I_1} \otimes \ldots \otimes S_r^{I_r}, 
  \end{equation}
  with $I_1, \ldots, I_r$ multiindices of degree $|I_j | = \lambda_j$
  constitute a generating system of the space of sections.
  
  Without loss of generality we may take a basis $\xi_1, \ldots, \xi_n$ of
  $\mathfrak{g}$, such that $\xi_1$ is represented by a diagonal hermitian
  matrix of spectral norm 1 in every fundamental representation. Since the
  operator norm is equal to the spectral norm, we wish to give an estimate for
  the maximal absolute value of an eigenvalue of $\xi_1$ in
  $\rho_{\ast,\lambda}$.
  
  But on the generating system of vectors given by
  (\ref{eq:tensors_in_representation}) the action is on each
  factor separately, so we have
  \begin{equation} 
    \| \rho_{\ast,\lambda} (\xi_1)\| \leqslant \lambda_1 + \ldots + \lambda_r = : |\lambda|. 
  \end{equation}
  Clearly, the same argument can be carried out for $\xi_2, \ldots,
  \xi_n$. So we have the following estimate
  \begin{equation}
    \| \rho_{\ast,m} (\xi_j)\| \leqslant m (\lambda_1 + \ldots + \lambda_r) = m|\lambda|
    \label{eigvestimate}
  \end{equation}
  for all $j = 1, \ldots, n$.
  
  Now, consider $\gamma = \sum_I a_I X^I$. Then
  \begin{equation} 
    \label{eq:operatornorm_estimation}
  \| \tilde{\rho}_{\ast,m} (i_{\rho_m} (\gamma) t)\|_{\tmop{End} (V_m)}
     \leqslant \sum_I \frac{1}{(\dim \rho_k)^{|I|}} \| \rho_{\ast,m}
     (\xi_1)\|_{\tmop{End} (V_m)}^{i_1} \cdot \ldots \cdot \| \rho_{\ast,m}
     (\xi_n)\|^{i_n}_{\tmop{End} (V_m)} . 
  \end{equation}
  Using the estimates given by (\ref{eigvestimate}) and
  Lemma \ref{lem:dim-estimate}, we see that
  \begin{equation}
    \label{eq:norm-estimate}
    \| \tilde{\rho}_{\ast,m} (i_{\rho_m} (\gamma) t)\|_{\tmop{End} (V_m)}     
    \leq \sum_I \frac{|a_I |}{C \cdot m^{q|I|}} m^{|I|} \cdot |\lambda|^{|I|} 
    = C' \sum_I |a_I | |\lambda|^{|I|} m^{|I|-q|I|}, 
  \end{equation}
  where $C$ and $C'$ are constants depending only on $\lambda$, which
  completes the proof.
\end{proof}

We use this lemma to prove the following theorem.
\begin{theorem}
  Consider the ray $(\rho_m:K\to\tmop{U}(V_m))_{m\in\mathbb{N}^\ast}$
  through an irreducible representation $\rho:K\to \tmop{U}(V)$ of
  highest weight $\lambda$ and assume $q>2$.

  Then for all $\xi_H=\eta + \sum_{|I|\geq 2} a_{I}\Xi^I \in
  \mathcal{U}(\mathfrak{g})$ with $\eta\in\mathfrak{g}\backslash\{0\}$ and
  $\xi_H^\dagger=\xi_H$
  \begin{equation}
    d_{KS}(\mu_{\rho_{\ast,m}(i_{\rho_m}(\xi_H))},\mu_{\rho_{\ast,m}(\eta)}) \to 0 \text{ as }m\to\infty.
  \end{equation}
\end{theorem}

\begin{proof}
  We claim, that
  \begin{equation}
    \label{eq:dimension_estimation_for_algebra}
    \lim_{m\to\infty}(\tmop{dim} V_m)\cdot \left\|\rho_{\ast,m}\left(i_{\rho_m} 
    \left( \sum_{|I|\geq 2} a_{I}\Xi^I \right) \right) \right\|_{\tmop{End(V_m)}}=0.
  \end{equation}
  This implies the theorem, because the nearest neighbor statistics
  for hermitian matrices are scaling invariant, i.e.
  \begin{equation}
    \mu_{(\tmop{dim} V_n)\rho_{\ast,m}(\xi_H)}=\mu_{\rho_{\ast,m}(\xi_H)}
  \end{equation}
  and $(\tmop{dim} V_m)\rho_{\ast,m}(i_{\rho_m}\eta)=\rho_{\ast,m}(\eta)$. Thus,
  \begin{equation}
    \lim_{m\to\infty}\| (\tmop{dim} V_n)\rho_{\ast,m}(i_{\rho_m}(\xi_H)) - 
    \rho_{\ast,m}(\eta) \|_{\tmop{End}(V_m)}=0.
  \end{equation}
  It remains to proof (\ref{eq:dimension_estimation_for_algebra}). But
  by (\ref{eq:norm-estimate}) we obtain
  \begin{equation}
     (\tmop{dim} V_m)\sum_{|I|\geq 2} a_{I}\left\|\rho_{\ast,m}\left(i_{\rho_m} 
    \left(\Xi^I \right) \right) \right\|_{\tmop{End(V_m)}}
    \leq C \sum_I |a_I | |\lambda|^{|I|} m^{|I|-(q-1)|I|}, 
  \end{equation}
  where $C$ is a constant. Since $\lambda$ is fixed and $q>2$ the
  right hand side converges to zero.
\end{proof}
 
So, we only have to study the convergence of $\mu_{\rho_{\ast,m}(\eta)}$
to gain information about the convergence of the nearest neighbor
distribution of the whole operator under rescaling by inverse
dimension. For example, we may use Theorem
\ref{thm:spectr-stat-simpl-ops}.

\subsection[Rescaling and $\exp$]{Rescaling and \boldmath$\exp$\unboldmath}
Rescaling can affect the limit measure of exponentiated operators as
shown in the following lemma.
\begin{lemma}
  \label{lem:mainlemma_algebra}Let $\rho_k : K \rightarrow
  \tmop{U} (V_k)$, $k \in \mathbb{N}$, be a sequence of irreducible,
  unitary representations and $\gamma \in \mathcal{U}(\mathfrak{g})$
  with $\gamma^\dagger=\gamma$.
  
  Let us assume that
  \begin{equation}
    \label{eq:normcondition}
    \lim_{n \rightarrow \infty} \| \tilde{\rho}_{\ast,k} (r_{\rho_k} (\gamma)
    t)\|_{\tmop{End} (V_k)} = 0 \text{ for all } t > 0,
  \end{equation}
  where $\| \cdot \|_{\tmop{End} (V_k)}$ denotes the operator norm on
  $\tmop{End} (V_k)$.
  
  Then $\mu_{\tmop{exp}(\rho_{\ast,k} (r_{\rho_k} (\gamma))t)}$ does not converge to 
  any Borel measure $\mu$ on the positive real line with
  \begin{equation}
    \int_0^{\frac{1}{2\pi}} d\mu < 1
  \end{equation}
  as $n$ goes to infinity for any $t>0$. In particular it does not
  converge to $\mu_{\tmop{Poisson}}$ or $\mu_{\tmop{CUE}}$.
\end{lemma}

\begin{proof}
  For simplicity set $\gamma_k = r_{\rho_k} (\gamma)$ and let $t>0$ be
  fixed. Now by (\ref{eq:normcondition}) we see that starting from a
  sufficiently large $k_0$ the spectrum of $\rho_{\ast,k}(\gamma_k)t$
  is in the interval $[-\pi,-\pi]$.

  Now we may consider a subsequence of $\rho_{k_j}$ such that the
  spectrum of $\rho_{\ast,k}(\gamma_k)t$ is in the interval
  $]-\frac{1}{2j},-\frac{1}{2j}[$. Analogously to the counterexample
  in Remark \ref{rm:exp_under_mu_counterexample} in Chapter
  \ref{chap:gen-level-spacings}, one proves that a limit measure must
  necessarily have the whole mass between $0$ and $1/2\pi$.
\end{proof}

The following theorem states that rescaling by inverse dimension will
destroy convergence to $\mu_{\tmop{Poisson}}$ in many cases.

\begin{theorem}
  \label{thm:dimestimation}
  Choose a fixed irreducible, unitary representation $\rho_{\lambda} :
  K \rightarrow \tmop{U} (V_{\lambda})$ with highest weight $\lambda$,
  where $\lambda = \lambda_1 f_1 + \ldots + \lambda_r f_r$ is the
  decomposition into fundamental weights with every $\lambda_j
  \geqslant 0$.

  Let $\tmop{rank} (\mathfrak{g}) \geqslant 2$ and
  assume that at least two $f_j$ are positive. Then for every $\gamma \in
  \mathcal{U}(\mathfrak{g})$ without constant term
  \begin{equation} \lim_{k \rightarrow \infty} \| \tilde{\rho}_{\ast,k} (i_{\rho_k} (\gamma)
     t)\|_{\tmop{End} (V_k)} = 0 \text{ for all } t > 0 \label{normcondition},
  \end{equation}
  where $\rho_k : K \rightarrow \tmop{U} (V_k)$ is an
  irreducible representation with highest weight $k \cdot \lambda$ and $\|
  \cdot \|_{\tmop{End} (V_k)}$ is the usual operator norm in $\tmop{End}
  (V_k)$.
\end{theorem}

\begin{proof}
  Apply Lemma \ref{lem:resc-spectr-stat} and note that the right
  hand side of (\ref{eq:resc-spectr-stat}) converges to zero. 
\end{proof}

\begin{corollary}
  Under the above assumptions $\mu_{\tmop{exp}(\rho_{\ast,k}
    (r_{\rho_k} (\gamma))t)}$ does not converge to the measures 
  $\mu_{\tmop{Poisson}}$ or $\mu_{\tmop{CUE}}$.
\end{corollary}

\begin{proof}
  This follows from Lemma \ref{lem:mainlemma_algebra}, since we proved that
  (\ref{eq:normcondition}) is fulfilled.
\end{proof}

\begin{remark}
  Note that there is an obvious counterexample to Theorem
  \ref{thm:dimestimation} if the rank of $\mathfrak{g}$ is 1. Namely, the
  irreducible representation of $\mathfrak{s}\mathfrak{l}(2, \mathbb{C})$ on
  the homogeneous polynomials in two indeterminates.
  
  Take $\xi = \tmop{diag} (1, - 1)$. Then $\| \rho_k (\xi)\|_{V_k} = k$
  where $V_k$ is the vector space of homogeneous polynomials of degree $k$.
  Therefore $\dim \rho_k (\xi) = k + 1$. We see that
  \begin{equation} 
     \| \rho_k (r_k (\xi))\|_{V_k} = \frac{k}{k + 1} \rightarrow 1. 
  \end{equation}
\end{remark}

The reader may wonder what happens in the case of the rescaling by
inverse parameter as in Chapter \ref{chap:rep-theo-limit}. There is no
analogue of Theorem \ref{thm:dimestimation} in this case, because
the denominator in (\ref{eq:operatornorm_estimation}) scales like the
numerator, so there is no convergence to zero.

In fact, the statements of this chapter can be made more general by
allowing rescaling maps which decrease operators faster than the
rescaling by inverse parameter. The theo\-rems will still be true in
this case, although some corrections to the constants will be
required.

%%% Local Variables:  
%%% mode: latex 
%%% TeX-master: "../thesis" 
%%% End:
% \chapter{Spectral statistics of generic Hamiltonian operators}
\chapter[Spectral Statistics of Generic Operators]{Spectral Statistics of Generic Hamiltonian Operators}
\label{chap:spec-stat-gen}
Having studied the spectral statistics of simple Hamiltonian
operators, i.e., simple ``polynomials'' of Lie algebra elements 
in irreducible representations, we are now interested in 
more complicated operators.

In Chapter \ref{chap:rep-theo-limit} ``polynomials'' in some basis of
the Lie algebra were considered, which gave rise to Hamiltonians. But
for a more analytic treatment of the matter, we investigate the
spectral statistics in a completion of the polynomial algebra. Note
that such a completion was already implicitly used in \cite{haakekus},
where the authors used the sine of a Lie algebra element.

Thereafter we will define the notion of a generic Hamiltonian operator and
prove that the irreducible representations of the flows through
the generic operators have spectral statistics converging to
$\mu_{Poisson}$ under special assumptions on the dimensions of the
representation spaces.

We will use the following notation throughout this chapter. Let $K$
denote a compact semi-simple Lie group with complexification $G$. The
corresponding Lie algebras are called $\mathfrak{k}$ and
$\mathfrak{g}$. Every representation of $K$ will be assumed to be
continuous, finite-dimensional and unitary. The $K$-invariant inner
product will be denoted by $\langle\cdot,\cdot\rangle$ without putting
the representation space into the notation. It will be clear by the
arguments or by the context which representation space is meant.

\section[Topology and Completion of $\mathcal{U}(\mathfrak{g})$]{Topology and Completion of \boldmath$\mathcal{U}(\mathfrak{g})$\unboldmath}
\label{sec:topol-compl-mathc}
In this section we introduce a topology on the
universal enveloping algebra $\mathcal{U} (\mathfrak{g})$ and
complete it to a Fr\'echet space. To do so, choose a basis $\xi_1,
\ldots, \xi_n$ of $\mathfrak{g}$. By the Poincar\'e-Birkhoff-Witt
Theorem we have a vector space isomorphism
\begin{equation} 
  \psi : \mathbb{C}[X_1, \ldots, X_n] \to \mathcal{U}
  (\mathfrak{g}) 
\end{equation} 
given by substituting $\xi_i$ for
$X_i$ in every polynomial $p$ in which we have ordered the
indeterminates in each monomial lexicographically. Note that this
ordering is necessary since $\psi$ is only a vector space
isomorphism, but not an algebra isomorphism.
% because $\mathcal{U}(\mathfrak{g})$ is non-abelian.

We use $\psi$ to give a topology to $\mathcal{U} (\mathfrak{g})$ by
the natural embedding of $\mathbb{C}[X_1, \ldots, X_n]$ into the
algebra of holomorphic functions $\mathcal{O}(\mathbb{C}^n)$.%
\nomenclature{$\mathcal{O}(\mathbb{C}^n)$}{The algebra of holomorphic
  functions on $\mathbb{C}^n$}%

It is a well-known fact that $\mathcal{O}(\mathbb{C}^n)$ is a
Fr\'echet space with respect to the topology of uniform convergence on
compact subsets of $\mathbb{C}^n$. If we change the basis of
$\mathfrak{g}$ to $\eta_1, \ldots, \eta_n$ we obtain \emph{a priori}
another completion of $\mathbb{C}[X_1, \ldots X_n]$. But changing the
basis is nothing more than a linear change of coordinates, yielding an
induced linear homeomorphism of Fr\'echet spaces. So, a different
choice of basis does not change the topology.

\begin{remark}
  \label{sup_of_coeffs}If a sequence of holomorphic functions on
  $\mathbb{C}^n$ converges to zero in the Fr\'echet topology, then the
  suprema of the coefficients in the Taylor expansion around the origin
  also converge to zero.
\end{remark}

\begin{proof}
  Let $(f_j)_{j \in \mathbb{N}}$ be a sequence of holomorphic functions with
  Taylor expansion $f_j = \sum_I a_I^{(j)} X^I$, where $I$ is a multiindex
  with the usual conventions. 
  
  By the general Cauchy integral formula in several variables we see that
  \begin{equation} a_I^{(j)} = \frac{1}{(2 \pi i)^n} \oint_{\zeta \in \partial P} \frac{f_j
     (\zeta)}{\zeta^{I + (1, \ldots, 1)}} d \zeta , \end{equation}
  where $P$ is the unit polycylinder in $\mathbb{C}^n$ and $\partial P$ its
  distinguished boundary. From this we obtain
  \begin{equation} |a_I^{(j)} | \leqslant \sup_{\zeta \in \partial P} |f_j (\zeta) |. \end{equation}
  The right hand side does not depend on $I$, so the inequality holds for the
  supremum of the $|a_I^{(j)}|$ for a fixed $j$, but the $f_j$ converge
  uniformly on compact sets, especially on $\partial P$.
\end{proof}

Let $\rho_\ast: \mathfrak{g} \rightarrow \tmop{End} (V)$ be an
irreducible representation on a finite-dimensional complex vector
space $V$. This map extends to an irreducible representation of $U
(\mathfrak{g})$, which we will again call $\rho_\ast$.

\begin{proposition}
  The map $\rho_\ast : \mathcal{U} (\mathfrak{g}) \rightarrow \tmop{End} (V)$
  extends to a continuous, surjective, linear map
  \begin{equation} \tilde{\rho}_\ast : \mathcal{O}(\mathbb{C}^n)
    \rightarrow \tmop{End} (V)
  \end{equation}
  with respect to the above completion of $\mathcal{U}
  (\mathfrak{g})$, where the topology on $\tmop{End}(V)$ is given by
  the operator norm with respect to some norm on $V$.
\end{proposition}

\begin{proof}
  Let $f = \sum a_I X^I  \in \mathcal{O}(\mathbb{C}^n)$ be given. We
  define 
  \begin{equation}
    \tilde{\rho}_\ast (f) = \sum a_I \rho_\ast (\xi_1)^{i_1} \ldots \rho_\ast
  (\xi_n)^{i_n}.
  \end{equation}
  By the basic inequality for the operator norm
  \begin{equation} 
    \|AB\| \leqslant \|A\| \cdot \|B\|\quad\forall A, B \in \tmop{End} (V)
  \end{equation}
  it follows that 
  \begin{equation} \|a_I \rho_\ast (\xi_1)^{i_1} \ldots \rho_\ast (\xi_n)^{i_n} \| 
     \leqslant |a_I |b_1^{i_1} \ldots b_n^{i_n} 
  \end{equation}
  for $b_i := \| \rho_\ast (\xi_i)\|$. This series is convergent since
  $f\in\mathcal{O}(\mathbb{C}^n)$. Moreover, $\tilde{\rho}_\ast$ is
  linear.  To show the continuity, it suffices to show that
  $\tilde{\rho}_\ast$ is continuous at zero.
  So let $(f_j)_{j \in \mathbb{N}}$ be a sequence of holomorphic functions on
  $\mathbb{C}^n$ converging to zero uniformly on compact subsets. We must
  show that
  \begin{equation} 
     \lim_{j \rightarrow \infty} \tilde{\rho}_\ast (f_j) = 0,
  \end{equation}
  but this is the claim that
  \begin{equation} 
     \| \sum a_I^{(j)} \rho_\ast (\xi_1)^{i_1} \ldots \rho_\ast (\xi_n)^{i_n} \|
     \rightarrow 0.
  \end{equation}
  Note that
  \begin{equation} 
     \| \sum a_I^{(j)} \rho_\ast (\xi_1)^{i_1} \ldots \rho_\ast (\xi_n)^{i_n} \|
     \leqslant \sum |a_I^{(j)} |b_1^{i_1} \ldots b_n^{i_n}  . 
  \end{equation}
  Again the right-hand side converges to zero because the $\xi_i$ can
  be chosen such that $|b_i | \leqslant \frac{1}{2}$ for all $i \in
  \{1, .,, n\}$, and the right-hand side is less or equal to
  \begin{equation} 
     \sup |a_I^{(j)} | \sum  \frac{1}{2^{|I|}}, 
  \end{equation}
  which converges to zero according to Remark \ref{sup_of_coeffs}. We
  can then scale back to the original $\xi_i$, which is just an
  isomorphism of Frech\'et spaces.
  
  To see that $\tilde{\rho}_\ast$ is surjective we use the Lemma of
  Burnside which states that $\rho_\ast : \mathcal{U}
  (\mathfrak{g}) \to \tmop{End} (V)$ is already surjective.
\end{proof}

\section[A Notion of Hermitian Operators for
$\mathcal{O}(\mathbb{C}^n)$]{A Notion of Hermitian Operators for \boldmath$\mathcal{O}(\mathbb{C}^n)$\unboldmath}

In the following a notion of self-adjointness or hermitian operators for
$\mathcal{O}(\mathbb{C}^n)$ will be required. For this we will extend
the definition of $\dagger$ on $\mathcal{U}(\mathfrak{g})$ by continuity.
\begin{lemma}
  The map $\dagger$ extends to a continuous involution of $\mathcal{O}(\mathbb{C}^n)$.
\end{lemma}

\begin{proof}
  We choose a basis of $\mathfrak{g}$ in the following way. First, 
  fix a maximal torus $\mathfrak{t}$ in $\mathfrak{g}$. Let $\tau_1,
  \ldots, \tau_r$ be a basis of the torus such that $\tau_i^{\dagger}
  = \tau_i$ for all $i$. Then choose a system $\Pi$ of positive roots
  and a basis $\xi_{\alpha}$ of the root spaces
  $\mathfrak{g}_{\alpha}$ for $\alpha \in \Pi$ such that
  \begin{equation} 
    \xi_{\alpha}^{\dagger} = \xi_{- \alpha} . 
  \end{equation}
  With this basis, $\dagger$ operates on the basis elements just by
  permutation.
  
  Let $f = \sum_I a_I X^I$ be in $\mathcal{O}(\mathbb{C}^n)$. We define
  $f^{\dagger} := \sum_I \bar{a}_I (X^I)^{\dagger} $. Clearly,
  $f^{\dagger}$ is again everywhere convergent because we just changed the
  order of the summation and conjugated each coefficient.
  
  Let $(f_j)_{j \in \mathbb{N}}$ be a
  sequence of holomorphic functions on $\mathbb{C}^n$ converging to zero
  uniformly on compact subsets. To show that $\dagger$ is continuous, 
  we must show that
  \begin{equation} 
     \lim_{j \rightarrow \infty} (f_j^{\dagger}) = 0. 
  \end{equation}
  But since in each $f_j^{\dagger}$ we have only changed the order of the
  summands and conjugated to coefficients, this is also a series of
  holomorphic functions converging uniformly on compact subsets.
  
  As stated before, the choice of basis has no effect on the topology.
\end{proof}

We define the notion of an abstractly hermitian operator as follows.
\begin{definition}
  $f \in \mathcal{O}(\mathbb{C}^n)$ is called an \emph{abstractly
    hermitian operator} if $f^{\dagger} = f$. The set of all abstract
  hermitian operators is denoted by $\mathcal{H}$.
\end{definition}

Note that this definition is compatible with the one given for the tensor
algebra in the Appendix. 

\begin{remark}
  $\mathcal{H}$ is a closed subspace of $\mathcal{O}(\mathbb{C}^n)$ and as
  such is a Fr\'echet space.
\end{remark}

\begin{proof}
  The linear map $\dagger - \tmop{id}_{\mathcal{O}(\mathbb{C}^n)}$ is
  continuous and $\mathcal{H}$ is its kernel.
\end{proof}

\begin{lemma}
  Let $\rho : K \to \tmop{U}(V)$ be an irreducible unitary representation and 
  $\rho_\ast: \mathcal{U} (\mathfrak{g}) \to \tmop{End} (V)$ the
  induced representation with extension $\tilde{\rho}_\ast :
  \mathcal{O}(\mathbb{C}^n) \rightarrow \tmop{End} (V)$. Then the restriction of
  $\tilde{\rho}_\ast$ to $\mathcal{H}$ is surjective onto the subspace of
  self-adjoint linear operators of $V$.
\end{lemma}

\begin{proof}
  For $A\in\tmop{End}(V)$ we denote by $A^\dagger$ the conjugate
  transpose of $A$.  We remark that by the definition of $\dagger$ we
  have
  \begin{equation}
    \rho_\ast (\xi)^{\dagger} = \rho_\ast (\xi^{\dagger}) \ \forall \xi
    \in \mathfrak{g}. \label{dagger-equation}
  \end{equation}
   Therefore
  \begin{equation} 
    \rho_\ast (\mathcal{H} \cap \mathcal{U}
    (\mathfrak{g})) \subset \text{self-adjoint operators in End(V)}
    .
  \end{equation}

  To show that the restriction is surjective, consider a self-adjoint operator
  $A \in \tmop{End} (V)$. Since $\rho_\ast$ is surjective, we find an $\alpha
  \in \mathcal{H} \cap U (\mathfrak{g})$, such that $\rho_\ast (\alpha) = A$. By
  (\ref{dagger-equation}) it follows that
  \begin{equation} 
     \rho_\ast (\alpha^{\dagger}) = \rho_\ast (\alpha)^{\dagger} = A^{\dagger} =
     A. 
  \end{equation}
  Therefore we see that
  \begin{equation} 
    \rho_\ast \left( \frac{1}{2} (\alpha + \alpha^{\dagger})\right) = 
    \frac{1}{2} \rho_\ast
    (\alpha) + \frac{1}{2} \rho_\ast (\alpha^{\dagger}) = \frac{1}{2} A +
    \frac{1}{2} A = A. 
  \end{equation}
  But
  \begin{equation} 
     \frac{1}{2} (\alpha + \alpha^{\dagger}) \in \mathcal{H} \cap U
     (\mathfrak{g}), 
  \end{equation}
  so the restriction of $\tilde{\rho}_\ast$ to $\mathcal{H}$ is surjective.
\end{proof}

\section{Examples of Convergence}
\label{sec:examples-convergence}
In this section we will give a class of examples for the convergence
of nearest neighbor statistics of abstractly hermitian operators in
suitable sequences of irreducible representations.

Before these examples are considered we briefly discuss the effect of
holomorphic maps on operators. Consider a holomorphic map
$f:\mathbb{C}\to \mathbb{C}$. It induces a map
\begin{equation}
  \tilde{f}: \mathcal{O}(\mathbb{C}^n)\to \mathcal{O}(\mathbb{C}^n),
  \ g\mapsto f \circ g.   
\end{equation}
Let $\rho:K\to \tmop{U}(V)$ be an irreducible representation and
$\xi\in \mathcal{O}(\mathbb{C}^n)$ be a fixed operator. We are
interested in the spectrum of $\tilde{\rho}_{\ast}(\tilde{f}(\xi))$.

\begin{remark}
  \begin{equation}
    \label{eq:specunderanalyticfct}
    \tmop{Spec} (\tilde{\rho}_{\ast}(\tilde{f}(\xi))) = 
    f(\,\tmop{Spec} (\tilde{\rho}_{\ast}(\xi))\, ).
  \end{equation}
\end{remark}

\begin{proof}
  Let $\sum_j b_j z^j$ be the power series expansion for $f$ at
  zero. Since $\tilde{\rho}_{\ast}$ is continuous, it follows that
  \begin{equation}
    \tilde{\rho}_{\ast}(\tilde{f}(\xi)) = \sum_j b_j \tilde{\rho}_{\ast}(\xi)^j.    
  \end{equation}
  Conjugating $\tilde{\rho}_{\ast}(\xi)$ to a diagonal matrix and
  inserting in the above equation gives then the desired result.
\end{proof}

\begin{theorem}
  \label{thm:approxpoisson}
  Let $(\rho_m:K \to \tmop{U}(V_m))_{m\in \mathbb{N}}$ be a sequence
  of irreducible representations with strictly increasing dimension.
  Assume that $\xi\in\mathcal{H}$ has the following properties:
  \begin{enumerate}
  \item Every eigenvalue of $\tilde{\rho}_{\ast,m}(\xi)$ has multiplicity one.
  \item $S:=\bigcup_{m\in \mathbb{N}}
    \tmop{Spec}(\tilde{\rho}_{\ast,m}(\xi))$ is a discrete subset of $\mathbb{R}$.
  \end{enumerate}
  Then for every absolutely continuous measure $\mu$ on $\mathbb{R}^+$
  with $\int_0^{\infty}xd\mu\in[0,1]$ there exists a function $f\in
  \tmop{Hol}(\mathbb{C})$ and a subsequence $(\rho_{m_k}:K \to
  \tmop{U}(V_{m_k}))_{k\in \mathbb{N}}$ such that $\eta:=f(\xi)$
  satisfies
  \begin{equation}
    d_{KS} (\mu_{\tilde{\rho}_{\ast,m_k}(\eta)},\mu) \to 0 \text{ as } k\to\infty.
  \end{equation}
\end{theorem}

\begin{proof}
  We begin by choosing a subsequence $\rho_{m_k}$ in the following
  way. First, we set $r_{m_1}=\rho_1$ and proceed inductively by
  requiring that
  \begin{equation}
    \label{eq:necessarygroth_poisson}
    N_{k+1}:=\tmop{dim}\,\rho_{m_k+1} \geq k (\tmop{dim}\,\rho_{m_k} + 2).
  \end{equation}
  Without loss of generality we assume that $N_1\geq 3$ and find an
  $N_1$-tuple $X_1$ such that
  \begin{equation}
    d_{KS}(\mu(X_1),\mu) \leq \frac{2}{N_1}.
  \end{equation}

  We now proceed inductively again, i.e. by Corollary
  \ref{cor:approxpoisson} in the Appendix, there is an $N_{k+1}$-tuple
  $X_{k+1}$ that contains the $N_{k}$-tuple $X_k$ as subset such that
  \begin{equation}
    d_{KS}(\mu(X_{k+1}),\mu) \leq \frac{N_k+2}{N_{k+1}} \leq \frac{1}{k},
  \end{equation}
  where the last inequality follows from
  (\ref{eq:necessarygroth_poisson}).

  Elementary complex analysis yields that there exists a holomorphic
  function $f:\mathbb{C}\to \mathbb{C}$, such that
  \begin{equation}
     f(\,\tmop{Spec} (\tilde{\rho}_{\ast,m_k}(\xi))\, ) = 
     X_k\ \forall k\in\mathbb{N}, 
  \end{equation}
  since $S$ is a discrete subset in $\mathbb{R}$ and each $X_k\subset
  X_{k+1}$.  By (\ref{eq:specunderanalyticfct}) it follows that
  \begin{equation}
    \tmop{Spec} (\tilde{\rho}_{\ast,m_k}(\tilde{f}(\xi))) =X_k 
    \ \forall k\in\mathbb{N}.
  \end{equation}
  Thus, $\eta=\tilde{f}(\xi)$ has the property
  \begin{equation}
    d_{KS} (\mu_{\tilde{\rho}_{\ast,m_k}(\eta)},\mu) \to 0 \text{ as } k\to\infty.
  \end{equation}
\end{proof}

Operators $\xi$ with the above properties will in general exist for
every ray of irreducible representations. One strategy of producing them
goes as follows:

Start with an operator $\xi$ of degree 2 that fulfills condition 1.
Such operators can be found for every simple group $K$ and should
exist in general. We now force condition 2 to hold by adding Casimir
operators to $\xi$.  Recall that Casimir operators act by scalar
multiplication so they just add these scalars to the eigenvalues.  If
these scalars increase quickly enough, the spectra of $\xi$ along the
irreducible representations will lie in disjoint intervals and
consequently condition 2 is satisfied.

The problem is that the operator $\xi$ depends on the group $K$ and we
do not know if there is an abstract way of giving examples. So we will
give here an example for $K=SU_n$ for the ray of irreducible
representations through the standard representation.

\begin{proposition}
  Let $(\rho_m: SU_n \to \tmop{U}(V_m))_{m\in\mathbb{N}}$ be the
  sequence of irreducible representations on the homogeneous
  polynomials of degree $m$ in $\mathbb{C}[x_1,\ldots,x_n]$.

  Then there exists an operator $\xi\in \mathcal{U}(\mathfrak{g})$ that
  satisfies the conditions of Theorem \ref{thm:approxpoisson}.
\end{proposition}
   
\begin{proof}
  Let $\alpha_j$ denote the $n\times n$-matrix with $1$ in the $j$-th
  diagonal component and $-1$ in the $(j+1)$-th diagonal
  component. Every other component should be equal to zero. These
  matrices form a basis for the standard maximal torus in
  $SL_n(\mathbb{C})=SU_n^{\mathbb{C}}$. They also define a system of
  simple roots (cf.\ the tables in Appendix C of \cite{knapp}). 

  The operation of $\alpha_j$ on the homogeneous polynomial
  $x_1^{a_1}\ldots x_n^{a_n}$ of degree $m$ is given by
  \begin{equation}
    \rho_{\ast,n}(\alpha_j) . x_1^{a_1}\ldots x_n^{a_n} = (a_j - a_{j+1}) x_1^{a_1}
    \ldots x_n^{a_n}.
  \end{equation}
  Therefore, the largest eigenvalue of $\rho_{\ast,n}{\alpha_j}$ is
  $m$ and the smallest $-m$ and every other eigenvalue is an integer
  number in-between these extremes. Now, we consider the operator $\xi =
  \sum c_j \alpha_j$, where the $c_j$ are real constants with $0 < c_j
  < \frac{1}{n}$ and which are linearly independent over $\mathbb{Q}$.
  Thus, $\rho_{\ast,m}(\xi)$ is represented as diagonal matrix and has
  eigenvalues with multiplicity greater than 1, since otherwise there
  would exist a linear relation between the $c_j$ over $\mathbb{Q}$.
  Note that by the choice of the $c_j$ the eigenvalues of $\xi$ are
  still in the interval $[-m,m]$.

  By now $\xi$ satisfies condition 1 of Theorem
  \ref{thm:approxpoisson} and we will now add the Laplace operator to
  $\xi$ to guarantee that condition 2 holds. For this, let $\Omega\in
  \mathcal{U}(\mathfrak{sl}_n(\mathbb{C}))$ be the Laplace operator
  associated to $\mathfrak{sl}_n(\mathbb{C})$. It acts on the
  homogeneous polynomials of degree $m$ by $r_{\Omega,m}:=\langle
  m\lambda, m\lambda + 2\delta \rangle_{\tmop{Kil}}$, where $\lambda$
  is the highest weight of the standard representation of $SU_k$,
  $\delta$ denotes the half sum of positive weights and
  $\langle\cdot,\cdot \rangle_{\tmop{Kil}}$ denotes the Killing
  form. It follows that
  \begin{equation}
    r_{\Omega,m+1}-r_{\Omega,m} = \langle \lambda,\lambda+ 2 \delta 
    \rangle_{\tmop{Kil}} + m\langle \lambda,\lambda\rangle_{\tmop{Kil}} 
    + m\langle \lambda,\lambda+2\delta \rangle_{\tmop{Kil}}.   
  \end{equation}
  Choosing a constant $b$ such that
  \begin{equation}
    b(r_{\Omega,m+1}-r_{\Omega,m}) \geq 2m \ \forall\ m\in\mathbb{N}^\ast
  \end{equation}
  yields that
  \begin{equation}
    \xi':= b\Omega + \sum_j c_j\alpha_j
  \end{equation}
  fulfills conditions 1 and 2 of Theorem \ref{thm:approxpoisson}.
\end{proof}

\section{Rational Independence of the Spectra in Representations}

In this section we give a notion of generic operators in
$\mathcal{H}$.

\begin{definition}
  An abstract hermitian operator $\alpha \in \mathcal{H}$ is called
  \textbf{generic} if for every irreducible representation $\rho$ the
  eigenvalues of $\tilde{\rho}$ \ are linearly independent over
  $\mathbb{Q}$. We denote the set of generic operators in
  $\mathcal{H}$ by $\mathcal{H}_{\tmop{gen}}$.
\end{definition}

We start with the following theorem. 
\begin{theorem}
  \label{nowheredense-theorem}The set of generic operators
  $\mathcal{H}_{\tmop{gen}}$ is dense in $\mathcal{H}$.
\end{theorem}

Before the prove is given, we need to fix the notation. 
The ordered tuple of eigenvalues with multiplicity of a hermitian
matrix $A$ will be denoted by $X(A)$ and the set of ordered $n$-tuples by
$\mathbb{R}_{\tmop{ord}}^n$.

\begin{lemma}
  \label{nowheredense-label}
  Let $V$ be a unitary vector space of dimension $n$ and $\tmop{Herm}
  (V)$ be the real subspace of hermitian endomorphisms of $V$. 
  For every $\lambda\in (\mathbb{Q}^n)^{\ast}$ the set 
  \begin{equation}
    S_\lambda:=\{A \in \tmop{Herm} (V) : \lambda(X(A))=0   \}
  \end{equation}
  is nowhere dense in $\tmop{Herm} (V)$.
\end{lemma}

\begin{proof}
  Let $\lambda \in (\mathbb{Q}^n)^{\ast}$ be a non-zero linear
  form. The set $\lambda^{- 1} (0)$ is a hyperplane in $\mathbb{R}^n$,
  thus nowhere dense.  In follows that the intersection of
  $\mathbb{R}_{\tmop{ord}}^n\cap\lambda^{-1}(0)$ is nowhere dense in
  $\mathbb{R}_{\tmop{ord}}$.
  
  Now, let us fix a given point $x \in \mathbb{R}^n_{\tmop{ord}}$.
  From linear algebra we know that the set of hermitian operators with
  spectrum $\{x_1, \ldots, x_n \}$ is just the $U (n)$ orbit under matrix
  conjugation through the diagonal matrix
  \begin{equation} 
    X = \tmop{diag} (x_1, \ldots, x_n) . 
  \end{equation}
  Therefore, the set $\mathbb{R}_{\tmop{ord}}^n$ can be identified
  with $\tmop{Herm} (V) / U(n)$ and the projection map $p:\tmop{Herm}
  (V) \to \tmop{Herm} (V)/U(n) = \mathbb{R}_{\tmop{ord}}^n$ is an open
  map.

  Because preimages of nowhere dense sets under open maps are nowhere
  dense, the lemma is proved. 
\end{proof}

\begin{proof}
  (Theorem \ref{nowheredense-theorem})\ \  
  Since $\rho_\ast : \mathcal{U} (\mathfrak{g})
  \rightarrow \tmop{End} (V)$ is an irreducible, finite-di\-men\-sio\-nal
  representation, the induced mapping $\tilde{\rho}_\ast : \mathcal{H}
  \rightarrow \tmop{Herm} (V)$ is a real linear, surjective mapping between
  Fr\'echet spaces. Therefore it is an open mapping by the open mapping
  theorem.
  
  So for any given non-zero linear from $\lambda \in (\mathbb{Q}^{\dim
    V})^{\ast}$, the set
  \begin{equation}
    M_{\lambda,\rho} := \{ \alpha\in\mathcal{H}\, :\, X( \tilde{\rho}_\ast(\alpha))
    \in \lambda^{-1}(0)   \}
  \end{equation}
  is nowhere dense in $\mathcal{H}$. Otherwise, we could find an inner
  point in this set, but because $\tilde{\rho}$ is an open mapping
  this would contradict Lemma \ref{nowheredense-label}.
  
  Thus, the set
  \begin{equation}
    M := \bigcup_{\rho \tmop{irrep.}\,, \lambda\in (\mathbb{Q}^{\dim
    V})^{\ast}} M_{\lambda,\rho} 
  \end{equation}
  contains no inner point by Baire's category theorem, i.e. its
  complement is dense. It follows that $\mathcal{H}_{\tmop{gen}}$ \ is
  dense.
\end{proof}

\section[Ergodic Properties of $\mathcal{H}_{\tmop{gen}}$]{Ergodic
  Properties of \boldmath$\mathcal{H}_{\tmop{gen}}$\unboldmath}

Before we come to the main point of this section, we have to recall some
terminology from ergodic theory. All details can be found in
\cite{sinai} or \cite{cornfeldfominsinai}. We follow the latter in
terminology.

Let $(X, \mu)$ be a measure space, where $\mu$ denotes the measure on
some $\sigma$-algebra in the power set of $X$. A measurable map $f : X
\rightarrow X$ is called an automorphism of the measure space $(X,
\mu)$, if $f$ is bijective, $f^{- 1}$ is measurable again, and for all
measurable sets $A \subset X$, we have
\begin{equation} 
  \mu (f (A)) = \mu (f^{- 1} (A)) = \mu(A) . 
\end{equation}
By a flow $(\varphi_t)_{t \in \mathbb{R}}$ of the measure space $(X, \mu)$,
we mean a 1-parameter group of automorphisms of $(X, \mu)$, i.e., a group
homomorphism of $\mathbb{R}$ into the group of all automorphisms of the
measure space $(X, \mu)$ such that $\varphi : \mathbb{R} \times X
\rightarrow X$ is measurable.

For us $X$ will be an $N$-dimensional torus, i.e., $X = [0, 1]^N
\tmop{mod} 1$ and the measure $\mu$ is the Haar measure on $X$, which
is equal to the Lebesgue measure here. We consider some $N$-tuple $x =
(x_1, .., x_n)$ such that $0 < x_i < 1$ for all $i \in \{1, \ldots,
N\}$ \ and the $x_i$'s are linearly independent over the rational
numbers. The map $\varphi_t : X \rightarrow X, z \mapsto z + t \cdot x
\tmop{mod} 1$ defines a group homomorphism $\mathbb{R} \rightarrow
\tmop{Diff} (X)$, $t \mapsto \varphi_t$, where $\tmop{Diff} (X)$
denotes the group of diffeomorphisms of $X$. It is a standard fact
from ergodic theory that $(\varphi_t)_{t \in \mathbb{R}}$ is a flow
of the measure space $(X, \mu)$ (cf.\ \cite{cornfeldfominsinai})
Chapter 3, \S 1, Theorem 1).

A flow is called ergodic if for every $t \neq 0$, the only invariant
sets of $\varphi_t$ have measure either $0$ or $1$. We make use of the
following

\begin{theorem}
  {\dueto{Birkhoff}}Let $(X, \mu)$ be a measure space with $\mu (X) =
  1$ and $(\varphi_t)_{t \in \mathbb{R}}$ be a flow of the measure
  space $(X, \mu)$.  Then for every integrable function $f : X
  \rightarrow \mathbb{R}$,
  \begin{equation} 
     \bar{f} (y) := \lim_{t \rightarrow \infty} \frac{1}{2 t} \int_{- t}^t f
     (\varphi_{\tau} (y)) d \tau  = \int_X f (x) dx 
  \end{equation}
  for almost all $y \in T$ with respect to $\mu$. 
\end{theorem}

It is a standard result of ergodic theory that $(\varphi_t)_{t \in
\mathbb{R}}$ is a uniquely ergodic flow, i.e., $\bar{f}$ is
constant, (cf.\ \cite{cornfeldfominsinai}) Chapter 3, \S 1, 
Theorem 2).

In this case, we obtain the formula for the
characteristic function $\chi_A$ of a measurable set $A$:
\begin{equation}
  \label{eq:ergodic-formula}
  \lim_{t \rightarrow \infty} \frac{1}{2 t} \int_{- t}^t \chi_A
  (\varphi_{\tau} (y)) d \tau  = \mu (A) \ \forall y\in X. 
\end{equation}
Let us now consider an element $\alpha \in \mathcal{H}_{\tmop{gen}}$
and the induced irreducible, finite-dimen\-sio\-nal representation
$\rho_\ast : \mathcal{U} (\mathfrak{g}) \rightarrow \tmop{End}
(V)$. Since $\tilde{\rho}_\ast (\alpha)$ is a self-adjoint operator,
it follows that $(\exp (2 \pi i \tilde{\rho} (\alpha) t))_{t \in
  \mathbb{R}}$ is a uniquely ergodic flow on the torus
\begin{equation}
  T (V) = \tmop{closure}
(\{\exp (2 \pi i \tilde{\rho}_{\ast} (\alpha) t) | t \in \mathbb{R}\}).
\end{equation}

This torus depends on the starting direction
$\tilde{\rho}_{\ast}(\alpha)$, but we will in the following always
assume that we have conjugated it into a diagonal matrix. There is no
loss of generality because we are only interested in the eigenvalues
and they do not change under conjugation. Thus, we will just write
$T_N$ for the $N$-dimensional torus, i.e.,
\begin{equation}
  T_N = \{ \tmop{diag}(e^{2\pi i\phi_1},\dots, e^{2\pi i \phi_N})  : \phi_j \in [0,1] \}.
\end{equation}

\section[The Sets $B_N$]{The Sets  \boldmath$B_N$\unboldmath}

In this section we will use the ergodic properties of
$\mathcal{H}_{\tmop{gen}}$ in combination with a theorem of Chapter
\ref{chap:poisson-spectral-stat} to connect the spectral properties
of an abstract hermitian operator with the Poisson-statistics. For
this we first need to fix some notation.

For a unitary automorphism $A \in U (V)$ of a finite-dimensional
unitary vector space $V$ of dimension $N$ we have the nearest neighbor
statistics $\mu_c(X(A))$ as defined in Definition
\ref{def:nearest-neighbor-unitary-matrix} of the Appendix. By
$\mu_{\tmop{Poisson}}$ we denote the absolutely continuous probability
measure on the positive real line with density function $\exp (- x)$
with respect to the Lebesgue measure. Finally, let us write
$d_{\tmop{KS}} (\mu_1, \mu_2)$ for the Kolmogorov-Smirnoff distance
(cf.\ (\ref{eq:def_dks}) in the Appendix).

The following theorem is analogous to the second main theorem of
{\cite{katzsarnak}} and is the main result of Chapter
\ref{chap:poisson-spectral-stat}.

\begin{theorem}
  \label{katz-sarnak-theorem}
  Let $\alpha > 0$ be given. Then there exists an natural number
  $N_0$ such that for every $N \geq N_0$
  \begin{equation}
    \int_{T_N} d_{KS} ( \mu_c(X(A)),\mu_{\tmop{Poisson}}  ) dA 
    < \frac{1}{e^{\alpha\sqrt{\log N}}}.
  \end{equation}
\end{theorem}

The rather technical proof is given in Chapter
\ref{chap:poisson-spectral-stat}, cf.~Theorem
\ref{thm:main-theorem-tn}.

\begin{corollary}
  For all $\alpha \in \mathbb{R}$
  with $\alpha > 0$ and any $N \geq N_0= N_0(\alpha)$ we have
  \begin{equation} 
     d_{\tmop{KS}} (\mu_A, \mu_{\tmop{Poisson}}) \leqslant e^{-\frac{1}{2}\alpha\sqrt{\log(N)}}
  \end{equation}
  for all $A$ in a set in $T_N$ of measure at least $1 - e^{-\frac{1}{2}\alpha\sqrt{\log(N)}}$.
\end{corollary}

\begin{proof}
  Let us assume the contrary, i.e., we assume that
  \begin{equation}
    d_{\tmop{KS}}
    (\mu_A, \mu_{\tmop{Poisson}}) >
    e^{-\frac{1}{2}\alpha\sqrt{\log(N)}}       
  \end{equation}
  on a set M of measure at least
  $e^{-\frac{1}{2}\alpha\sqrt{\log(N)}}$. Then
  \begin{equation} 
     \int_M d_{\tmop{KS}} (\mu_A, \mu_{\tmop{Poisson}}) d \tmop{Haar} (A)
     > e^{-\frac{1}{2}\alpha\sqrt{\log(N)}} e^{-\frac{1}{2}\alpha\sqrt{\log(N)}}  
     = e^{-\alpha\sqrt{\log(N)}}. 
  \end{equation}
  Since the integrand is always positive, this is a contradiction to
  Theorem \ref{katz-sarnak-theorem}.
\end{proof}
This motivates the following definition.
\begin{definition}
  Let $\alpha > 0$ be given. The set $B_N$ is given by
  \footnote{The letter $B$ in $B_N$ is not an abbreviation for
    big. In fact these sets are small.}
  \begin{equation} B_N := \left\{B \in T_N : d_{\tmop{KS}} (\mu_B, \mu_{\tmop{Poisson}})
     \geq e^{-\frac{1}{2}\alpha\sqrt{\log(N)}} \right\}. 
  \end{equation}
\end{definition}%
\nomenclature{$B_N$}{The set of unitary $N\times N$-matrices with nearest neighbor
  statistics not close to $\mu_{\tmop{Poisson}}$}%

It is clear that $B_N$ depends on the choice of $\alpha$. However, for
reasons of simplicity we suppress this fact in the notation.  In the
following we will always assume that the $N$ are so large that Theorem
\ref{katz-sarnak-theorem} is valid, i.e.\ $N\geq N_0\geq 2$.

Let us now collect some properties of $B_N$. First of all, $B_N$ is not
empty because the identity matrix $E_N$ is in $B_N$. For this just recall
that $\int_0^c\mu_{\tmop{Poisson}}$ is close to zero for small $c$ and
that $\int_0^c d \mu_{E_N} = 1$ for every non-negative $c$, so
$d_{\tmop{KS}} (\mu_{E_n}, \mu_{\tmop{Poisson}}) = 1$.

Due to the fact that the map $A \mapsto d_{\tmop{KS}} (\mu_A,
\mu_{\tmop{Poisson}})$ is continuous (cf. Lemma
\ref{lem:dks-continous}), $B_N$ is closed and the identity matrix is an
inner point as a consequence of continuity.

Moreover, $B_N$ is invariant under scalar multiplication with $z =
e^{i \lambda}$, where $\lambda \in \mathbb{R}$, cf.\ Chapter
\ref{chap:gen-level-spacings}.

\begin{figure}
  \centering
  \includegraphics[width=8cm]{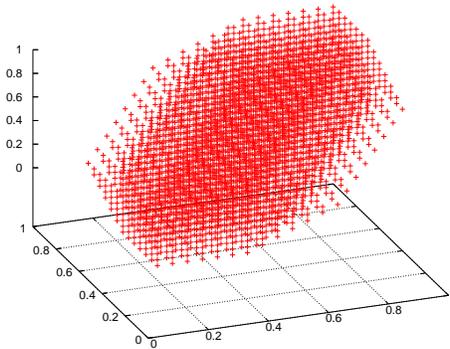}
  \caption{A picture of $B_3$.}
  \label{fig:bn3d}
\end{figure}

\begin{figure}
  \centering
  \includegraphics[width=6cm]{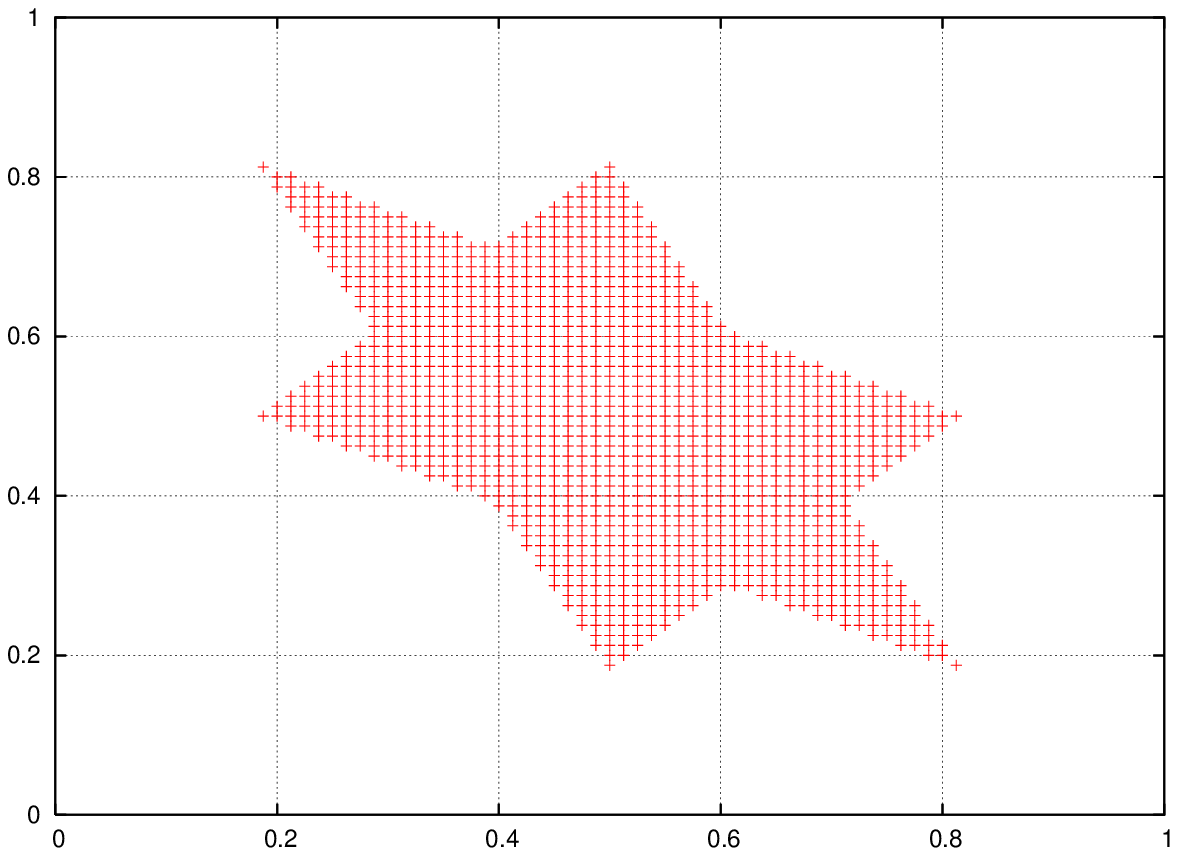}
  \includegraphics[width=6cm]{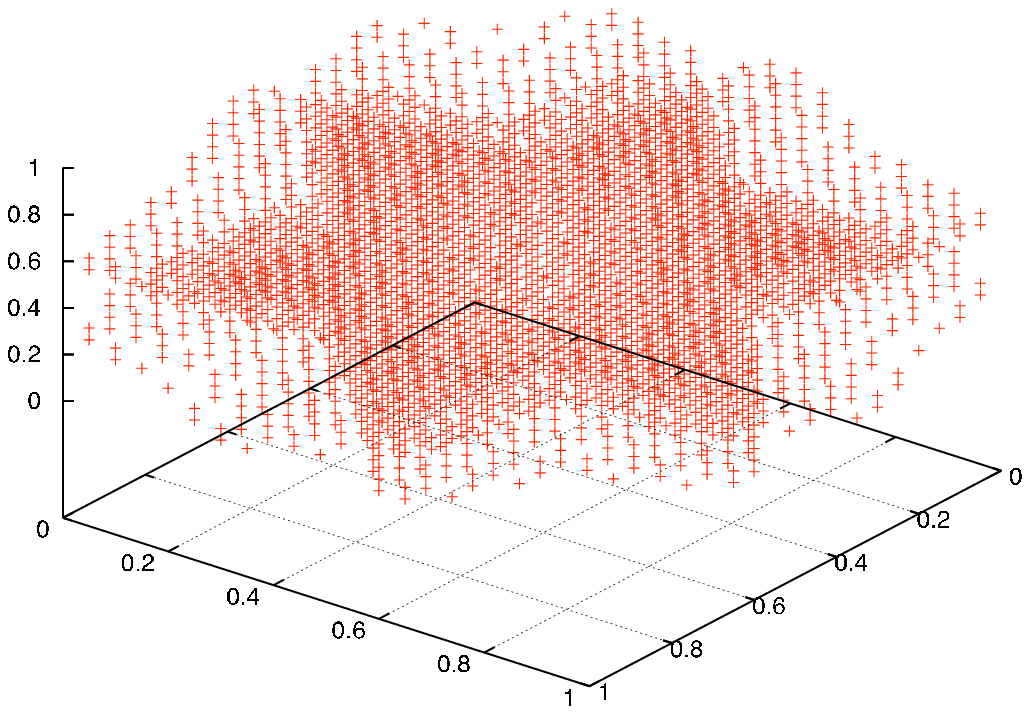}
  \caption{Pictures of $B_3$ and $B_4$ intersected with the hyperplane normal to the diagonal.}
  \label{fig:bn3d_schnitt}
\end{figure}

The set $B_3$ for $\alpha=\frac{4}{3}$ is visualized by Figure
\ref{fig:bn3d}. For the drawing, we have discretized the torus $T_3$
into a cubical lattice with $20\times 20 \times 20$ points and
calculated a discretized version of $d_{KS}$ for a grid size of $20$
points. The axises show the coordinates $\phi_1,\phi_2$ and $\phi_3$.
The intersection of $B_3$ and the cubical grid is the drawn set of
points. The definition of the discretized version is given as
Definition \ref{def:d_ks_mgrid}. One can see the invariance under
multiplication with $e^{i\phi}$ here as invariance under diagonal
shifts.

Thus, it is enough to know the sets $B_N$ only on that hyperplane
which is normal to the diagonal and contains the point
$\frac{1}{2}(1,\dots,1)$, i.e., the hyperplane given by
\begin{equation}
   a_1\phi_1+ \dots + a_N\phi_N = N/2.
\end{equation}
Figure \ref{fig:bn3d_schnitt} shows these hyperplanes for $N=3,\,4$ parametrized by $\phi_1, \dots, \phi_{N-1}$.

We now use the ergodic properties of $\mathcal{H}_{\tmop{gen}}$ to formulate
our key lemma.

\begin{lemma}
  Let $\gamma \in \mathcal{H}_{\tmop{gen}}$ and $\rho : K
  \rightarrow \tmop{U} (V)$ be an irreducible, finite-dimensional, unitary
  representation with $\dim V = N$ and denote the characteristic function
  of the set $B_N$ by $\chi$. \ Then
  \begin{equation}
    \lim_{t \rightarrow \infty} \frac{1}{2 t} \int_{- t}^t \chi (\exp (2 \pi i
    \tilde{\rho}_\ast (\gamma) \tau))  d \tau = \tmop{vol}_{T_N} (B_N),
    \label{bn-equation}
  \end{equation}
  where $\tmop{vol}_{T_N} (B_N)$ denotes the measure of $B_N$ with respect to
  the Haar measure on $T_N$.
\end{lemma}

\begin{proof}
  This is just the ergodic property of equation (\ref{eq:ergodic-formula}).
\end{proof}

We would like to emphasize the role of $t$ in the above lemma. Consider
the set $R(N)$ defined by
\begin{equation}
  R (N) =\{t \in \mathbb{R}: \exp (2 \pi i
  \tilde{\rho}_\ast (\gamma) t \in B_N \}.
\end{equation}

\begin{corollary}
  Under the assumptions of the above lemma
  \begin{equation} d_{\tmop{KS}} (\mu_{\exp (2 \pi i \tilde{\rho}_\ast (\gamma) t},
     \mu_{\tmop{Poisson}}) < e^{-\frac{1}{2}\alpha\sqrt{log(N)}} 
  \end{equation}
  for every $t \not\in R (N)$.
  
  Moreover $\mathbb{R} \backslash R (N)$ has infinite measure and he have
  the following estimation on the size of $R (N)$
  \begin{equation}
    0 < \lim_{t \rightarrow \infty} \frac{1}{2 t} \int_{- t}^t \chi_{R (N)} (\tau)
    d \tau < e^{-\frac{1}{2}\alpha\sqrt{log(N)}}, \label{rn-equation}
  \end{equation}
  where $\chi_{R (N)}$ denotes the characteristic function of $R (N)$. 
\end{corollary}

\begin{proof}
  By virtue of equation (\ref{bn-equation}) we obtain the corollary.
\end{proof}

\section[Convergence to $\mu_{\tmop{Poisson}}$]{Convergence 
to \boldmath $\mu_{\tmop{Poisson}}$\unboldmath}

From now on, consider a sequence $(\rho_k)_{k \in \mathbb{N}}$ of
irreducible, unitary representations $\rho_k : K \rightarrow \tmop{U} (V_k)$ 
such that $d_k := \dim(V_k)$ is increasing. Before the main result can
be stated, it is necessary to introduce two rather technical conditions.

\begin{definition}
  A sequence $(\rho_k)_{k\in\mathbb{N}}$ is said to be of \textbf{admissible
  growth}, if there exists an $\alpha>0$ such that
\begin{equation}
  \sum_{k=0}^\infty e^{-\frac{1}{2}\alpha\sqrt{\log(d_k)}} < \infty.
\end{equation}
\end{definition}

\begin{definition}
  A generic hermitian operator $\gamma\in\mathcal{H}_{\tmop{gen}}$ is
  said to be \textbf{admissible of width} \boldmath
  $\epsilon$\unboldmath\ for the sequence $(\rho_k)_{k\in\mathbb{N}}$,
  where $0<\epsilon<1$ if there exists a $k_0$ and a $t_0$ such that
  for all $t\geq t_0$ and all $k\geq k_0$ the inequality
  \begin{equation}
   \label{eq:admissible_generic_operator}
     \left| \tmop{vol}(B_{d_k})-\frac{1}{2t}\int_{-t}^{t} \chi_{B_{d_k}}
     (\tmop{exp}(2\pi
     i\rho_{\ast,k}(\gamma)\tau) d\tau \
     \right| < \epsilon
  \end{equation}
  holds. Here $\chi_{B_{d_k}}$ denotes the characteristic function of
  the set $B_{d_k}$ as defined above.
\end{definition}

Let us briefly discuss these definitions. As will become clear in the
following theorem the first describes a condition on the growth of the
dimensions $d_k$. By a direct calculation we see that the condition
requires $d_k$ to grow faster than
$e^{\left(\frac{2\log(k)}{\alpha}\right)^2}$. We will come back to
this later.
 
The second definition guarantees that we are outside the sets $B_N$ in
each representation. For fixed $k$ the condition can be
fulfilled for every $\epsilon$ by Birkhoff's ergodic theorem. But we
require here that $t_0$ as a function of $k$ is bounded. So, condition
(\ref{eq:admissible_generic_operator}) only fails, if
\begin{equation}
  \left| \tmop{vol}(B_{d_k})-\frac{1}{2t}\int_{-t}^{t}
    \chi_{B_{d_k}}(\dots) d\tau \ \right| \to 1
\end{equation}
is true. This will happen if the leaving time, i.e., the supremum
of all $t$, such that $\tmop{exp}(2\pi i\rho_{\ast,k}(\gamma)\tau)\in
B_{d_k}$, converges too rapidly to infinity as function of $k$. In
Lemma \ref{lem:mainlemma_algebra} we saw this kind of behavior. 
The reader may wonder if operators of width
$\epsilon$ do exist at all. But in Section
\ref{sec:examples-convergence} we saw examples of operators $\gamma$
whose nearest neighbor statistics converge to a given measure
$\mu$. Although the situation is a little different here, because of
the exponentiation, we could use the proof of Theorem
\ref{thm:approxpoisson} to construct operators $\gamma$ such that
$\exp(2\pi i\tilde{\rho}_{\ast,k})$ has nearest neighbor statistics
which converge to $\mu_{\tmop{Poisson}}$. These $\gamma$ have a leaving
time less than 1 by construction.

Now we state our key theorem in this chapter.
\begin{theorem}
  \label{thm:maintheorem}Let $\gamma \in \mathcal{H}_{\tmop{gen}}$ be
  admissible of width $\epsilon$ for a sequence
  $(\rho_k:K\to\tmop{U}(V_k))_{k\in\mathbb{N}}$ of irreducible,
  unitary representations which is of admissible growth. Then for
  every $\epsilon'>0$ there exists a set $R=R(\epsilon')$ in
  $\mathbb{R}$, such that
  \begin{equation} 
      \label{eq:R_percentage}
      \lim_{r \rightarrow \infty} \frac{1}{2r}\int_{- r}^r \chi_R (x) d x \leq \epsilon+\epsilon'
  \end{equation}
  and 
  \begin{equation} \mu_{\exp (2 \pi i \tilde{\rho}_{\ast,k} (\gamma) t)} \rightarrow
     \mu_{\tmop{Poisson}} \text{ as } k \rightarrow \infty
  \end{equation}
  for all $t \not\in R$.
\end{theorem}

Before we prove the theorem, let us discuss the claim about the
measure of $R$.  Any bounded set $R$ is of this type, or any set of
measure 0. But from the point of view of percentage of real numbers,
we prove that a fraction of $(1-\epsilon-\epsilon')$ of the real numbers
yields convergence to $\mu_{\tmop{Poisson}}$ for the subsequence.

\begin{proof}
  According to the condition 
  of (\ref{eq:admissible_generic_operator}), we find a $t_0$ such that
  for all $t\geq t_0$ and all $k\geq k_0$ 
  \begin{equation}
    \left| \tmop{vol}(B_{d_k})-\frac{1}{2t}\int_{-t}^{t} \chi_{B_{d_k}}(\tmop{exp}(2\pi
   i\rho_{\ast,k}(\gamma)\tau) d\tau \
  \right| < \epsilon.
  \end{equation}
  By the definition of admissible growth it follows that
  \begin{equation} 
      \sum_{k = 1}^{\infty}  e^{-\frac{1}{2}\alpha\sqrt{log(d_k)}} < \infty. 
  \end{equation}
  Thus for every $\epsilon_1 > 0$ we find a natural number $N_0 = N_0
   (\epsilon_1)$ such that
  \begin{equation}
    \sum_{k = N_0}^{\infty} e^{-\frac{1}{2}\alpha\sqrt{log(d_k)}} < \epsilon_1 . \label{dn-sum}
  \end{equation}
  Now set
  \begin{equation} 
    R_{\epsilon_1} = \bigcup_{k = N_0}^{\infty} R (d_k),
  \end{equation}
  where $R(d_k) = \{t \in \mathbb{R}: \exp (2 \pi i \tilde{\rho} (\gamma) t
  \in B_{d_k} \}$. 
  We set $Q_{\epsilon_1} = \mathbb{R} \backslash R_{\epsilon_1}$ and note
 that for all $t\in Q_{\epsilon_1}$
  \begin{equation}
    \mu_{\exp (2 \pi i \tilde{\rho}_{\ast,k} (\gamma) t)} \rightarrow
    \mu_{\tmop{Poisson}} \text{ as } k \rightarrow \infty.
  \end{equation}
  Now we have to show that $Q_{\epsilon_1} \neq \emptyset$.

  By enlarging $N_0$ if necessary, we may also assume that $k_0\leq N_0$. We 
  fix an interval $[-t,t]$, where $t\geq t_0$, and obtain
  \begin{equation}
     2t(\tmop{vol}(B_{d_k})-\epsilon) \leq \tmop{vol}(R(d_k) \cap [-t,t]) \leq 2t(\tmop{vol}(B_{d_k})+\epsilon).
  \end{equation}
  Summing over all $k\geq N_0$ and applying (\ref{dn-sum}) it follows that
  \begin{equation}
    \label{eq:R_percentage_estimate}
    \tmop{vol}(R_{\epsilon_1}\cap [-t,t])=\tmop{vol}(\bigcup R(d_k) \cap [-t,t])  \leq   
    2t(\epsilon_1+\epsilon).
  \end{equation}
  We can choose $\epsilon_1$ so small that
  $\epsilon_1+\epsilon<1$. This yields
  \begin{equation}
   Q_{\epsilon_1}\cap[-t,t] \neq \varnothing.
  \end{equation}
  It remains to show (\ref{eq:R_percentage}). But since
  (\ref{eq:R_percentage_estimate}) holds for all $t\geq t_0$:
  \begin{equation}
    \tmop{vol}(R_{\epsilon_1}\cap [-t,t])= \int_{-t}^t \chi_{R_{\epsilon_1}}(s) ds \leq 2t(\epsilon_1+\epsilon).
  \end{equation}
  This completes the proof of the theorem.
\end{proof}

Let us briefly discuss this theorem. For every generic, admissible
operator one has convergence of the nearest neighbor distributions for
all $t\not\in R$. But the reader may wonder how restrictive the
condition of admissible width is. This will depend on the geometric
structure of the sets $B_N$. If they are regular enough, the condition
of admissible width should be automatically fulfilled for most generic
operators. Unfortunately, we do not know enough about this structure
yet, although in low dimensions the sets $B_N$ are very regular 
(cf.~Figures \ref{fig:bn3d} and \ref{fig:bn3d_schnitt}).

%%% Local Variables:  
%%% mode: latex 
%%% TeX-master: "../thesis" 
%%% End:

% LocalWords:  diffeomorphisms Haar Birkhoff Kolmogorov Birkhoff's Weyl's Borel
% LocalWords:  Weil Weyl Baire's complexification

% \chapter{The Poisson spectral statistics for tori}
\chapter{The Poisson Spectral Statistics for Tori}
\label{chap:poisson-spectral-stat}
In this chapter we give a proof for the convergence of the
nearest neighbor statistics of a real torus $T(N)$ to the Poisson spectral
statistics in the sense of the Kolmogorov-Smirnov distance.

We follow the structure of the proof in \cite{katzsarnak} for the CUE
case but will try to make this chapter as self-contained as possible,
citing only some combinatorial lemmas and some facts about measures.

\section{Some Combinatorics}
We give here the basic definitions of $\tmop{Sep}$, $\tmop{Cor}$,
$\tmop{Clump}$ and so on from \cite{katzsarnak} again. To do this let
$f : \mathbb{R} \rightarrow \mathbb{R}$ be a function, $a$ be a
non-negative integer called the \emph{separation} and $X$ be an
$N$-tuple of real numbers in increasing order.

We define
\begin{equation} \tmop{Clump} (a, f, N, X) = \sum_{1 \leq t_1
    \leq \ldots \leq t_{a + 2} \leq N} f (x_{t_{a + 2}}
  - x_{t_1}) 
\end{equation}
and
\begin{equation} \tmop{Sep} (a, f, N, X) = \sum_{1 \leq t_1
    \leq \ldots \leq t_{a + 2} \leq N, t_{j + 1} - t_j
    = 1 \text{ for all } j} f (x_{t_{a + 2}} - x_{t_1})
  . 
\end{equation}

Let us briefly discuss what these definitions signify, first taking a
closer look at $\tmop{Clump}$. We sum over all $(a+2)$-tuples 
$(t_1,\dots,t_{a+2})$ with increasing entries such that the last
entry is smaller or equal than $N$, thereby evaluating the function $f$
at the differences between $x_{t_{a+2}}-x_{t_1}$. If $a+2>N$ then there 
are no tuples to sum over, so $\tmop{Sep}$ and $\tmop{Clump}$ vanish 
identically.

Formally we can think of this as integrating the function $f$ over a
sum of Dirac measures at the points $x_{t_{a+2}}-x_{t_1}$. The same
applies to the function $\tmop{Sep}$ with the restriction that we sum
only about the $(a+2)$-tuples of the form $(t_1,t_1+1,\dots,t_1+a+1)$. 

If we consider $a=0$, then we evaluate $f$ exactly at the nearest
neighbor spacings. This may give a clear motivation why we are
interested in $\tmop{Sep}$. The point in the definition of
$\tmop{Clump}$ will become clear later on. For the moment, let us
just indicate that there will be a combinatorial identity expressing
$\tmop{Sep}$ as alternating sum over some versions of $\tmop{Clump}$.    
 
By now, $\tmop{Sep}$ and $\tmop{Clump}$ are defined over increasing
$N$-tuples $X$. We extend this definition to all $N$-tuples by first
ordering the tuple $X$.

\begin{equation} \tmop{Clump} (a, f, N, \cdot) : \mathbb{R}^N
  \rightarrow \mathbb{R}, X \rightarrow \tmop{Clump} (a, f, N, X
  \tmop{ordered}) 
\end{equation} 
and
\begin{equation} \tmop{Sep} (a, f, N, \cdot) : \mathbb{R}^N
  \rightarrow \mathbb{R}, X \mapsto \tmop{Sep} (a, f, N, X
  \tmop{ordered}) . 
\end{equation} 

$\tmop{Sep}$ and $\tmop{Clump}$ are special cases of a certain 
class of functions which we will deal with in the following. We define 
this class in the following way:

\begin{definition}
  Let $N \geq 2$ be an integer.
  
  A function $f : \mathbb{R}^N \rightarrow \mathbb{R}$ is called a
  {\textbf{function of class}} $\mathbf{\mathcal{T}(N)}$ if $f$ is
  Borel measurable, $S_N$-invariant and invariant under additive
  diagonal translations 
  \begin{equation} 
    (x_1, \ldots, x_N) \mapsto (x_1 + t, \ldots, x_N + t),
  \end{equation}
  with $t \in \mathbb{R}$.

  A function $f : \mathbb{R}^N \rightarrow \mathbb{R}$ is called a
  {\textbf{function of class}} $\mathbf{\mathcal{T}_0 (N)}$ if f is a
  function of class $\mathcal{T}(N)$ and f vanishes outside the set
  $\{(x_1, \ldots, x_N) \in \mathbb{R}^N : \max_{i, j} |x_i - x_j |
  \leq \alpha\}$ for some $\alpha>0$. We abbreviate this condition by
  \begin{equation} 
      \tmop{supp} f \leq \alpha . 
  \end{equation}
\end{definition}

The following lemma lists some basic properties of $\tmop{Sep}$ and
$\tmop{Clump}$. This is Lemma 2.5.11 of \cite{katzsarnak}.

\begin{lemma}
  \label{lem:properties-of-sep-and-clump}For $a \in \mathbb{N}$ and
  $f : \mathbb{R} \rightarrow \mathbb{R}$ Borel measurable and $N
  \geq 2$.
  
  Then $\tmop{Sep} (a, f, N, \cdot)$ and $\tmop{Clump} (a, f, N,
  \cdot)$ are functions of class $\mathcal{T}(N)$. If $f$ is
  continuous, then $\tmop{Sep} (a, f, N, \cdot)$ and $\tmop{Clump} (a,
  f, N, \cdot)$ are also continuous.
  
  If $f$ vanishes outside the interval $[- \alpha, \alpha]$,
  $\tmop{Sep} (a, f, N, \cdot)$ and $\tmop{Clump} (a, f, N, \cdot)$
  are of class $\mathcal{T}_0 (N)$ and
  \begin{equation}
    \tmop{supp} \tmop{Sep} (a, f, N, \cdot) \leq \alpha \text{ and }
    \tmop{supp} \tmop{Clump} (a, f, N, \cdot) \leq \alpha . 
  \end{equation}
\end{lemma}

\begin{proof}
  See \cite{katzsarnak} p.52.
\end{proof}

Using $\tmop{Clump}$ we define a third function for an integer $k$, $k
\geq a$.
\begin{equation} \tmop{TClump} (k, a, f, N, \cdot) : \mathbb{R}^N
  \rightarrow \mathbb{R}, X \mapsto \binom{k}{a} \tmop{Clump} (k, f,
  N, X) . 
\end{equation} 
Note that this definition may seem a bit superfluous, but it is added
here to show the parallels to \cite{katzsarnak}. If we were working
with multiple neighbor statistics, i.e. $r > 1$ in terms of
\cite{katzsarnak}, then $\tmop{TClump}$ would be a more complicated
sum.

We now relate this functions on $\mathbb{R}^N$ to functions on the
torus $T (N)$. Again following \cite{katzsarnak}, we name these functions
$\tmop{Int}$ for ``integral'', Cor for ``correlation'' and
$\tmop{TCor}$ for ``total correlation''.

These are defined as functions from $T (N)$ to $\mathbb{R}$ which map
$A\in T(N)$ as follows
\begin{eqnarray*}
  \tmop{Int} (a, f, T (N), A) & := & \frac{1}{N} \tmop{Sep} \left( a, f,
    N, \frac{N}{2 \pi} X (A) \right)\\
  \tmop{Cor} (a, f, T (N), A) & := & \frac{1}{N} \tmop{Clump} \left( a,
    f, N, \frac{N}{2 \pi} X (A) \right)\\
  \tmop{TCor} (k, a, f, T (N), A) & := & \frac{1}{N} \tmop{TClump} \left(
    k, a, f, N, \frac{N}{2 \pi} X (A) \right),
\end{eqnarray*}
where $X (A)$ is  $-i$ times the component-wise logarithm of $A$, i.e. for the matrix $A =
\tmop{diag} (e^{i \varphi_1}, \ldots, e^{i \varphi_N})$ with $0
\leq \varphi_j < 2 \pi$ for all $j$, we have $X (A) = (\varphi_1,
\ldots, \varphi_N)$. It is now, obvious why we study these
objects because
\begin{equation} \tmop{Int} (a, f, T (N), A) = \int_{\mathbb{R}} f d
  \mu (\tmop{naive}, A, T (N), a). 
\end{equation} 
It is exactly this $\mu (\tmop{naive}, A, T (N),a)$ we want to study
for a=0. For $a \geq 1$ we may take the above equation as definition
of $\mu (\tmop{naive}, A, T (N),a)$. In the notation of Chapter 1 this
measure is given as
\begin{equation}
  \mu (\tmop{naive}, A, T (N), 0)  = \frac{1}{N}\int_A
  \sum_{j=1}^{N-1} \delta\left(y-\frac{N}{2\pi}\cdot(\varphi_{j+1}-\varphi_j) \right) dy
\end{equation}
if $a=0$, which is almost identical to $\mu_c(X)(A)$ but the wrapped
eigenangle between $x_N$ and $x_1$ is missing. Therefore it is called
``naive'' in \cite{katzsarnak}.

If we think of $\tmop{Int}$, Cor and $\tmop{TCor}$ as random
variables, we may calculate their expectation value. But instead of
writing $E (\tmop{Int} (a, f, T (N), A))$ we use capital letters:
\begin{eqnarray*}
  \tmop{INT} (a, f, T (N), A) & := & \int_{T (N)} \tmop{Int} (a, f, T
  (N), A) d A, \\
  \tmop{COR} (a, f, T (N), A) & := & \int_{T (N)}  \tmop{Cor} (a,
  f, T (N), A) d A,\\
  \tmop{TCOR} (k, a, f, T (N), A) & := & \int_{T (N)} \tmop{TCor} (k, a,
  f, T (N), A)  d A.
\end{eqnarray*}
There are numerous relations between these functions, but we will stop
the combinatorics here, coming back when we need it.

\section[The Random Variable $Z(n, F, T (N)) $ ]{The Random Variable \boldmath $ Z [n, F, T (N)] $ \unboldmath}
Define the random variable $Z [n, F, T (N)]$ by
\begin{equation} Z [n, F, T (N)] (A) = \frac{1}{N} \sum_{\#T = n} F
  \left( \frac{N}{2 \pi} \tmop{pr} (T) X (A) \right) , \end{equation}
where $\tmop{pr} (T)$ is the projection from $T (N)$ to $T (n)$,
$(x_1, \ldots, x_N) \mapsto (x_{t_1}, \ldots, x_{t_n})$ for a subset
$T \subset \{1, \ldots, N\}$ of cardinality $n$ and $X (A)$ is the
vector of angles for $A$.

We will later use this random variable with $ F = \tmop{TCor}$, but
for the start we formulate our version of Theorem 4.2.2 of
\cite{katzsarnak}.

The following theorem should be thought of as a very special limit
theorem for measures on the tori $T(N)$ as $N$ goes to infinity. We
fix a small torus of dimension $n$ and sum over all projections of
$T(N)$ to $T(n)$. In doing so we obtain induced measures on $T(n)$ and the
statement of the following theorem can be interpreted as stating that
these induced measures on $T(n)$ have a converging expectation value
and decreasing variance.

\begin{theorem}
  \label{thm:Z-estimation}Consider $n \in \mathbb{N}, n \geq 2$
  and $F \in \mathcal{T}_0 (n)$ with $\tmop{supp} F < \alpha$ for
  $\alpha > 0$.  Assume furthermore $F \geq 0$.
  \begin{enumerate}
  \item The sequence $E (Z [n, F, T (N)])$ converges for $N
    \rightarrow \infty$ to a limit $E (n, F, \tmop{univ})$ and 
    the estimation
    \begin{equation} 
      |E (Z [n, F, T (N)]) - E (n, F, \tmop{univ}) |
      \leq \|F\|_{\sup} \frac{1}{N} \frac{\alpha^{n - 1}}{(n - 2)
        !} . 
    \end{equation} 
    is true for all $N \geq 2$.

  \item For all $N \geq 2$ the expectation is bounded as follows:
    \begin{equation} |E (Z [n, F, T (N)]) | \leq \|F\|_{\sup}
      \frac{\alpha^{n - 1}}{(n - 1) !} . 
    \end{equation}

  \item For all $N \geq 2$ the variance is bounded as follows:
    \begin{equation} 
       \tmop{Var} (Z [n, F, T (N)]) \leq
           \frac{\|F\|_{\sup}^2}{N} \max\{1,(2 \alpha)^{2 n - 2}\} \frac{2
           n^2}{  \left( \tmop{floor} \left( \frac{n}{2} \right) ! \right)^2}, 
    \end{equation} 
    where $\tmop{floor}$ denotes the
    function rounding a real number down to the next integer.
  \end{enumerate}
\end{theorem}

\begin{proof}
  We start with the proof of statement 2.
  
  By a direct calculation we see that
  \begin{equation} E (Z [n, F, T (N)]) = \frac{1}{N} \int_{[0, N]^n}
    \binom{N}{n} \frac{1}{N^n} F (x) d x_1 \ldots d x_n.
  \end{equation} 
  Since $\tmop{supp} F < \alpha$, we consider the
  set $\Delta (n, \alpha) =\{x \in \mathbb{R}^n : \sup_{i, j} |x_i -
  x_j | < \alpha\}$. By Lemma 5.8.3 of \cite{katzsarnak}, we know
  \begin{equation} \label{eq:volume_of_Delta}
    \frac{1}{N} \tmop{Vol} (\Delta (n, \alpha) \cap [0, N]^n) \leq n
    \alpha^{n - 1}.
  \end{equation}
  Applying this to the above, it follows that
  \begin{equation}
    \begin{aligned}
      |E (Z [n, F, T (N)]) | & \leq \frac{1}{N} \|F\|_{\sup}
      \frac{1}{N^n}
      \binom{N}{n} N \alpha^{n - 1} n \\
      & \leq \|F\|_{\sup} \frac{\alpha^{n - 1}}{(n - 1) !} \frac{N}{N}
      \cdot \frac{N - 1}{N} \cdot \ldots \cdot \frac{N - n + 1}{N} .
    \end{aligned}
  \end{equation}
  The following inequality
  \begin{equation}
      \prod_{\nu = 1}^k  \left( 1 -
      \frac{\nu}{N} \right) \leq 1 - \frac{k}{N}, \label{induction-formula}
  \end{equation}
  gives
  \begin{equation} |E (Z [n, F, T (N)]) | \leq \|F\|_{\tmop{sup}}
    \frac{\alpha^{n - 1}}{(n - 1) !} \left( 1 - \frac{n - 1}{N}
    \right) \leq \|F\|_{\tmop{sup}} \frac{\alpha^{n - 1}}{(n - 1)
      !} . 
  \end{equation} 
  Thus, statement 2 has been proven.
  
  Now we wish to prove the first statement. For this, recall that $F
  \in \mathcal{T}_0 (n)$ means, that $F$ is $S_n$-invariant and
  invariant under diagonal addition. So
  \begin{equation} E (Z [n, F, T (N)]) = \frac{1}{N} \binom{N}{n}
    \frac{1}{N^n} n! \int_{[0, N]^n (\tmop{ordered})} F (x) d x_1
    \ldots d x_n . 
  \end{equation}
  This is true since the tuples with two or more equal components are
  a zero set and can be neglected. Substituting
  \begin{equation} y_1 = x_1, y_2 = x_2 - x_1, \ldots, y_n = x_n -
    x_1 
  \end{equation}
  yields
  \begin{equation}
    \begin{aligned}
      E (Z &[n, F, T (N)]) \\
      & = \frac{1}{N} \binom{N}{n} \frac{n!}{N^n} \int_0^N \left(
        \int_{[0, N - y_1]^{n - 1} (\tmop{ordered})} F (0, y_2,
        \ldots, y_n) d y_2 \ldots d y_n \right) dy_1.\
    \end{aligned}
  \end{equation}

  We call the inner integral $g (y_1)$ and assume that $\alpha<N$. Note that
  \begin{equation} g (y_1) \leq \|F\|_{\sup} \int _{[0,
      \alpha]^{n - 1} (\tmop{ordered})} d y_2 \ldots d y_n
    =\|F\|_{\sup} \frac{\alpha^{n - 1}}{(n - 1) !} 
  \end{equation}
  since $\tmop{supp} F < \alpha$. Therefore the integral extends
  from 0 to $\min (\alpha, N - y_1)$. We now set
  \begin{equation} E (n, F, \tmop{univ}) : = g_{\alpha} := \int_{[0,
      \alpha]^{n - 1} (\tmop{orderhned})} F (0, y_2, \ldots, y_n) d y_2
    \ldots d y_n
  \end{equation}
  and consider the difference
  \begin{equation} \label{eq:def_of_D}
    D = |E (Z [n, F, T (N)]) - E (n, F, \tmop{univ})| =
    \left|\frac{1}{N} \binom{N}{n} \frac{n!}{N^n} \int_0^N g (y_1) d y_1 -
    g_{\alpha}\right| . 
  \end{equation} 
  By splitting the integral into two parts, i.e., integrating from 0 to
  $N - \alpha$ and from $N - \alpha$ to $N$, it follows that in the
  first case $g (y_1) = g_{\alpha}$ because of the $\tmop{supp} F <
  \alpha$ condition, and thus
  \begin{equation} D = \left|\frac{1}{N} \binom{N}{n} \frac{n!}{N^n}
    \int_0^{N - \alpha} g_{\alpha} d y_1 + \frac{1}{N} \binom{N}{n}
    \frac{n!}{N^n} \int_{N - \alpha}^N g (y_1) d y_1 -
    g_{\alpha}\right|. 
  \end{equation} 
  Therefore we have
  \begin{equation} D \leq \left( \binom{N}{n} \frac{n!}{N^n}
      \frac{N - \alpha}{N} - 1 \right) g_{\alpha} + \frac{\alpha}{N}
    \binom{N}{n} \frac{n!}{N^n} g_{\alpha} 
  \end{equation} 
  which leads to
  \begin{equation} D \leq \left( \binom{N}{n} \frac{n!}{N^n} - 1
    \right) g_{\alpha} \leq \frac{n - 1}{N} g_{\alpha} \leq
    \|F\|_{\sup} \frac{\alpha^{n - 1}}{(n - 2) !}
    \frac{1}{N}. 
  \end{equation} 
  This proves statement 1, if $N>\alpha$.

  If $N\leq \alpha$, we define $g_\alpha$ as above, but immediately see that
  \begin{equation}
    E (n, F, \tmop{univ}) = g_{\alpha} = \int_{[0,
      N]^{n - 1} (\tmop{ordered})} F (0, y_2, \ldots, y_n) d y_2
    \ldots d y_n.
  \end{equation}
  Inserting this into (\ref{eq:def_of_D}) it follows that statement 1 is 
  fulfilled in this case as well.
  
  For the proof of the last statement let us first look at
  \begin{equation}
    E (Z [\dots]^2) = \frac{1}{N^2} \int_{T (N)} \sum_{\#T = n,\#S =
      n} F \left( \frac{N}{2 \pi} \tmop{pr} (T) X (A) \right) F \left(
      \frac{N}{2 \pi} \tmop{pr} (S) X (A) \right) d A,
  \end{equation}
  where the sum extends over all subsets $S$ and $T$ of cardinality
  $n$ of the set $\{1, \ldots, N\}$ and we write $E (Z [\dots]^2)$ for
  $E (Z [n, F, T (N)]^2)$.  This can be written as
  \begin{equation}
    \begin{aligned}
      E (Z [\dots]^2) & = \frac{1}{N^2} \sum_{l = n}^{2 n}
      \binom{N}{l}\frac{1}{N^l} \binom{l}{n} \binom{n}{l - n}\times\\
      &\times \int_{[0, N]^l} F (x_1, \ldots, x_n) F (x_{l -n + 1},
      \ldots, x_l) d x_1 \ldots d x_l,
    \end{aligned}
  \end{equation}
  as can seen by writing the double sum over $T$ and $S$ as a
  sum over the cardinality of $S \cup T$ and an inner sum. Using the
  $S_n$-invariance of $F$ one obtains the above formula. Now we
  consider the summands with $l < 2 n$. This means that
  \begin{equation} \tmop{supp} F (x_1, \ldots, x_n) F (x_{l - n +1},
    \ldots, x_l) \leq 2 \alpha 
  \end{equation} 
  because $|x_j - x_i | \leq |x_j - x_n | + |x_n - x_i |$ and $\tmop{supp} F
  \leq \alpha$. Using again (\ref{eq:volume_of_Delta}) and 
  (\ref{induction-formula}) we obtain
  \begin{equation}
    \begin{aligned}
      E (Z [\dots]^2) & \leq \frac{1}{N} \sum_{l = n}^{2 n - 1} \left(
        1 - \frac{l - 1}{N} \right) (2 \alpha)^{l - 1} \frac{l}{(l -
        n) !
        (l - n) ! (2 n - l) !} \|F\|_{\sup}^2\\
      & + \frac{1}{N^2} \binom{N}{2 n} \frac{1}{N^{2 n}} \binom{2
        n}{n} \left( \int_{[0, N]^n} F (x_1, \ldots, x_n) d x_1 \ldots
        d x_n \right)^2.
    \end{aligned}
  \end{equation}
  Now, compute the variance:
  \begin{multline}
    \tmop{Var}(E (Z [n, F, T (N)])) = E (Z [n, F, T (N)]^2) -
    E (Z [n, F, T (N)])^2\\
    \begin{aligned}
      & \leq \frac{\|F\|_{\sup}^2}{N} \sum_{l = n}^{2 n - 1} (2
      \alpha)^{l - 1} \frac{l}{\left( (l - n) ! \right)^2 (2 n - l) !}
      \\
      & + \frac{1}{N^2} \left( \int_{[0, N]^n} F (x_1, \ldots, x_n) d
        x_1 \ldots d x_n \right)^2 \left( \binom{2 n}{n} \frac{1}{N^{2
            n}} \binom{N}{2 n} - \left( \binom{N}{n} \frac{1}{N^n} \right)^2
\right)\\
      & \leq \frac{\|F\|_{\sup}^2}{N} \sum_{l = n}^{2 n - 1} (2
      \alpha)^{l - 1} \frac{l}{\left( (l - n) ! \right)^2 (2 n - l) !}\\
      & + \left( \alpha^{n - 1} \right)^2 n^2 \left( \binom{2 n}{n}
        \frac{1}{N^{2 n}} \binom{N}{2 n} - \left( \binom{N}{n}
          \frac{1}{N^n} \right)^2 \right) \|F\|_{\sup}^2.
    \end{aligned}
  \end{multline}
  But the last summand is negative:
  \begin{equation}
    \begin{aligned}
      \binom{2 n}{n} & \frac{1}{N^{2 n}} \binom{N}{2 n} - \left(
        \binom{N}{n}
        \frac{1}{N^n} \right)^2 \\
      & = \frac{1}{(n!)^2} \left( \prod_{\nu = 0}^{2 n - 1} \left( 1 -
          \frac{\nu}{N} \right) - \prod_{\nu = 0}^{n - 1} \left( 1 -
          \frac{\nu}{N} \right)^2   \right)\\
      & \leq 0
    \end{aligned}
  \end{equation}
  which gives the result
  \begin{equation}
    \tmop{Var} (E (Z [n, F, T (N)])) \leq \frac{\|F\|_{\sup}^2}{N} \max\{(2
    \alpha)^{2 n - 2},1\}\sum_{p = 0}^{n - 1} \frac{n + p}{(p!)^2 (n - p) !}\,.
    \end{equation}
  For simplicity we estimate further
  \begin{equation} 
    \sum_{p = 0}^{n - 1} \frac{n + p}{(p!)^2 (n - p) !}
    \leq 2 n \sum_{p = 0}^{n - 1} \frac{1}{n!} \frac{n!}{p!p! (n
      - p) !} \leq 2 n \sum_{p = 0}^{n - 1} \frac{1}{n!p!}
    \binom{n}{p}. 
  \end{equation} 
  For $n$ even this yields
  \begin{equation} 
      \binom{n}{p} \leq \frac{n!}{\left( \frac{n}{2}
      \right) ! \left( \frac{n}{2} \right) !}. 
  \end{equation} 
  and for $n$ odd
  \begin{equation} 
    \binom{n}{p} \leq \frac{n!}{\left( \frac{n +
          1}{2} \right) ! \left( \frac{n - 1}{2} \right) !}. 
  \end{equation} 
  So we have the following estimation for the variance:
  \begin{equation}
    \tmop{Var} (Z [n, F, T (N)]) \leq \frac{\|F\|_{\sup}^2}{N} 
      \max\{(2 \alpha)^{2 n - 2},1\} \frac{2
      n^2}{\left( \tmop{floor} \left( \frac{n}{2} \right) ! \right)^2}.
  \end{equation}
  Combining everything finishes the proof of the last statement.
\end{proof}

\section[Moving the Estimates to $\tmop{TCor} (k, a, f, T (N))$]{Moving the 
Estimates to \boldmath$\tmop{TCor} (k, a, f, T (N))$\unboldmath}
The reader may wonder how the above theorem is related to spectral
statistics. The answer is given by the following theorem which
transfers the above estimation on $Z [n, F, T (N)]$ to estimations
about $\tmop{TCor}$.

\begin{theorem}
  \label{thm:tcor-estimates}Let $f : \mathbb{R} \rightarrow
  \mathbb{R}_{\geq 0}$ be a bounded, non-negative,
  Borel-measurable function with upper bound $\alpha$ and $a, k \in
  \mathbb{N}$ with $k \geq a$.
  \begin{enumerate}
  \item The sequence $\tmop{TCOR} (k, a, f, T (N))$ converges for $N
    \rightarrow \infty$ to a limit which is denoted by $\tmop{TCOR} (k, a,
    f, \tmop{univ})$, and the following estimation
    \begin{equation} | \tmop{TCOR} (k, a, f, T (N)) - \tmop{TCOR} (k,
      a, f, \tmop{univ}) | \leq \binom{k}{a} \|f\|_{\sup}
      \frac{1}{N} \frac{\alpha^{k + 1}}{k!} 
    \end{equation} 
    holds for all $N \geq 2$.
  \item For all $N \geq 2$ the expectation is bounded as follows:
    \begin{equation} \tmop{TCOR} (k, a, f, T (N)) \leq \binom{k}{a}
      \|f\|_{\sup} \frac{\alpha^{k + 1}}{(k + 1) !} . \end{equation}
  \item For all $N \geq 2$ the variance is bounded as follows:
    \begin{multline}
      \tmop{Var} (A \mapsto \tmop{TCOR} (k, a, f, T (N), A) \text{ on
      } T
      (N)) \\
      \leq \left. \binom{k}{a} \right.^2 \frac{\|f\|_{\sup}^2}{N} 
      \max\{(2 \alpha)^{2 k + 2},1\} \frac{2 (k + 2)^2}{\left( \tmop{floor} \left(
            \frac{k}{2} + 1 \right) ! \right)^2}.
    \end{multline}
  \end{enumerate}
\end{theorem}

\begin{proof}
  This is Proposition 4.2.3 of \cite{katzsarnak}. For self-containtedness
  we give the proof here again.
  
  The idea is to use Theorem \ref{thm:Z-estimation} for the function
  $F (X) = \tmop{TClump} (k, a, f, k + 2, X)$, where $k = n + 2$. We
  claim that
  \begin{equation} Z [k + 2, F, T (N)] (A) = \tmop{TCor} (k, a, f, T
    (N), A) 
  \end{equation} 
  then. This can be seen by unwinding
  the definitions and using a combinatorial identity for equation
  (\ref{eq:tclump-tcor-in-proof}):
  \begin{eqnarray}
    Z [k + 2, F, T (N)] (A) & = & \frac{1}{N} \sum_{\#T = k + 2} F \left(
      \frac{N}{2 \pi} \tmop{pr} (T) X (A) \right) \nonumber\\
    & = & \frac{1}{N} \sum_{\#T = k + 2} \tmop{TClump} \left( k, a, f, k + 2,
      \frac{N}{2 \pi} \tmop{pr} (T) X (A) \right) \nonumber\\
    & = & \frac{1}{N} \tmop{TClump} \left( k, a, f, N, \frac{N}{2 \pi} X (A)
    \right)  \label{eq:tclump-tcor-in-proof}\\
    & = & \tmop{TCor} (k, a, f, T (N), A) . \nonumber
  \end{eqnarray}
  The theorem now follows from the fact that proved $F \in
  \mathcal{T}_0 (n)$ and $\|f\|_{\sup} \binom{k}{a} \geq
  \|F\|_{\sup}$. But these are direct consequences of the definition of
  $\tmop{TClump}$ as $\binom{k}{a} \tmop{Clump}$ and Lemma
  \ref{lem:properties-of-sep-and-clump}.
\end{proof}

\section[The Weak Convergence of $\mu (\tmop{naive}, U (N), 1)$ to
  the Poisson Distribution]{The Weak Convergence of \boldmath$ \mu (\tmop{naive}, U (N), 1)$\unboldmath\ to
  the Poisson Distribution}

We only cite a part of Proposition 2.9.1 of \cite{katzsarnak} here without
repeating the proof.

\begin{theorem}
  \label{thm:weak-convergence}{\dueto{Katz, Sarnak}}Assume that $a \in
  \mathbb{N}$ is fixed. If for every $k \in \mathbb{N}$, $k
  \geq a$ and every $f : \mathbb{R} \rightarrow \mathbb{R}$
  which is bounded, Borel measurable, non-negative and of compact
  support,
  \begin{equation} 
    \lim_{N \rightarrow \infty} \tmop{TCOR} (k, a, f, T
    (N)) = : \tmop{TCOR} (k, a, f, \tmop{univ}) 
  \end{equation} 
  exists and moreover
  \begin{equation} 
    \sum_{k \geq a} \tmop{TCOR} (k, a, f,
    \tmop{univ}) < \infty,
  \end{equation}
  then the limit measure $\mu (\tmop{naive}, a)$ exists and
  \begin{equation} 
    \int_{\mathbb{R}} f d \mu (\tmop{naive}, a) =
    \sum_{k \geq a} (- 1)^{k - a} \tmop{TCOR} (k, a, f,
    \tmop{univ}) . 
  \end{equation}
\end{theorem}

\begin{proof}
  \cite{katzsarnak}, p.~58 and following.
\end{proof}

Since we are dealing with $T (N)$, it is possible to give an explicit
formula for the Lebesgue density of $\mu (\tmop{naive}, 0)$.

\begin{theorem}
  The limit measure $\mu (\tmop{naive}, a)$ exists and for $a = 0$ it
  has the probability density $e^{- x}$.
\end{theorem}

\begin{proof}
  Here Theorem \ref{thm:weak-convergence} will be applied to prove the
  convergence result. By statement 1 of Theorem \ref{thm:tcor-estimates},
  we know the existence of $\tmop{TCOR} (k, a, f, \tmop{univ})$ and by
  statement 2 we see that
  \begin{equation} 
    \tmop{TCOR} (k, a, f, \tmop{univ}) \leq
    \binom{k}{a} \|f\|_{\sup} \frac{\alpha^{k + 1}}{(k + 1) !}
    \leq \|f\|_{\sup} \frac{(2 \alpha)^{k + 1}}{(k + 1)
      !}. 
  \end{equation} 
  Thus
  \begin{equation} 
    \sum_{k \geq a} \tmop{TCOR} (k, a, f,
    \tmop{univ}) \leq \|f\|_{\sup} \sum_{k = 0}^{\infty} \frac{(2
      \alpha)^{k + 1}}{(k + 1) !}  \leq \|f\|_{\sup} e^{2 \alpha}
    < \infty . 
  \end{equation} 
  It remains to prove the explicit form for
  $a = 0$. For this it suffices to calculate to
  \begin{equation} \int_{\mathbb{R}} f d \mu (\tmop{naive},
    0) 
  \end{equation} 
  for the characteristic functions of intervals
  of the form $[0, p]$. But this integral can be  calculated directly
  \begin{equation} \int_0^p d \mu (\tmop{naive}, 0) = \lim_{N
      \rightarrow \infty} \int_0^p d \mu (\tmop{naive}, T (N), 0)
    , 
  \end{equation} 
  where
  \begin{multline}    
      \int_0^p d \mu (\tmop{naive}, T (N), 0) = \frac{1}{N} \sum_{j
        = 1}^{N-1} N! \int_{T (N) (\tmop{ordered})} f \left( \frac{N}{2
          \pi} (x_{j + 1}
        - x_j) \right) dA\\
      = (N - 1) !  \sum_{j = 1}^{N-1} \int_0^1 \int_0^{x_N} \ldots
      \int_0^{x_{j + 2}} \times
      \times \int_{x_{j + 1} - p / N}^{x_{j + 1}} \int_0^{x_j}
      \ldots \int_0^{x_2} d x_1 \ldots d x_N.
  \end{multline}
  The desired result follows by evaluating the right-hand side. For this
  let us define the integrands $I_j$ by
  \begin{equation}
    \label{eq:mu_naive_weak}
    \int_0^p d \mu (\tmop{naive}, T (N), 0)  = (N-1)! \sum_{j=1}^{N-1} I_J, 
  \end{equation}
  By a direct calculation we derive the recursion formula
  \begin{equation}
   I_{j+1} = I_j - (-1)^j \frac{1}{N!}\binom{N}{j+1}\left( 
   \frac{p}{N}\right)^{j+1}
  \end{equation}
  and thus the explicit formula for the $I_j$
  \begin{equation}
     I_j = \frac{1}{N!} \sum_{k=1}^{j} \binom{N}{k} \left( 
     \frac{p}{N}\right)^k (-1)^{k+1}.
  \end{equation}
  Now, we insert this into (\ref{eq:mu_naive_weak}) and compare
  it to the power series for $1-\exp(-p)$
  \begin{multline}
     \int_0^p d \mu (\tmop{naive}, T (N), 0)  = (N-1)! \sum_{j=1}^{N-1} \frac{1}{N!} \sum_{k=1}^{j} 
     \binom{N}{k} \left( \frac{p}{N}\right)^k (-1)^{k+1} \\
     \begin{aligned}
       & = \frac{1}{N} \sum_{j=1}^{N-1}\sum_{k=1}^{j} \binom{N}{k}\left( \frac{p}{N} \right)^k (-1)^{k+1}
         = - \frac{1}{N} \sum_{j=1}^{N-1} \sum_{k=1}^{j} \frac{(-p)^k}{k!} 
         \prod_{\nu=1}^{k-1}\left(1-\frac{\nu}{N} \right) \\
       & = - \frac{1}{N} \sum_{k=1}^{N-1} \sum_{j=k}^{N-1} \frac{(-p)^k}{k!} 
         \prod_{\nu=1}^{k-1}\left(1-\frac{\nu}{N} \right)
         = \sum_{k=1}^{N-1} \frac{(-p)^k}{k!} \frac{N-k}{N}  
         \prod_{\nu=1}^{k-1}\left(1-\frac{\nu}{N} \right) \\
       & = \sum_{k=1}^{N-1} a_k(N) \frac{(-p)^k}{k!},
    \end{aligned}
  \end{multline}
  where the $a_k(N)$ are the coefficients defined above. For fixed $k$
  \begin{equation}
    a_k(N) =  \prod_{\nu=1}^{k}\left(1-\frac{\nu}{N} \right) \to 1 
    \text{ as } N \to \infty,
  \end{equation}
  which completes the proof.
\end{proof}

\section[The $M$-grid]{The \boldmath$M$\unboldmath-grid}

We would like to study the Kolmogorov-Smirnov distance $d_{\tmop{KS}}$
for the nearest neighbor measures $\mu (\tmop{naive}, A, T (N),
a)$. This is a quite complicated matter and therefore we discretize on
the so-called $M$-grid.

For this let $M$ be a (big) positive natural number. Divide the
interval $[0, 1]$ into pieces of length $\frac{1}{M}$. This defines a
grid
\begin{equation} - \infty = s (0) < s (1) < \ldots < s (M - 1) < s (M)
  = + \infty, 
\end{equation} 
where
\begin{equation} \int_{- \infty}^{s (j)} \mu (\tmop{naive}, a) =
  \frac{j}{M} \text{ for all } 1 \leq j \leq M -1. 
\end{equation}
\begin{definition}
  \label{def:d_ks_mgrid}
  Define the $M$-grid version of the Kolmogorov-Smirnov distance
  to be
  \begin{equation} d_{M, \tmop{KS}} (\mu, \nu) = \max_{i = 1, \ldots,
      M - 1} | \int_{s (1)}^{s (i)} d \mu - \int_{s (1)}^{s (i)} d \nu
    | . 
  \end{equation}
\end{definition}

\begin{lemma}
  For any Borel measure of total mass $\leq 1$ we have the
  inequality
  \begin{equation} d_{\tmop{KS}} (\nu, \mu (\tmop{naive}, a))
    \leq \frac{5}{M} + 2 \cdot d_{M, \tmop{KS}} (\nu, \mu
    (\tmop{naive}, a)) . \end{equation}
\end{lemma}

\begin{proof}
  See \cite{katzsarnak} p.81.
\end{proof}

\section{The Key Lemma}

For simplicity we cite here Lemma 3.2.16 of \cite{katzsarnak}.

\begin{lemma}
  Let $f \geq 0$ be a bounded, Borel measurable function with
  compact support and $L \geq a$ be an integer, then the
  following basic inequality
  \begin{multline}
    | \tmop{INT} (a, f, \tmop{univ}) - \tmop{Int} (a, f, T (N), A) | \\
    \begin{aligned}
      & \leq \sum_{L \geq k \geq a} | \tmop{TCOR} (k, a, f, T
      (N)) - \tmop{TCor} (k, a, f, T (N), A) | \\
      & + \sum_{L \geq k \geq a} | \tmop{TCOR} (k, a, f, T (N)) -
      \tmop{TCOR} (k, a, f, \tmop{univ})  |\\
      & + \tmop{TCOR} (L, a, f, \tmop{univ}) + \tmop{TCOR} (L + 1, a,
      f, \tmop{univ})
    \end{aligned}
  \end{multline}
  holds.
\end{lemma}

\begin{proof}
  See \cite{katzsarnak} p.83.
\end{proof}

\begin{notation}
  Since it is too cumbersome to write the measure $\mu (\tmop{naive},
  A, T (N), a)$, we abbreviate in the following
  \begin{equation} 
    \mu = \mu (\tmop{naive}, a) \text{ and } \mu_A =
    \mu (\tmop{naive}, A, T (N), a) 
  \end{equation} 
  for fixed $a$.
\end{notation}

\begin{corollary}
  Let $R \subset [s (1), s (M - 1)]$ be a Borel measurable. Then
  \begin{eqnarray*}
    | \mu (R) - \mu_A (R) | & \leq & \sum_{L \geq k \geq a} |
    \tmop{TCOR} (k, a, \chi, T (N)) - \tmop{TCor} (k, a, \chi, T (N), A) |
    \\
    & + & \binom{L}{a} \frac{\alpha^{L + 1}}{(L + 1) !} + \binom{L + 1}{a}
    \frac{\alpha^{L + 2}}{(L + 2) !} + \frac{1}{N} \sum_{L \geq k
      \geq a} \binom{k}{a} \frac{\alpha^{k + 1}}{k!}\\
    & < & \sum_{L \geq k \geq a} | \tmop{TCOR} (k, a, \chi, T (N))
    - \tmop{TCor} (k, a, \chi, T (N), A) | \\
    & + & \frac{(2 \alpha)^{L + 1}}{(L + 1) !} + \frac{(2 \alpha)^{L + 2}}{(L
      + 2) !} + \frac{1}{N} \sum_{L \geq k \geq a}  \frac{(2
      \alpha)^{k + 1}}{k!}
  \end{eqnarray*}
  where $\chi$ is the characteristic function of $R$ and $\alpha =
  \tmop{diam} (R)$.
\end{corollary}

\begin{proof}
  Apply the above Lemma and the $\tmop{TCOR}$ estimations.
\end{proof}

\begin{corollary}
  Set $\beta = s (M - 1) - s (1)$. The following estimation holds:
  \begin{equation}
    \begin{aligned}
      d_{M, \tmop{KS}} (\mu, \mu_A) & < \max_i \{ \sum_{L \geq k
        \geq a} | \tmop{TCOR} (k, a, \chi_{[s (1), s (i)]}, T (N)) \\
      & -\tmop{TCor} (k, a, \chi_{[s (1), s (i)]}, T (N), A) |  \} \\
      & + \frac{(2 \beta)^{L + 1}}{(L + 1) !} + \frac{(2 \beta)^{L +
          2}}{(L + 2) !} + \frac{1}{N} \sum_{L \geq k \geq
        a} \frac{(2 \beta)^{k + 1}}{k!},
    \end{aligned}
  \end{equation}
  where $\chi_R$ denotes the characteristic function of the interval
  $R$.
\end{corollary}

\begin{proof}
  This is clear from the definition of $d_{M, \tmop{KS}}$.
\end{proof}

\begin{lemma}
  \begin{multline}
    \int_{T (N)}  | \tmop{TCOR} (k, a, \chi_{[s (1), s (i)]}, T (N)) -
    \tmop{TCor} (k, a, \chi_{[s (1), s (i)]}, T (N), A) |  d A  \\
    \leq \binom{k}{a}  \sqrt{\frac{2}{N}} \max \{(2s (i) - 2s (1))^{k + 1},1\} \frac{k +
      2}{\tmop{floor} ( \frac{k}{2} + 1) !}.  
  \end{multline}
\end{lemma}

\begin{proof}
  This is just the Cauchy-Schwarz inequality
  \begin{equation}
    \int_{T(N)} |h(A)| dA \leq \sqrt{\int_{T(N)} |h(A)|^2 dA}\,,
  \end{equation}
  where $h(A) = \tmop{TCOR} (k, a, \chi_{[s (1), s (i)]}, T (N)) -
    \tmop{TCor} (k, a, \chi_{[s (1), s (i)]}, T (N), A) $\\ 
  combined with the statement about the variance of the 
  $\tmop{TCOR}$ estimations.
\end{proof}

\begin{theorem}
  For $N \rightarrow \infty$
  \begin{equation} \int_{T (N)} d_{\tmop{KS}} (\mu, \mu_A) d A
    \rightarrow 0.
  \end{equation}
\end{theorem}

\begin{proof}
  Putting everything together, we obtain
  \begin{multline}
    \int_{T (N)} d_{\tmop{KS}} (\mu, \mu_A) d A < \frac{5}{M} + 2
    \cdot \left( \frac{(2 \beta)^{L + 1}}{(L + 1) !} + \frac{(2 \beta)^{L +
          2}}{(L + 2) !} + \frac{1}{N} \sum_{L \geq k \geq 0} 
      \frac{(2 \beta)^{k + 1}}{k!} \right)\\
    +  \sqrt{\frac{2}{N}} \max_i \left\{ \sum_{L \geq k
        \geq a}\binom{k}{a}  \frac{k + 2}{\tmop{floor} (
        \frac{k}{2} + 1) !} \max\{1, ( 2s (i) - 2s (1))^{k + 1}\} \right\}
  \end{multline}
  where $\beta=s(M-1) - s(1)$ as above.

  Now we combine the summands to make more explicit estimations
  \begin{equation} \sum_{L \geq k \geq 0} \frac{(2 \beta)^{k
        + 1}}{k!} < (2 \beta) e^{2 \beta} 
  \end{equation} 
  and if $s(i) - s(1) \geq \frac{1}{2}$  we have the following estimation for the second sum
  \begin{equation}
    \begin{aligned}
      \sum_{L \geq k \geq 0} (4 (s (i) - s (1)))^{k + 1}
      \frac{k + 2}{\tmop{floor} ( \frac{k}{2} + 1) !} & \leq 8 (s
      (i) - s (1)) \sum_{L \geq k \geq 0} \frac{(64 (s (i) -
        s (1))^2)^{k /
          2}}{\tmop{floor} ( \frac{k}{2} + 1) !}\\
      & < 8 (s (i) - s (1)) e^{64 (s (i) - s (1))^2}.
    \end{aligned}
  \end{equation}
  If $s(i) - s(1) < \frac{1}{2} $ we may estimate the sum as
  \begin{equation}
    \sum_{L \geq k \geq 0} 2^k \frac{k+2}{\tmop{floor} ( \frac{k}{2} + 1) !}
    \leq \sum_{L \geq k \geq 0} 3 \cdot 2^k \leq 3L \cdot 2^L.
  \end{equation}
  Applying this to the above it follows that
  \begin{equation}
    \begin{aligned}
      \int_{T (N)} d_{\tmop{KS}} (\mu, \mu_A) dA < & \frac{5}{M} + 2
      \left( \frac{(2 \beta)^{L + 1}}{(L + 1) !} + \frac{(2 \beta)^{L
            + 2}}{(L + 2) !}  \right)  \\ & + \frac{1}{\sqrt{N}} 
       \left(
        \frac{1}{\sqrt{N}} (2 \beta) e^{2
          \beta} + 8 \beta \sqrt{2} e^{64 \beta^2} + 6L \cdot 2^L \right) .
    \end{aligned}
  \end{equation}
  It is clear that $\beta$ depends only on $M$. So given $\varepsilon
  > 0$, we first choose $M$ so large, that
  \begin{equation} 
     \frac{5}{M} < \frac{\varepsilon}{3},
  \end{equation}
  then we can choose $L$ so large, that
  \begin{equation} 
      \frac{(2 \beta)^{L + 1}}{(L + 1) !} + \frac{(2
      \beta)^{L + 2}}{(L + 2) !}  <
    \frac{\varepsilon}{6} 
  \end{equation} 
  and finally $N$ so large that
  \begin{equation} \frac{1}{\sqrt{N}} \left( \frac{1}{\sqrt{N}} (2
      \beta) e^{2 \beta} + 8 \beta \sqrt{2} e^{64 \beta^2} 
      + 6L \cdot 2^L\right) <
    \frac{\varepsilon}{3} .
  \end{equation}
\end{proof}

\section{The Final Estimation}
In this last section we will give the final form of the
estimation. But before we do so we state a series of lemmas which we
will combine to give the main estimation.

We start by fixing two positive constants $\alpha,\gamma\in
\mathbb{R}_{>0}$. Set the grid size $M$ to be the largest integer
smaller than $e^{\alpha \sqrt{\log N} }$ and the cut-off $L$ to be the
largest integer such that
\begin{equation}
  (L-1)! \leq N^{\gamma^2} \leq L!\,.
\end{equation}
Then
\begin{equation}
 \log M \leq \alpha\sqrt{\log N} \leq \log (M+1)       
\end{equation}
and
\begin{equation}
  \log (L-1)! \leq \gamma^2 \log N \leq \log L!\,.
\end{equation}
Thus, we see that
\begin{equation}
  \label{eq:logMbylogN}
  \log M \leq \frac{\alpha}{\gamma}\sqrt{\log L!}\,.
\end{equation}

The following lemma is a useful corollary of Stirling's formula. 
\begin{lemma}
  \label{lem:sterling}
  Given $\epsilon>0$ and $c>0$, there exists a $k_0$ such that for
  all $k\geq k_0$:
  \begin{enumerate}
  \item $(\log k!)^{k+2} \leq (k!)^{1+\epsilon}$.
  \item $c^{k+2} \leq (k!)^{\epsilon/2}$.
  \end{enumerate}
\end{lemma}

\begin{proof}
  The proof can be found in \cite{katzsarnak} p.93.
\end{proof}

Next, note that $\beta=s(M-1)-s(0)) < \log M$ by construction. We will
now give estimations for each summand in (\ref{eq:collected-estimate}).

\begin{lemma}
  The following estimation holds:
  \begin{equation}
    \frac{(2\beta)^{L+1}}{(L+1)!} + \frac{(2\beta)^{L+2}}{(L+2)!}
  \leq \frac{1}{N^{\gamma^2-\epsilon\gamma^2}}.  
  \end{equation}
\end{lemma}

\begin{proof}
  By lemma \ref{lem:sterling} we see that
  \begin{equation}
    \begin{aligned}
    \frac{(2\beta)^{L+1}}{(L+1)!} + \frac{(2\beta)^{L+2}}{(L+2)!}
    & \leq \frac{2 (2\beta)^{L+2}}{(L+1)!} 
    \leq \frac{2}{L+1}\frac{(2\log M)^{L+2}}{L!} \\
%\intertext{by (\ref{eq:logMbylogN})}    
    & \leq \frac{2}{L+1}\frac{\sqrt{L!}^{1+\epsilon}}{L!}
    \left(\frac{2\alpha}{\gamma}\right)^{L+2} \\
    &\leq (L!)^{\epsilon-1} \leq
    \frac{1}{N^{\gamma^2-\epsilon\gamma^2}}, 
    \end{aligned}
  \end{equation}
  which is the desired result.
\end{proof}

\begin{lemma}
  \begin{equation}
    \frac{1}{N}\sum_{L\geq k\geq 0} \frac{(2\beta)^{k+1}}{k!} 
    \leq \frac{1}{\sqrt{N}} \sum_{L\geq k\geq 0} 
    \frac{(4\log M)^{k+1}(k+1)}{\tmop{floor}\left(
     \frac{k}{2} + 1\right)!}.
  \end{equation}
\end{lemma}

\begin{proof}
  This follows by direct calculation.
\end{proof}

\begin{lemma}
  \begin{equation}
   \frac{1}{\sqrt{N}} \sum_{L\geq k\geq 0} 
    \frac{(4\log M)^{k+1}(k+2)}{\tmop{floor}\left(
     \frac{k}{2} + 1\right)!} 
   \leq \frac{1}{\sqrt{N}} N^{\gamma^2+2\gamma^2\epsilon}
  \end{equation}
\end{lemma}

\begin{proof}
  \begin{equation}
    \begin{aligned}
      \frac{1}{\sqrt{N}} \sum_{L\geq k\geq 0} 
    \frac{(4\log M)^{k+1}(k+2)}{\tmop{floor}\left(
     \frac{k}{2} + 1\right)!}
      & \leq \frac{1}{ \sqrt{N} } \sum_{L\geq k\geq 0} (8\log M)^{k+1}
      \\
      & \leq \frac{1}{\sqrt{N}} L (8\log M)^{L+1} \leq
      \frac{1}{\sqrt{N}} (16\log M)^{L+1} \\
      & \leq \frac{1}{\sqrt{N}}
      \left(\frac{16\alpha}{\beta}\right)^{L+1} \sqrt{\log L!}^{L+1}
      \leq \frac{1}{\sqrt{N}} \left(L!\right)^{\frac{1}{2}+\epsilon} \\
      & \leq \frac{ \left( N^{2\gamma^2}\right)^{(\frac{1}{2}+\epsilon)}
}{\sqrt{N}},
    \end{aligned}
  \end{equation}
  where  in the last line we used that
  \begin{equation}
    L! \leq L N^{\gamma^2} \leq N^{2\gamma^2}.
  \end{equation}
\end{proof}

\begin{lemma}
  \begin{equation}
      \frac{1}{\sqrt{N}}\sum_{L \geq k \geq 0} \binom{k}{a} \frac{k+2}{\tmop{floor} ( \frac{k}{2} + 1) !} 
   \leq 3 N^{2\gamma^2-\frac{1}{2}} \text{ for sufficiently large }L. 
  \end{equation}
\end{lemma}
   
\begin{proof}
  \begin{equation}
    \sum_{L \geq k \geq 0} \binom{k}{a} \frac{k+2}{\tmop{floor} ( \frac{k}{2} + 1) !} 
    \leq 3L\cdot 2^L \leq 3 L! \text{ for sufficiently large }L.
  \end{equation}

\end{proof}

Now, we want to combine these estimations. Starting with equation
(\ref{eq:collected-estimate}) 

\begin{multline}\label{eq:collected-estimate}
    \int_{T (N)} d_{\tmop{KS}} (\mu, \mu_A) d A < \frac{5}{M} + 2
    \cdot \left( \frac{(2 \beta)^{L + 1}}{(L + 1) !} + \frac{(2 \beta)^{L +
          2}}{(L + 2) !} + \frac{1}{N} \sum_{L \geq k \geq 0} 
      \frac{(2 \beta)^{k + 1}}{k!} \right)\\
    +  \sqrt{\frac{2}{N}} \max_i \left\{ \sum_{L \geq k
        \geq a}\binom{k}{a}  \frac{k + 2}{\tmop{floor} (
        \frac{k}{2} + 1) !} \max\{1, ( 2s (i) - 2s (1))^{k + 1}\} \right\}
  \end{multline}

the following intermediary result is a consequence of the above lemmas:

\begin{equation}
  \int_{T(N)} d_{KS}(\mu,\mu_A) dA \leq \frac{5}{M} +
  \frac{2}{N^{\gamma^2-\epsilon\gamma^2}}
  + \frac{2\sqrt{2}}{N^{\frac{1}{2}-\gamma^2-2\gamma^2\epsilon}}
  + \frac{3\sqrt{2}}{N^{\frac{1}{2}-2\gamma^2}}.
\end{equation}

The summand $\frac{5}{M}$ decreases like $\exp(-\alpha\sqrt{\log N})$
as $N$ goes to infinity. The other summands decrease much faster.
Therefore we may neglect them, i.e. for N sufficiently large, the
left-hand side is smaller than $\frac{6}{M}$. If we substitute $\alpha$
from the beginning by $\alpha/2$ the constant $6$ can also be neglected.
Thus, we have proved the main theorem of this chapter.
\begin{theorem}
  \label{thm:main-theorem-tn}
  Let $\alpha$ be a positive constant. Then  the following estimation
  \begin{equation}
    \int_{T(N)} d_{KS} ( \mu, \mu_A ) dA < \frac{1}{e^{\alpha\sqrt{\log N}}}
  \end{equation}
  holds for $N$ sufficiently large.
\end{theorem}

%%% Local Variables:  
%%% mode: latex 
%%% TeX-master: "../thesis" 
%%% End: 

% \chapter{Appendix}
\chapter{Appendix}
\label{chap:gen-level-spacings}

In this Appendix the fundamental results from representation theory
and momentum geometry which are used in the main body of the text are
stated in detail. With few exceptions, for the proofs only references
to the literature are given. We close this Appendix with elementary
observations about nearest neighbor statistics.

\section{Representation Theory}
Throughout this text we are concerned with the representation theory
of compact Lie groups. For the standard facts we refer the reader to
\cite{broecker:compactlie} and \cite{knapp}.

We will always use the following conventions: $K$ denotes a semi-simple,
compact Lie group with Lie algebra $\mathfrak{k}$. %
\nomenclature{$K$}{A semi-simple, compact Lie group}%
\nomenclature{$\mathfrak{k}$}{The Lie algebra of $K$}%
Further, let $G$ denote the complexification of $K$ and $\mathfrak{g}$ be the
Lie algebra of $G$. %
\nomenclature{$G$}{The complexification of $K$}%
\nomenclature{$\mathfrak{g}$}{The Lie algebra of $G$}%
Furthermore for any unitary vector space $V$ the symbol $\tmop{U}(V)$ is used
for the set of unitary automorphisms.%
\nomenclature{$\tmop{U}(V)$}{The set of unitary operators on the unitary
  vector space $V$} % 

\subsection{Representations of Compact Lie Groups}
Fix a maximal torus $T$ in $K$ with Lie algebra $\mathfrak{t}$, i.e.~$T$
is a maximal, connected, commutative subgroup of $K$, %
\nomenclature{$T$}{A maximal torus of $K$}%
\nomenclature{$\mathfrak{t}$}{The Lie algebra of $T$}%
and every irreducible representation can be decomposed into one
dimensional representations of $T$. On each of these $T$ acts by
scalar multiplication, i.e., we are given a group homomorphism $f:T\to
S^1\subset{\mathbb{C}^*}$. We make the following definition.
\begin{definition}
  Let $\rho:K\to U(V)$ be an irreducible, unitary representation of
  $K$ on some finite dimensional vector space $V$. Then a \textbf{weight} of
  $\rho$ is an element $\lambda\in\mathfrak{t}^*$ such that there
  exists a non-trivial subspace $V_\lambda$ of $V$ with
\begin{equation}
  \label{eq:weightdefintion}
  d_e\rho(t).x=2\pi i\lambda(t) x\ \forall\ x\in V_\lambda,\ t\in \mathfrak{t}.
\end{equation}
\end{definition}
Note that these weights are sometimes called real infinitesimal
weights.

\begin{proposition}
  The set of weights (with multiplicity) of an irreducible
  representation determines the representation uniquely.
\end{proposition}

\begin{proof}
  This is a very weak form of the Theorem 5.110 in \cite{knapp}.  
\end{proof}

Moreover, one can order the set of all weights such that every
irreducible representation has a unique highest weight. We will define
such an ordering here, but we have to elaborate on the weights first.

Recall that the adjoint representation $Ad:K\to GL(\mathfrak{k})$ is
given by $k\mapsto d_e\tmop{int}(k)$, where $\tmop{int}:K\to
\tmop{Aut}(K), k \mapsto (g \mapsto kgk^{-1})$. The complexified
weights of the adjoint representation are called \textbf{roots}.

Of all $Ad$-invariant scalar products on $\mathfrak{g}$ the most
important one is the so called \textbf{Killing form} $\langle\cdot ,
\cdot\rangle_{\tmop{Kil}}$%
\nomenclature{$\langle\cdot , \cdot\rangle_{\tmop{Kil}}$}{The
  Killing-Form on $\mathfrak{k}$ or $\mathfrak{g}$}%
, which is defined by
\begin{equation}
  \langle \xi, \eta \rangle_{\tmop{Kil}} = \tmop{trace}({ad(\xi)\circ ad(\eta)}),
\end{equation}
where $\xi,\eta\in\mathfrak{g}$ and $ad:\mathfrak{g}\to \tmop{End} 
(\mathfrak{g}), \xi \mapsto [\xi,\cdot]$. 

Let us denote the set of roots by $\Delta$. %
\nomenclature{$\Delta$}{The set of roots of $K$}%
The following lemma summarizes some properties of the roots.
\begin{lemma}
    The set $\Delta$ has the following properties:
    \begin{enumerate}
    \item $\{\alpha\in\Delta\}$ generates $\mathfrak{t}^\ast$.
    \item $\alpha\in\Delta$ if and only if $-\alpha\in\Delta$. 
    \item There exists  a set of simple roots, i.e.\ a smallest subset
      $\Delta'$ of $\Delta$, such that every $\alpha\in\Delta$ is an
      integer combination of simple roots.
    \item In the integer combination either all coefficients are
      non-negative or all are non-positive. 
    \item The simple roots form a basis for $\mathfrak{t}^\ast$.
    \item The non-negative linear combinations over $\mathbb{R}$ of
      the simple roots give a closed convex cone in
      $\mathfrak{t}^\ast$.
  \end{enumerate}
\end{lemma}

\begin{proof}
  Cf.~\cite{knapp} Chapter II.5.
\end{proof}
The cone in the lemma above is usually called the \textbf{Weyl
  chamber} with respect to the system of simple roots.  Identifying
$\mathfrak{t}^\ast$ with $\mathfrak{t}$ via an $Ad$-invariant scalar
product we can think of this cone as a subset of $\mathfrak{t}$.

Since every root is an integer combination of the simple roots, where
all coefficients are either non-negative or non-positive, we divide
the set $\Delta$ into the set of \textbf{positive roots}
\begin{equation}
  \Pi_+ = \{ \alpha\in \Delta \,:\, \alpha\text{ is non-negative
    combination of simple roots} \} 
\end{equation}%
\nomenclature{$\Pi_+$}{The set of positive roots of $K$}%
and \textbf{negative roots}
\begin{equation}
  \Pi_- = \{ \alpha\in \Delta \,:\, \alpha\text{ is non-positive
    combination of simple roots} \}.
\end{equation}%
\nomenclature{$\Pi_-$}{The set of negative roots of $K$}%
The simultaneous eigenspace of a root $\alpha$ is denoted by
$\mathfrak{g}_{\alpha}$, i.e.\
\begin{equation}
  \mathfrak{g}_{\alpha} = \{\xi \in \mathfrak{g}: \alpha(\tau)\xi =
    [\tau,\xi]\}.   
\end{equation}
This yields a direct sum decomposition of the Lie algebra
$\mathfrak{g}$
\begin{equation}
  \mathfrak{g} = \mathfrak{u}_+ \oplus \mathfrak{t}^{\mathbb{C}}
                 \oplus \mathfrak{u}_-, 
\end{equation}
where
\begin{equation} 
  \mathfrak{u}_- := \bigoplus_{\alpha\in \Pi_-}
  \mathfrak{g}_{\alpha}\quad \text{and}\quad 
  \mathfrak{u}_+ := \bigoplus_{\alpha\in \Pi_+}
  \mathfrak{g}_{\alpha}. 
\end{equation}%
\nomenclature{$\mathfrak{u}_-$}{The Lie algebra of negative roots}%
\nomenclature{$\mathfrak{u}_+$}{The Lie algebra of positive roots}% 
 \begin{definition}
  The group generated by the reflections on the faces of the Weyl
  chamber is called the \textbf{Weyl group}.
\end{definition}
We denote the Weyl group by $W$%
\nomenclature{$W$}{The Weyl group of $K$ with respect to fixed $\Pi_+$}% 
and remark that it is a finite group.
\begin{definition}
  The \textbf{ordering of weights} is given by
  \begin{equation}
    \label{eq:1}
    \lambda \leq \mu :\Leftrightarrow 
    \tmop{Conv}(W.\lambda)\subset \tmop{Conv}(W.\mu),
  \end{equation} 
  where $\lambda,\mu$ are weights.
\end{definition}
\begin{lemma}
  Every weight is equivalent to a weight in the Weyl chamber under the
  action of the Weyl group. 
\end{lemma}
\begin{proof}
  Cf.~\cite{knapp} Corollary 2.68.
\end{proof}

The main statement about weights is called the Theorem of the Highest
Weight.
\begin{theorem}
  Every irreducible representation has a unique highest weight in the
  Weyl chamber. Moreover, two irreducible representations are
  equivalent if and only if the highest weights are equal.
\end{theorem}

\begin{proof}
  Cf.~\cite{knapp} Theorem 5.110.
\end{proof}

Connected to the above definitions are special complex subgroups of
$G$, which are introduced subsequently.
\begin{definition}
  A \textbf{Borel subgroup} of $G$ is a maximal, connected, solvable,
  complex subgroup of $G$. A \textbf{parabolic subgroup} is a complex
  subgroup which contains a Borel subgroup.
\end{definition}

Given a fixed torus and a notion of positivity of roots, we have two
natural Borel subgroups, which are called $B_+$ and $B_-$. These can be
obtained as follows:
\begin{equation}
  B_-:=\exp(\mathfrak{u}_-\oplus \mathfrak{t}^{\mathbb{C}}) \quad
  \text{and} \quad 
  B_+:=\exp(\mathfrak{u}_+\oplus \mathfrak{t}^{\mathbb{C}}).
\end{equation}%
\nomenclature{$B_-$}{The Borel subgroup of negative roots of $G$}%
\nomenclature{$B_+$}{The Borel subgroup of positive roots of $G$}%

\subsection{The Universal Enveloping Algebra}
Let $\mathcal{T}(\mathfrak{g})$ %
\nomenclature{$\mathcal{T}(\mathfrak{g})$}{The full tensor algebra of
  $\mathfrak{g}$}%
denote the full tensor algebra of $\mathfrak{g}$, .i.e.\
$\mathcal{T}(\mathfrak{g}) = \oplus_{j\in\mathbb{N}}
(\otimes^j\mathfrak{g}).$ The \textbf{universal enveloping algebra}
$\mathcal{U}(\mathfrak{g})$ %
\nomenclature{$\mathcal{U}(\mathfrak{g})$}{The universal enveloping
  algebra of $\mathfrak{g}$}% 
of $\mathfrak{g}$ is given by the quotient
algebra
\begin{equation}
  \mathcal{U}(\mathfrak{g}) = \mathcal{T}(\mathfrak{g})/\mathcal{I},
\end{equation}
where $\mathcal{I}$ is the ideal generated by all $\langle \xi\otimes \eta-\eta\otimes
\xi-[\xi,\eta]\rangle$ for $\xi,\eta\in\mathfrak{g}$.

One directly checks that $\mathcal{U}(\mathfrak{g})$
is an associative algebra. 

\begin{theorem}
  The universal enveloping algebra $\mathcal{U}(\mathfrak{g})$ has the
  following properties:
  \begin{enumerate}
    \item $\mathfrak{g}$ is embedded in $\mathcal{U}(\mathfrak{g})$ by
      $X\mapsto X + \mathcal{I}$.
    \item Every Lie algebra representation $\rho_\ast: \mathfrak{g}
      \to \tmop{End}(V)$ has a continuation as a homomorphism of
      associative algebras $\rho_\ast: \mathcal{T} (\mathfrak{g}) \to
      \tmop{End}(V)$. The kernel of $\rho_\ast$ contains $\mathcal{I}$
      so this yields an induced homomorphism of associative algebras
      $\rho_\ast: \mathcal{U} (\mathfrak{g}) \to \tmop{End}(V)$.%
      \nomenclature{$\rho_\ast$}{The Lie algebra representation associated
        to a Lie group representation $\rho$}%
    \item (Lemma of Burnside) Let $\rho_\ast: \mathfrak{g} \to
      \tmop{End}(V)$ be an irreducible Lie algebra representation on a
      finite dimensional vector space. Then $\rho_\ast: \mathcal{U}
      (\mathfrak{g}) \to \tmop{End}(V)$ is surjective.
    \item (Theorem of Poincare-Birkhoff-Witt) Let $\xi_1,\dots,\xi_n$
      be a basis of $\mathfrak{g}$. Then the map
      \begin{equation}
        \psi :\mathbb{C}[X_1,..,X_n] \to \mathcal{U}(\mathfrak{g}),\ 
        \sum_I a_IX^I \to \sum_I a_I \xi^I
      \end{equation}
      is an isomorphism of vector spaces, where it is assumed that 
      every monomial in $\mathbb{C}[X_1,..,X_n]$ %
      \nomenclature{$\mathbb{C}[X_1,..,X_n]$}{The Ring of polynomials in
        $n$ indeterminates with coefficients in $\mathbb{C}$}%
      is ordered lexicographically. 
  \end{enumerate}
\end{theorem}

\begin{proof}
  The proof of the Lemma of Burnside can be found in \cite{farenick}
  Chapter 3.3. The rest is proved in \cite{knapp} Chap.\ III.
\end{proof}

Note that $\psi$ is not an isomorphism of algebras since
$\mathbb{C}[X_1,..,X_n]$ is commutative and $\mathcal{U}
(\mathfrak{g}) $ is not.  

In the text a notion of hermitian operators on the tensor algebra and
on the universal enveloping algebra is needed.

\begin{definition}\label{def:formaladjoint_on_tg}
  The $\mathbb{R}$-linear map 
  $\dagger:\mathcal{T}(\mathfrak{g})\to \mathcal{T}(\mathfrak{g})$ 
  defined by
  \begin{enumerate}
  \item $(z\alpha_1\otimes\dots\otimes\alpha_n)^\dagger = \bar{z}\alpha_n^\dagger
    \dots \alpha_1^\dagger\ \ \forall
    \alpha_1,\dots,\alpha_n\in\mathfrak{g},z\in\mathbb{C}$
  \item $\xi^\dagger=-\xi\ \ \forall \xi\in\mathfrak{k}$
  \end{enumerate}
  and, extended by $\mathbb{R}$-linearity to $\mathcal{T}(\mathfrak{g})$, is
  called the \textbf{formal adjoint}. %
  \nomenclature{$\dagger$}{The formal adjoint on
    $\mathcal{T}(\mathfrak{g})$ or $\mathcal{U}(\mathfrak{g})$}%
  An operator $\alpha\in \mathcal{T}(\mathfrak{g})$ is called
  \textbf{abstractly self-adjoint} or \textbf{abstractly hermitian},
  if $\alpha^\dagger=\alpha$.
\end{definition}
Note that the formal adjoint is not complex linear because of the
conjugation involved in condition 1.
\begin{remark}
  The map $\dagger$ is compatible with $\rho$ in the following sense:
  \begin{equation}
    \label{eq:compatibiliyofdagger}
    \rho_\ast(\xi^\dagger)=\rho_\ast(\xi)^\dagger.
  \end{equation}
\end{remark}

\begin{lemma}
  The map $\dagger$ induces a $\mathbb{R}$-linear map $\mathcal{U}
  (\mathfrak{g}) \to \mathcal{U} (\mathfrak{g})$, which we also call
  $\dagger$.
\end{lemma}

\begin{proof}
  We have to show that the ideal $\mathcal{I}$ in $\mathcal{T}
  (\mathfrak{g})$ is fixed by $\dagger$. For this let $\xi = \xi_1 + i
  \xi_2$ and $\eta = \eta_1 + i \eta_2$ with $\xi_1, \xi_2, \eta_1$
  and $\eta_2 \in \mathfrak{k}$ be given. We calculate
  \begin{equation}
    \begin{aligned}
     & \left( (\xi_1 + i \xi_2) (\eta_1 + i \eta_2) - (\eta_1 + i \eta_2) (\xi_1
     + i \xi_2) - [\xi_1 + i \xi_2, \eta_1 + i \eta_2] \right)^{\dagger}\\
     & = (\xi_1 \eta_1 - \eta_1 \xi_1 - [\xi_1, \eta_1])^{\dagger} + (i
    (\xi_2 \eta_1 - \eta_1 \xi_2 - [\xi_2, \eta_1]))^{\dagger} \\
    & + (i (\xi_1 \eta_2 - \eta_2 \xi_1 - [\xi_1, \eta_2]))^{\dagger}
     - (\xi_2 \eta_2 - \eta_2 \xi_2 - [\xi_2, \eta_2])^{\dagger}\\
    & = (\eta_1^{\dagger} \xi_1^{\dagger} - \xi_1^{\dagger}
    \eta_1^{\dagger} - [\xi_1, \eta_1]^{\dagger}) - i (\eta_1^{\dagger}
    \xi_2^{\dagger} - \xi_2^{\dagger} \eta_1^{\dagger} - [\xi_2,
    \eta_1]^{\dagger}) - i (\eta_2^{\dagger} \xi_1^{\dagger} - \xi_1^{\dagger}
    \eta_2^{\dagger} - [\xi_1, \eta_2]^{\dagger})\\
    & - (\eta_2^{\dagger} \xi_2^{\dagger} - \xi_2^{\dagger}
    \eta_2^{\dagger} - [\xi_2, \eta_2]^{\dagger})\\
    & = (\eta_1 \xi_1 - \xi_1 \eta_1 - [\eta_1, \xi_1]) - i (\eta_1 \xi_2 -
    \xi_2 \eta_1 - [\eta_1, \xi_2]) - i (\eta_2 \xi_1 - \xi_1 \eta_2 -
    [\eta_2, \xi_1])\\
    & - (\eta_2 \xi_2 - \xi_2 \eta_2 - [\xi_2, \eta_2]).
    \end{aligned}
  \end{equation}
  This proves the lemma. 
\end{proof}

\subsection{The Laplace Operator}
A Casimir operator is by definition an element of the center
$\mathcal{Z}(\mathfrak{g})$ of $\mathcal{U}(\mathfrak{g})$. If we
consider an irreducible representation
$\rho_{\ast}:\mathcal{U}(\mathfrak{g})\to \tmop{End}(V)$, then due to
Schur's Lemma every Casimir operator has to act by scalar
multiplication.

The most important example of a Casimir operator is the Laplace
operator $\Omega$. %
\nomenclature{$\Omega$}{The Laplace operator in
  $\mathcal{U}(\mathfrak{g})$}%
Sometimes it is even called \emph{the} Casimir
element, e.g. in \cite{knapp}. We do not give an explicit formula for
the Laplace operator here, but just state that it is an operator of
degree two in the basis elements of $\mathfrak{g}$.

Let $\delta$ denote half the sum of positive roots.
\begin{lemma}
  \label{lem:laplace-operator}
  The Laplace operator $\Omega$ operates by the scalar
  $\langle\lambda,\lambda+2\delta\rangle_{\tmop{Kil}}$ in an
  irreducible representation of $\mathfrak{g}$ of highest weight
  $\lambda$.
\end{lemma}

\subsection{The Theorem of Borel-Weil and the Embedding Of Line Bundles}
Let $H\subset G$ be a closed complex subgroup and $\rho:H\to
\tmop{End}(V)$ be a holomorphic representation. The fiber product
$F:=G\times_H V$ is the quotient space of $G\times V$ by the equivalence
relation
\begin{equation}
  (g_1,v_1) \sim (g_2,v_2), \text{ if } g_1=g_2h^{-1},
  v_1=\rho(h)v_2\text{ for some }h\in H.
\end{equation}

The projection $p: F \to G/H, [(g,v)]\mapsto gH$ is holomorphic and it can
be shown by a direct calculation that $p: F \to G/H$ is a vector
bundle with typical fiber $V$. We define a $G$-action on $F$ by
\begin{equation}
  x.[(g,v)] := [(xg,v)]. 
\end{equation}
This action induces a representation of $G$ on the vector space of
holomorphic sections\footnote{Since we only deal with holomorphic
  sections, we write $\Gamma(G/H,F)$ instead of $\Gamma_{hol}(G / H,
  F)$ for the rest of this chapter.} $\Gamma(G/H,F)$. For our purpose
it is useful to give this representation in the context of $H$-invariant
functions. Therefore, we identify the sections of $F\to G/H$ with
the $H$-invariant functions $f:G\to V$, i.e.,
\begin{equation}
  f(gh^{-1})=\rho(h)f(g) \quad \forall\ h\in H,g\in G.
\end{equation}
The $G$-action on these functions is given by
\begin{equation}
  x.f(g):=f(x^{-1}g) \quad \forall\ g,x\in G.
\end{equation}
In our context $H$ will be a Borel subgroup of $G$. 

After this preparation, we can formulate a weak version of the
Borel-Weil Theorem. For a more complete version we refer to
{\cite{huck:actions}} and {\cite{akhiezer}} for a treatment from the
complex analytic point of view. An algebraic approach can be found in
{\cite{wallachgoodman}}.

\begin{theorem}{\dueto{Borel-Weil}}
  \label{thm:borel-weil}
  Let $\rho : G \rightarrow \tmop{End} (V)$ be an irreducible
  representation with highest weight $\lambda$ and $B_-$ the Borel
  subgroup of the negative roots. Then $B_-$ acts by multiplication on
  $V_{\lambda}$ with character $\chi : B_- \rightarrow
  \mathbb{C}^{\ast}$, where $\left. \left. d_e \chi
    \right|_{\mathfrak{t}} = 2 \pi i \lambda \right.$ and the
  representation on $\Gamma(G / B_-, G\times_{B_-}\mathbb{C})$ is
  isomorphic to $\rho$.
\end{theorem}

\begin{proof}
  Cf.~\cite{akhiezer} Chap.~4.3.
\end{proof}

We now follow the classical construction of embedding a $G$-line
bundle into the dual of the vector space of its sections.  For this,
set $L=G\times_{B_-}\mathbb{C}$ and fix a basis $s_0, \ldots, s_N$ of
$\Gamma (G / B_-, L)$ and the corresponding dual basis
$s_0^\ast,\ldots, s_N^\ast$.

Let $\mathcal{Z}$ be the zero section of $L$. In the view of $L = G
\times \mathbb{C}/\!\sim$, the zero section is exactly given by the
elements of the form $(g, 0)$ for $g \in G$. We claim that we obtain an
equivariant, holomorphic map of $L \backslash \mathcal{Z}$ into $\Gamma (G /
B_-, L)^{\ast}$ by the following construction. We think of the $s_i$'s
as $B_-$-equivariant functions $G \rightarrow \mathbb{C}$ and define
\begin{equation}
  \label{eq:def_of_varphi}
  \varphi : L \backslash \mathcal{Z} \rightarrow \Gamma (G / B_-, L)^{\ast},
  [(g, z)] \mapsto \frac{1}{z} \sum_{j = 0}^N   s_j (g) s_j^{\ast} .
\end{equation}
This is reasonable because $z$ is not 0, otherwise we would have $[(g,
0)] \in \mathcal{Z}$. Moreover, $\varphi$ is well-defined.  Indeed, if
we take another representative $(gb^{- 1}, \chi (b) z)$, we get
\begin{equation}
  \sum_{j = 0}^N \frac{s_j (gb^{- 1})}{\chi (b) z}   s_j^{\ast} =
  \sum_{j = 0}^N   \frac{\chi (b)}{\chi (b)} \frac{s_j (g)}{z}
  s_j^{\ast}
\end{equation}
because the $s_j$ are equivariant under $B_-$, i.e.\
\begin{equation}
  s_j ( g b )  = \chi(b)^{-1} s_j(g).
\end{equation}

Next, we have to show the equivariance of $\varphi$ with respect to
the left action of $G$ on $L$ and the dual representation on $\Gamma
(G / B_-, L)^{\ast}$. For this let $x^{- 1} .s_j = \sum_{i = 0}^N a_i
s_i $ for a fixed $x \in G$ and we calculate
\begin{eqnarray*}
  x. \varphi ([g, z]) (s_j) & = & \frac{1}{z} \left( x. \sum_{i = 0}^N s_i (g)
    s_i^{\ast} \right)   (s_j)\\
  & = & \frac{1}{z} \sum_{i = 0}^N s_i (g) s_i^{\ast} (x^{- 1} .s_j)\\
  & = & \frac{1}{z} \sum_{i = 0}^N a_i s_i (g)\\
  & = & \frac{1}{z} (x^{- 1} .s_j) (g)\\
  & = & \frac{1}{z} s_j (xg)\\
  & = & \frac{1}{z} \sum_{i = 0}^N s_i (xg) s_i^{\ast}   (s_j)\\
  & = & \varphi ([xg, z]) (s_j) .
\end{eqnarray*}

Now, we claim that a vector of maximal weight is in the image of
$\varphi$. For this, consider the mapping
\begin{equation}
  \label{eq:def_of_j}
  j : G / B_- \rightarrow \mathbb{P}(\Gamma (G /B_-, L)^\ast),
  x \mapsto [s_0 (x) : \ldots : s_N (x)]  
\end{equation}
where the coordinates on the right hand side are the $s_j^\ast$. It is
an equivariant, holomorphic map of $G / B_-$ into the projective space
of $\Gamma (G / B_-, L)^\ast$. Thus, the image is a closed orbit in
$\mathbb{P}(\Gamma (G /B_-, L)^\ast)$. But the orbit of the projection
of a maximal weight vector is the only such orbit
(cf.~\cite{huck:actions}). By comparison of (\ref{eq:def_of_varphi})
and (\ref{eq:def_of_j}) we obtain that a vector $v_{\tmop{max}}$ of
maximal weight is in the image of $\varphi$. Actually, every $c\cdot
v_{\tmop{max}},\,c\neq 0$, is in the image then. By equivariance, we
conclude that the whole $U_-$-orbit through every vector of maximal
weight is contained in the image of $\varphi$.

We state the following lemma.

\begin{lemma}
  \label{lem:norm-lemma--for-bundles}
  Let $\varphi : L \backslash \mathcal{Z} \rightarrow \Gamma (G / B_-,
  L)^{\ast}$ be the equivariant embedding described above. Then any
  $K$-invariant unitary structure on $\Gamma (G / B_-, L)^{\ast}$ induces a
  $K$-invariant hermitian bundle metric which is unique up to
  multiplication by a constant.
\end{lemma}

\begin{proof}
  First, we recall that the $K$-action on $G / B_-$ is transitive
  (cf.~{\cite{huck:actions}}), so every $K$-invariant bundle metric
  is the same up to a constant factor and we have completed the proof
  once we find the induced bundle metric is indeed $K$-invariant.

  For $[g, z_1], [g, z_2] \in L$ we define
  \begin{equation} h_g (z_1, z_2) = \left\{ \begin{array}{ll}
        \frac{1}{\langle \varphi ([g, z_1]), \varphi ([g, z_2])
          \rangle} &
        \text{if } z_1, z_2 \neq 0\\
        0 & \text{otherwise}.
      \end{array} \right. 
  \end{equation}
  By the relation
  \begin{equation}
    \frac{1}{\langle \varphi ([g, z_1]), \varphi ([g, z_2]) \rangle} =
    \overline{z_1} z_2 \frac{1}{\langle \sum_{i = 0}^n f_i (g) f_i^{\ast},
      \sum_{i = 0}^n f_i (g) f_i^{\ast} \rangle} = \overline{z_1} z_2
    \frac{1}{\| \varphi ([g, 1])\|^2} \label{buendelmetrik}
  \end{equation}
  we obtain a hermitian inner product at every point, since $\varphi$
  is well-defined and has only values different from zero. We claim
  that $h_g$ is \ a smooth bundle metric. We see that $h_g$ is
  continuous and smooth outside the zero section. Recall the standard
  fact that such a bundle metric is then smooth everywhere (cf.\
  \cite{lang:hyperbolic} p.96).  This metric is $K$-invariant because
  $\langle \cdot, \cdot \rangle$ is $K$-invariant and $\varphi$ is
  equivariant.
\end{proof}

\begin{lemma}\label{tensorbuendellemma}
  Let $L_1 \rightarrow G / B_-$ and $L_2 \rightarrow G / B_-$ be
  homogeneous complex line bundles that realize the representations
  corresponding to the highest weights $\lambda_1$ and $\lambda_2$. 

  The representation of highest weight $\lambda_1 + \lambda_2$ is then
  realized by $\Gamma(G/B_-,L_1\otimes L_2)$.
\end{lemma}

\begin{proof}
  This is a corollary to the Theorem of the Highest Weight as written
  in {\cite{huck:actions}} Chap.~7.1.
\end{proof}

\section{Symplectic geometry and momentum maps}
In this section the basic definitions of symplectic manifolds and 
momentums maps are given.

By definition a symplectic manifold $(M,\omega)$ is a real manifold $M$
with a non-degenerate two-form $\omega$. 

An action of a Lie group $H$ on $M$ is said to be symplectic if
\begin{equation}
  h^\ast \omega = \omega\ \forall\ h\in H.
\end{equation}
Before we define the notion of a momentum map, let us fix the notation.

The induced vector field of the flow $\exp(-\xi t)$ on $M$ is denoted
by $X_\xi$ and the Lie derivative along $X_\xi$ by
$\mathcal{L}_{X_\xi}$.  For a smooth map $\mu:M\to \tmop{Lie}(H)^\ast$ we
obtain an induced map $\mu^\xi:\tmop{Lie}(H)\to C^\infty(M)$ by
\begin{equation}
  \mu^\xi(x) := \mu(x)(\xi)\ \forall\ \xi\in\tmop{Lie}(H).
\end{equation}

\begin{definition}
  Let $(M,\omega)$ be a symplectic manifold on which $H$ acts by
  symplectic transformations.

  A \textbf{momentum map} is an equivariant, smooth map $\mu:M\to
  \tmop{Lie}(H)^*$ such that
  \begin{equation}
    d(\mu^\xi)= \omega(X_\xi,\cdot),
  \end{equation}
  where the action on $\tmop{Lie}(H)^\ast$ is the coadjoint action.
\end{definition}

We will use the momentum map only in the context of representations of
a compact Lie group. Let $\rho:K\to \tmop{U}(V)$ be a unitary
representation of the compact Lie group $K$ on a finite-dimension
vector space $V$. This representation induces an action of $K$ on
$\mathbb{P}(V)$ which is symplectic with respect to the Fubini-Study
metric on $\mathbb{C}$. Recall that the Fubini-Study metric is given
by the imaginary part of the form $\frac{i}{2}\partial\bar{\partial}
\log ||\cdot||^2$ pushed down from $V\backslash\{0\}$ to $\mathbb{P}(V)$.

\begin{theorem}
  Let $\rho:K\to \tmop{U}(V)$ be an irreducible representation of
  highest weight $\lambda$.

  The map $\mu:\mathbb{P}(V)\to\mathfrak{k}^*$ given by
  \begin{equation}
  \mu^\xi([v])= -2i\frac{\langle v,\rho_*(\xi).v\rangle}{\langle v,
    v\rangle}\forall\xi\in\mathfrak{k}, v\in\mathbb{P}(V)
  \end{equation}
  is the unique momentum map and
  \begin{equation}
    \mu([v_{max}])=\lambda
   \end{equation}
   for any vector $v_{max}$ of highest weight. 
\end{theorem}

\begin{proof}
  Cf.~\cite{huck:actions} Chap.~IV.7.   
\end{proof}

\section{Generalities on Level Spacings}
In this section we summarize the foundational facts on level spacings.

\subsection{The Nearest Neighbor Distribution}
\label{sec:near-neighb-distr}
The material in this subsection applies to arbitrary $N$-tuples of real
numbers, $N>1$. Later on it will be used only for eigenvalues of hermitian
matrices.

\begin{definition}
  Let $X=(x_1,\dots,x_N)\in\mathbb{R}^N$ be an $N$-tuple of real
  numbers, ordered by increasing value
  \begin{equation}
    x_1 \leq x_2 \leq \dots \leq x_N.
  \end{equation}
  The \textbf{nearest neighbor distribution} of $X$ is the Borel
  measure on $\mathbb{R}$ given by
  \begin{equation}
    \label{eq:def_nearestneighbordistribution}
    \mu(X)(A) = \int_A \frac{1}{N}\sum_{i=1}^{N-1} 
    \delta\left(y-\frac{N}{x_N-x_1}\cdot(x_{j+1}-x_j)\right) dy,
  \end{equation}%
  \nomenclature{$\mu(X)$}{The nearest neighbor statistics of the tuple $X$}%
  if $x_1\neq x_N$, and
  \begin{equation}
    \mu(X)(A) = \frac{N-1}{N}\int_A \delta(y) dy, 
  \end{equation}
  if $x_1=\ldots=x_N$, where $A$ is a Borel set in $\mathbb{R}$ and
  $\delta(y-p)$ denotes the Dirac measure with mass one at the point
  $p$.
\end{definition}
Thus, if $x_1\neq x_N$, $\mu(X)$ is a measure of total
mass\footnote{Note that for this reason it is common to use the factor
  $\frac{1}{N-1}$ in front of the sum and $N-1$ instead of $N$ inside
  the $\delta$ measures, but we will see that this is of no importance for
  questions of convergence.}  $1-\frac{1}{N}$ with expectation value
\begin{equation}
  \begin{aligned}
    E(\mu(X))&=\int_{\mathbb{R}}y\,d\mu(X)(y)\\
    &=\frac{1}{N}\sum_{j=1}^{N-1}\frac{N}{x_N-x_1}\cdot(x_{j+1}-x_j)\\
    &=\frac{1}{x_N-x_1}\cdot (x_N-x_1)=1.
  \end{aligned}
\end{equation}
\begin{remark}
  Note, that $\mu(X)$ does not change under scalar multiplication, i.e., 
  \begin{equation}
    \mu(aX)=\mu(X)\ \forall a\in\mathbb{R},a\neq 0
  \end{equation}
  nor under diagonal addition
  \begin{equation}
    \mu((x_1+a,\dots,x_n+a)) = \mu((x_1,\dots,x_n)).
  \end{equation}

\end{remark}
\indent If we know \emph{a priori} that our $N$-tuple $X$ is contained in 
 $[a,b]^N \tmop{mod} 1$, it is customary to measure the wrapped
around distance between $x_N$ and $x_1$:
\begin{equation}
  b-a-x_N+x_1
\end{equation}
and to replace $x_N-x_1$ in the denominator by $b-a$:
\begin{equation}
  \label{eq:def_nearestneighbordistribution_torus}
  \begin{aligned}
    \mu_w(X)(A) & = \frac{1}{N} \int_A  
    \delta\left(y-\frac{N}{b-a}(b-a-x_N+x_1)\right) \\
      & + \sum_{i=1}^{N-1} \delta\left(y-\frac{N}{b-a} \cdot(x_{j+1}-x_j)
\right) dy.
  \end{aligned}  
\end{equation}
Note that the total mass of this measure is $1$ and the expectation
value is also $1$.

Our main example for the above measure on the torus is given by
the logarithms of eigenvalues of a unitary matrix $U$. Here $a=0$ and
$b=2\pi$ and we obtain the following definition
\begin{equation}
  \label{eq:def_nearestneighbordistribution_circle}
  \mu_c(X)(A) = \frac{1}{N}\int_A
  \delta\left(y-\frac{N}{2\pi}(2\pi-x_N+x_1)\right)
  +\sum_{i=1}^{N-1} \delta\left(y-\frac{N}{2\pi}\cdot(x_{j+1}-x_j) \right) dy
\end{equation}%
\nomenclature{$\mu_c(X)$}{The nearest neighbor statistics of the angles
  $X$ with wrapping at $2\pi$}%
where $X=(x_1,\dots,x_N)$ is the set of ordered logarithms of the
eigenvalues with multiplicities, i.e. $\tmop{spec}(U) =
\{e^{ix_1},\dots,e^{ix_N}\}$. Here the differences $x_{j+1}-x_j$ are
the angles between the eigenvalues and $2\pi-x_N+x_1$ is the angle
between the first and the last eigenvalue.

In physical models such a wrapping occurs naturally
because the only physical data is encoded in the difference of the
arguments of the $e^{ix_j}$. Thus, the choice of the branch of the logarithm
is artificial, i.e.\ the position of zero cannot be measured.

\begin{note}
The measures $\mu_c(X_N)$ and $\mu_w(X_N)$ are no longer invariant under scalar 
multiplication.
\end{note}

\subsection{The Kolmogorov-Smirnov Distance}
\label{sec:kolom-smirn-dist}

Since we want to discuss convergence of measures on the real line,
we need a precise notion of the type of convergence we are dealing
with. For us only two types of convergence are important: the weak
convergence of distribution functions and the $\sup$-norm convergence
of distribution functions.

Recall that a sequence of measure $\mu_n$ is said to converge weakly
to a measure $\mu$ if for every bounded, continuous function $f$ the
following holds:
\begin{equation}
  \lim_{n\to\infty} \int fd\mu_n = \int f d\mu.
\end{equation}

\begin{definition}
  Let $\mu,\nu$ be Borel measures on $\mathbb{R}$ of finite
  mass. The \textbf{Kolmogorov-Smirnov distance} $d_{KS}$ of
  $\mu$ and $\nu$ is
  \begin{equation}
    \label{eq:def_dks}
    d_{KS}(\mu,\nu)=\sup_{t\in\mathbb{R}}\left|\int_{-\infty}^td\mu -
\int_{-\infty}^td\nu\right|,
  \end{equation}%
  \nomenclature{$d_{KS}$}{The Kolmogorov-Smirnov distance}%
  which is the $\sup$-norm for the difference of the
  cumulative distribution functions.

  We say a sequence of Borel measures $\mu_N$ converges to $\mu$ if
  $d_{KS}(\mu_N,\mu)$ converges to zero.
\end{definition}

\begin{remark}
  Convergence with respect to the Kolmogorov-Smirnov distance implies
  weak convergence. 
\end{remark}

\begin{proof}
  The convergence in the Kolmogorov-Smirnov distance implies the
  pointwise convergence of the cumulative distribution functions. But
  this implies weak convergence by a standard result of measure
  theory (cf.~\cite{elstrodt} chap. 8 Theorem 4.12 ).
\end{proof}

We now show that the scaling, with $N-1$ instead of $N$ which is
common in the literature (cf.~\cite{mehta}), gives the same results.

\begin{lemma}
  Let $(X_N)_{N\in\mathbb{N}}$ be a sequence of $N$-tuples such that
  $X_N\in\mathbb{R}^N$ and let $\nu$ be a Borel measure on $\mathbb{R}^+$
  with continuous density function $p(x)$ with respect to the Lesbesgue
  measure. Then the following are equivalent:
  \begin{enumerate}
  \item $\lim_{N\to\infty}\mu(X_N)=\nu$.
  \item $\lim_{N\to\infty}\mu_1(X_N)=\nu$, where
  \begin{equation}
    \mu_1(X)(A) = \frac{1}{N-1}\int_A \sum_{i=1}^{N-1}
    \delta\left(y -\frac{N}{x_N-x_1}\cdot(x_{j+1}-x_j)\right)dy.
  \end{equation}
  \item $\lim_{N\to\infty}\mu_2(X_N)=\nu$, where
  \begin{equation}
    \mu_2(X)(A) = \frac{1}{N-1}\int_A \sum_{i=1}^{N-1}
    \delta\left(y -\frac{N-1}{x_N-x_1}\cdot(x_{j+1}-x_j)\right)dy.
  \end{equation}
  \end{enumerate}
\end{lemma}

\begin{proof}
  The equivalence of 1.\ and 2.\ is clear, since
  $d_{KS}(\mu_1(X)(A),\mu_2(X)(A))=\frac{1}{N}$. For the proof of
  the equivalence of 2.\ and 3.\ we note that
  \begin{equation}
    \mu_1(X)([0,y])=\frac{1}{N-1}\tmop{card}\left\{j :
    \frac{x_{j+1}-x_j}{x_N-x_1}\cdot N\leq y \right\}
  \end{equation}
  and
  \begin{equation}
    \mu_2(X)([0,y])=\frac{1}{N-1}\tmop{card}\left\{j :
    \frac{x_{j+1}-x_j}{x_N-x_1}\cdot (N-1)\leq y \right\}.
  \end{equation}
  Therefore we see that
  \begin{equation}
    \mu_2(X)([0,y])=\mu_1(X)\left(\left[0,\frac{N-1}{N}\cdot y\right]\right).
  \end{equation}
  Now suppose 2.\ is true. Then
  \begin{equation}
    \begin{aligned}
    \left|\nu([0,y])-\mu_2(X)([0,y])\right|
    &\leq \left| \mu_1(X)\left(\left[0,\frac{N-1}{N}\cdot
          y\right]\right) - \nu([0,y]) \right|\\
    &\leq \left|\mu_1(X)\left(\left[0,\frac{N-1}{N}\cdot
          y\right]\right) -\nu\left(\left[0,\frac{N-1}{N}\cdot
          y\right]\right)\right| \\
    &+\left| \nu\left(\left[0, \frac{N-1} {N}\cdot y\right]\right) 
          - \nu([0,y])\right|.
    \end{aligned}
  \end{equation}
  Since $p(x)$ is continuous, the cumulative density function of $\nu$
  is uniformly continuous and the lemma follows from the estimation by
  a direct $\frac{\epsilon}{2}$ proof. Therefore, 2.\ implies 3.\ and,
  analogously, we see that the converse is true. 
\end{proof}

In the literature one often comes across histograms with densities
plotted into them for the nearest neighbor statistics. Compare Figure
\ref{fig:histogram} in the introduction, where we see a histogram
containing two curves.

Let us briefly discuss how the histogram in Figure \ref{fig:histogram}
was built. We start with an $N$-tuple $X=(x_1,\dots,x_N)$ of
non-decreasing real numbers as input and consider the $N-1$ rescaled
nearest neighbor distances
\begin{equation}
  \phi_j = \frac{N-1}{x_N-x_1}\cdot(x_{j+1}-x_j)\ \forall j=1,\dots,N-1.
\end{equation}
Now, we divide the real line into bins of some fixed width $w$ and
count the number of $\phi_j$ in each bin. At last we scale the height
of the boxes with a common factor such that the total area of the
histogram is one.

One usually has some measures with a continuous density function with
which to compare the histogram. In Figure \ref{fig:histogram}
two such densities are plotted.

This can be thought of as a visualization of the
$d_{\tmop{KS}}$-convergence in the following sense:  As the width $w$
becomes smaller and the $N$-tuples become larger, the histogram should
approach the density of the limit measure. This can be made precise in
the following way. Fix $p\geq 0$ and think of the histogram restricted
to $[0,p]$ as a Riemannian sum, which should converge to the integral
of the density over $[0,p]$.    

Unfortunately, this depends on the ratio of $N$
and $w$. Being a bit sloppy we can say that at the locus $\phi_j$ we
obtain a contribution of mass $1/(N-1)$ if the width is small
enough. This is exactly the point of the definition of $\mu(X)$.
Thus, a visualization of the convergence is obtained, although it is not
without problems because of the new dependence on the parameter
$w$. 

\subsection[Approximating $N$-tuples]{Approximating \boldmath$N$-tuples\unboldmath}
The following lemma shows how to construct approximating
$N$-tuples for any absolutely continuous measure.

\begin{lemma}
  \label{lem:approxpoisson}
  Let $\mu$ be a measure on $\mathbb{R}_{\geq 0}$ with continuous
  density $f$ such that
  \begin{equation}
    \int_0^{\infty}xf(x)dx \in [0,1].
  \end{equation}
  For every $N\geq 3$ there exists an $N$-tuple
  $X=(x_1,\ldots,x_N)$, $x_1\leq \ldots \leq x_n$ such that
  \begin{equation}
    d_{KS}(\mu(X),\mu) \leq \frac{2}{N-1}.
  \end{equation}
  Moreover, $x_1$ can be chosen to be $0$.
\end{lemma}

\begin{proof}
  First, define $y_j$ by the requirement
  \begin{equation}
    \frac{j}{N}=\int_0^{y_j}d\mu \quad \forall\,j=1,\ldots,N-1.
  \end{equation}  
  If we could choose $X$ in such a way that $\mu(X)$ has mass
  $\frac{1}{N}$ exactly at the $y_j$, i.e.,
  \begin{equation}
    \label{eq:bestapproxpoisson}
    y_j \stackrel{!}{=} \frac{N}{x_N-x_1}(x_{j+1}-x_j) \quad \forall\,j=1,\ldots,N-1,
  \end{equation}
  then
  \begin{equation}
     |\int_0^y d\mu-\int_0^y d\mu(X)| \leq \frac{1}{N-1} 
  \end{equation}
  since the cumulative distribution functions agree at the $y_j$ by
  construction and differ only by at most $\frac{1}{N-1}$ as indicated
  in the following picture for a certain measure.

  \begin{figure}[ht]
    \centering
    \includegraphics[width=12cm]{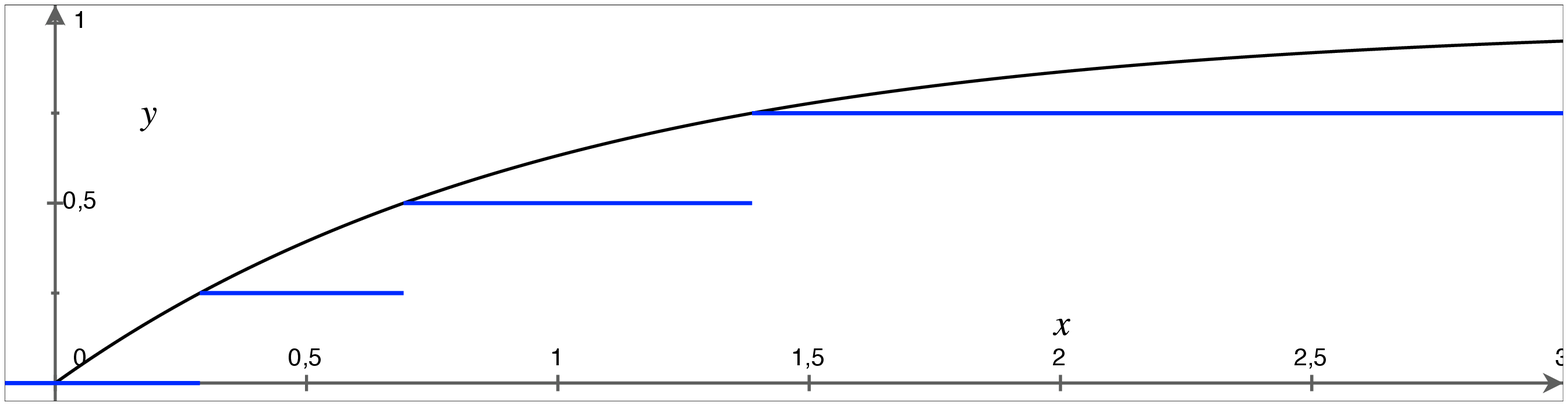}
    \caption{Approximation of $\mu_{\tmop{Poisson}}$.}
    \label{fig:approxpoisson}
  \end{figure}
  
  Unfortunately, the system (\ref{eq:bestapproxpoisson}) might have no
  solution since
  \begin{equation}
    \sum_{j=1}^{N-1} y_j \neq N = \sum_{j=1}^{N-1} 
    \frac{N}{x_N-x_1}(x_{j+1}-x_j). 
  \end{equation}
  Thus, we redefine $y_{N-1}$ in the following way
  \begin{equation}
    y_{N-1} = N - \sum_{j=1}^{N-2} y_j.
  \end{equation}
  We claim that $y_{N-1}$ is non-negative, i.e.,
  \begin{equation}
    N \geq \sum_{j=1}^{N-2}y_j.
  \end{equation}
  This follows at once from the inequality
  \begin{align}
    \int_0^{\infty} xf(x) dx&=\int_0^{y_1}xf(x)dx + \ldots +
    \int_{y_{N-2}}^{\infty}xf(x) dx \notag\\
    & \geq 0\cdot\int_0^{y_1}f(x)dx + y_1\cdot\int_{y_1}^{y_2}f(x)
    dx\ldots +
    y_{N-2}\int_{y_{N-2}}^{\infty}f(x) dx \\
    & =y_1\frac{1}{N} + \ldots y_{N-3}\frac{1}{N} + y_{N-2}\frac{2}{N}
    \geq \frac{1}{N} \sum_{j=1}^{N-2}y_j,\notag
  \end{align}
  because $\int_0^{\infty}xf(x)dx \in [0,1]$ by assumption.

  Writing (\ref{eq:bestapproxpoisson}) as a linear system
  \begin{equation}
    (x_N-x_1)y_j - N(x_{j+1}-x_j) = 0
  \end{equation}
  and calculating the space of solutions, we see that the solutions
  depend on real parameters $a$ and $b$:
  \begin{equation}
    x_j = a + \frac{b}{N} \sum_{1}^{j-1} y_j \quad \forall\, j=1,\ldots,N.
  \end{equation}
  Note that no solution with $b=0$ solves the original
  problem (\ref{eq:bestapproxpoisson}).

  For any solution $X$ with $b\neq 0$ the estimate  
  \begin{equation}
      d_{KS}(\mu(X),\mu) \leq \frac{2}{N-1}    
  \end{equation}
  holds since $y_{N-1}$ is not in the optimal position any
  more. Hence, we have to adjust by the factor $\frac{2}{N-1}$.
\end{proof}

\begin{corollary}
  \label{cor:approxpoisson}
  Let $\mu$ be an absolutely continuous measure on $\mathbb{R}_{\geq
    0}$ with $\int_0^\infty xd\mu\in[0,1]$ and $p>0$ be a fixed
  integer. Then for any $p$-tuple $(z_1,\ldots,z_p)$ and any $N\geq
  p+2$ there is an $N$-tuple $X=(x_1,\ldots,x_N)$ such that every
  $z_j$ occurs as one of the $x_k$ and
  \begin{equation}
    d_{KS}(\mu(X),\mu) \leq \frac{2+p}{N-1}.
  \end{equation}
\end{corollary}

\begin{proof}
  By Lemma \ref{lem:approxpoisson} we find an $N$-tuple $X$ such that
  \begin{equation}
    d_{KS}(\mu(X),\mu) \leq \frac{2}{N-1}.
  \end{equation}
  Due to the invariance of $\mu(X)$ under scalar multiplication and
  diagonal addition, we may assume that
  \begin{equation}
    x_2\leq z_j \leq x_N-1 \ \forall\ j=1,\ldots,p.  
  \end{equation}
  Now, insert the $z_j$ into the ordered sequence $x_1\leq \ldots \leq
  x_N$ at the corresponding positions and remove the closest $x_k$ for
  each $z_j$ inserted, as long as $x_k$ is neither $x_1$ nor $x_N$. In
  this case take the closest $x_k$ in the middle.  The resulting
  sequence is called $\tilde{X}$.

  In the picture of Figure \ref{fig:approxpoisson} we have changed $p$
  points of the jump loci of the approximating staircase function. Thus,
  we have to add an extra $\frac{p}{N-1}$ to the estimation. 
\end{proof}

\subsection[The Nearest Neighbor Statistics under $\tmop{exp}$]{The
  Nearest Neighbor Statistics under \boldmath$\tmop{exp}$\unboldmath}
\label{sec:near-neighb-stat-exp}
In this work the most important examples of sequences $(X_N)$ of
non-decreasing $N$-tuples are given by the spectra of sequences of
Hamiltonian operators on finite-di\-men\-sio\-nal Hilbert spaces or,
equally important, by the restrictions of Hamiltonian operators to
finite dimensional subspaces of some infinite-dimensional Hilbert
space such that the dimension of the finite-dimensional parts is
approaching infinity.

In the setting of general finite-dimensional Hilbert spaces the
Hamiltonians are just skew self-adjoint operators. The space of these
operators is again a finite-dimensional vector space. If $A$
is hermitian, the one-parameter group
\begin{equation}
  \{\tmop{exp}(iAt):t\in\mathbb{R}\}  
\end{equation}
is a subgroup of the unitary
group of this Hilbert space and the exponential mapping
$A\mapsto\tmop{exp}(iA)$ is surjective but not injective.

Now, we consider the spectrum of a unitary operator $\tmop{exp}(iA)$
and take the nearest neighbor statistics $\mu_c$ of the eigenangles,
i.e., the $a_j$ in the eigenvalue $e^{ia_j}$, where $0 \leq a_j<
2\pi$.

\begin{definition}
  \label{def:nearest-neighbor-unitary-matrix}
  Let $U \in U(N)$ be a unitary matrix, whose eigenvalues are given as
  $e^{2i\pi\phi_1},\dots,e^{2i\pi\phi_N}$, and
  $X(U)=(\phi_1,\dots,\phi_N)$.  The nearest neighbor statistics of
  the unitary matrix $U$ is $\mu_c (X(U))$.
\end{definition}
Frequently, we will write $\mu_A$ as an abbreviation for $\mu(X(A))$
and $\mu_U$ as abbreviation for $\mu_c (X(U))$.%
  \nomenclature{$\mu_U$}{The nearest neighbor statistics of the
    eigenangles of unitary matrix $U$}%
  \nomenclature{$\mu_A$}{The nearest neighbor statistics of the
    eigenvalues of hermitian matrix $A$}%

The nearest neighbor statistics of $\exp(iA)$ will not agree with the
nearest neighbor statistics of $A$ for two reasons, the first being
the wrapping discussed above and the second and more important is the
problem of reordering.

The eigenvalues $x_1\leq\dots\leq x_N$ of $A$ give the $\phi_j$ only
modulo $2\pi$, i.e.\
\begin{equation}
  a_j = \phi_j \tmop{ mod }2\pi 
\end{equation}
and it may happen that there are $j_1$ and $j_2$, such that $\phi_{j_1}<
\phi_{j_2}$ but $a_{j_1}>a_{j_2}$.

If, however, all eigenvalues are sufficiently close to each other,
meaning that they all lie in an interval of width $2\pi$, one does not
have to reorder, if choosing a different branch of the logarithm or
just by adding a constant to all eigenvalues such that the smallest
eigenvalue is 0.

\begin{lemma}
  Let $(X_N)_{N\in\mathbb{N}}$ be a sequence of non-decreasing
  $N$-tuples such that $X_N\in [0,2\pi[^N$ and let $\nu$ be a Borel
  measure on $\mathbb{R}^+$ with continuous density with respect to
  the Lebesgue measure. Assume that the difference between the largest
  and the smallest eigenvalue converges to $2\pi$. Then the following
  are equivalent:
\begin{enumerate}
 \item $\lim_{N\to\infty}\mu(X_N)=\nu$.
 \item $\lim_{N\to\infty}\mu_c(X_N)=\nu$.
\end{enumerate}
\end{lemma}

\begin{proof}
  Since the differences between the largest and the smallest
  eigenvalue converge to $2\pi$, the $\mu(X_N)$ come arbitrarily close
  to the $\mu_c(X_N)$ as is evident by their definitions.
\end{proof}

\begin{remark}
  \label{rm:exp_under_mu_counterexample}
  The above lemma is false if we drop the assumption on the largest
  and smallest eigenvalues.  Indeed, assume that every $X_N$ is
  contained in the subinterval $[0,\frac{1}{N}[$ with smallest
  eigenvalue 0 and largest eigenvalue $\frac{1}{N}$, then $\mu_C(X_N)$
  has the wrapping eigenangle given by
  \begin{equation}
    a_N = \frac{N}{2\pi}(2\pi-\frac{1}{N}) = N-\frac{1}{2\pi},
  \end{equation}
  all other eigenangles are less than or equal to $1/2\pi$. Therefore
  $\mu_c$ can only converge to a measure whose cumulative distribution
  function is 1 for all $t\in\mathbb{R}, t\geq 1/2\pi$.
\end{remark}

To summarize, care has to be taken if considering the nearest neighbor
statistics under $\tmop{exp}$. It is not enough to ensure that the
eigenvalues of a hermitian operator are in an interval $[0,2\pi]$ but
one must also ensure that the difference between the smallest and the
largest eigenvalue approaches $2\pi$.

\subsection{The Nearest Neighbor Statistics and the CUE Measure}
\label{sec:near-neighb-stat}

As discussed above we are mainly interested in the nearest neighbor
statistics associated to unitary matrices. We give some more details about these statistics here. In this section $\mu_c(X(A))$ will be abbreviated by $\mu_A$. 

The following lemma is necessary in certain of our applications.
\begin{lemma}
  \label{lem:dks-continous}
  If $\nu$ is an absolutely continuous probability measure on
  $\mathbb{R}$, then the map
  \begin{equation}
    U(N)\to [0,1],\ A\mapsto d_{KS}(\,\nu,\mu_A\,)    
  \end{equation}
  is continuous.
\end{lemma}
\begin{proof}
  Cf.\ \cite{katzsarnak} where the proof is given in lemma 1.0.11. and
  1.0.12.
\end{proof}

Since $U(N)$ is a compact group, functions on $U(N)$ can be
averaged. It is also possible to average the map $A\mapsto 
\mu_A$. This can be done in the following way. Let $\mu(U(N))$ 
denote the Borel measure given by
\begin{equation}
   \mu(U(N))(X) := \int_{U(N)} (\mu_A(X))\ d\tmop{Haar}(A) 
\end{equation}
for any Borel-measurable set $X$.

We now state Lemma 1.2.1 of \cite{katzsarnak}.
\begin{lemma}
  There exists an absolutely continuous probability measure $\nu$ on
  $\mathbb{R}$ with real analytic cumulative distribution function
  such that 
  \begin{equation}
    \mu(U(N))\to \nu \text{ weakly, as } n\to \infty.    
  \end{equation}
\end{lemma}

We call this measure \nomenclature{$\mu_{\tmop{CUE}}$}{The limit measure
  of the nearest neighborhood statistics of the unitary group}%
$\mu_{\tmop{CUE}}$. In \cite{katzsarnak} the following theorem is
given in a more general form as Lemma 1.2.6.
\begin{theorem}
  \label{lem:katzsarnak}
  For every $\epsilon>0$ there is a natural number $N_0$ such that
  \begin{equation}
    \label{eq:mean-dks-on-un}
    \int_{U(N)} d_{KS}(\mu_{\tmop{CUE}},\mu_A) d\tmop{Haar} \leq N^{\epsilon-1/6}
  \end{equation}
  for all $N\geq N_0$. 
\end{theorem}

The complete proof of the lemma and the theorem is given in all
detail in \cite{katzsarnak}, where it takes the first half of the
book, so it cannot be given here.

More details on $\mu_{\tmop{CUE}}$ can be found in \cite{mehta} and again in
\cite{katzsarnak}.

% This is the first theorem which states the limiting behavior of
% the nearest neighbor statistics of unitary matrices in a sequence of
% representations with increasing dimensions. Here the representations
% are the standard representations on $U(N)$ on $\mathbb{C}^N$. In
% chapter \ref{chap:poisson-spectral-stat} we prove a
% similar result.

% For this, let $T(N)$ be a maximal torus in $U(N)$ and $\mu(T(N))$ 
% denote the Borel measure given by
% \begin{equation}
%    \mu(T(N))(X) := \int_{T(N)} (\mu_A(X))\ d\tmop{Haar}(A), 
% \end{equation}
% for any Borel-measurable set $X$. Finally, let
% $\mu_{\tmop{Poisson}}$%
% \nomenclature{$\mu_{\tmop{Poisson}}$}{measure with density $\exp(-x)$ on positive real line}%
% denote the probability measure with density $\exp(-x)$ on the positive real line.

% \begin{lemma}
%   \label{lem:katzsarnak-tn}
%   The measure  $\mu_{\tmop{Poisson}}$ is the weak limit of averaged nearest
%   neighbor statistics over the Haar measure on $T(N)$, i.e. 
%   \begin{equation}
%     \mu(U(N))\to \nu \text{ weakly, as } n\to \infty.    
%   \end{equation}
%   Moreover, there is a natural number $N_0$, such that
%   \begin{equation}
%     \label{eq:mean-dks-on-tn}
%     \int_{U(N)} d_{KS}(\mu_{\tmop{Poisson}},\mu_A) d\tmop{Haar} \leq
%     e^{ \sqrt{ \log N }}
%   \end{equation}
%   for all $N\geq N_0$.
% \end{lemma}

%%% Local Variables:  
%%% mode: latex 
%%% TeX-master: "../thesis" 
%%% End: 

% LocalWords:  surjective

%\printglossary

\bibliography{Literatur}

\begin{thebibliography}{GHK00}

\bibitem[Akh91]{akhiezer}
Dmitri~N. Akhiezer.
\newblock {\em Lie Group Actions in Complex Analysis}.
\newblock Vieweg, 1991.

\bibitem[BtD85]{broecker:compactlie}
Theodor Br\"ocker and Tammo tom Dieck.
\newblock {\em Representations of Compact Lie Groups}.
\newblock Springer, New York, 1985.

\bibitem[CFS82]{cornfeldfominsinai}
I.~P. Cornfeld, S.V. Fomin, and Ya.~G. Sinai.
\newblock {\em Ergodic Theory}.
\newblock Springer, 1982.

\bibitem[Els04]{elstrodt}
J.~Elstrodt.
\newblock {\em {Ma\ss\ und Integrationstheorie}}.
\newblock {Springer}, forth edition, 2004.

\bibitem[Far01]{farenick}
Douglas~R. Farenick.
\newblock {\em Algebras of Linear Transformations}.
\newblock Springer, 2001.

\bibitem[GHK00]{haakekus}
Sven Gnutzmann, Fritz Haake, and Marek Kus.
\newblock Quantum chaos of {$SU_3$}-observables.
\newblock {\em J. Phys. A:Math. Gen.}, 33:143--161, 2000.

\bibitem[GK98]{gnutzmannkus}
Sven Gnutzmann and Marek Kus.
\newblock Coherent states and the classical limit on irreducible {$SU_3$}
  representations.
\newblock {\em J. Phys. A:Math. Gen.}, 31:9871--9896, 1998.

\bibitem[Gnu00]{gnutzmann}
Sven Gnutzmann.
\newblock {\em Klassischer Grenzfall, Semiklassik und Quantenchaos bei
  kollektiv gekoppelten n-Niveau-Atomen}.
\newblock PhD thesis, Universit\"at Essen, 2000.

\bibitem[Haa99]{haake}
Fritz Haake.
\newblock {\em Quantum Signatures of Chaos}.
\newblock Springer, second edition, 1999.

\bibitem[Huc91]{huck:actions}
Alan Huckleberry.
\newblock Introduction to group actions in symplectic and complex geometry.
\newblock In {\em Infinite Dimensional Lie Groups}, volume~51 of {\em DMV
  Seminarberichte}. Birkh\"auser, 1991.

\bibitem[Kna02]{knapp}
Anthony~W. Knapp.
\newblock {\em Lie Groups Beyond an Introduction}.
\newblock Birkh\"auser, second edition, 2002.

\bibitem[KS99]{katzsarnak}
Nicolas~M. Katz and Peter Sarnak.
\newblock {\em Random Matrices, Frobenius Eigenvalues and Monodromy}.
\newblock American Mathematical Society, 1999.

\bibitem[Lan87]{lang:hyperbolic}
Serge Lang.
\newblock {\em Introduction to Complex Hyperbolic Spaces}.
\newblock Springer, New York, 1987.

\bibitem[Meh91]{mehta}
Medan~L. Mehta.
\newblock {\em Random Matrices}.
\newblock Academic Press, 1991.

\bibitem[Per86]{perelomov}
Askold Perelomov.
\newblock {\em Generalized Coherent States and Their Applications}.
\newblock Springer, 1986.

\bibitem[Sin94]{sinai}
Ya.~G. Sinai.
\newblock {\em Topics in Ergodic Theory}.
\newblock Princeton University Press, 1994.

\bibitem[Sna00]{snaith}
Nina~Claire Snaith.
\newblock {\em Random Matrix Theory and Zeta Functions}.
\newblock PhD thesis, University of Bristol, 2000.

\bibitem[WG99]{wallachgoodman}
N.~Wallach and R.~Goodman.
\newblock {\em Representations and Invariants of the Classical Groups}.
\newblock Cambridge University Press, 1999.

\bibitem[Woo97]{woodhouse}
N.M.J. Woodhouse.
\newblock {\em {Geometric quantization}}.
\newblock {Oxford University Press}, second edition, 1997.

\end{thebibliography}
\bibliographystyle{alpha}
\end{document}